\documentclass[11pt,fleqn,leqno]{article}

\usepackage{fullpage}
\usepackage{amsmath,amsthm,amssymb,covington,paralist}
\usepackage[T1]{fontenc} 
\usepackage{titlesec}
\titleformat*{\section}{\large\bfseries}
\titleformat*{\subsection}{\large\bfseries}

\usepackage[authoryear]{natbib}
\setcitestyle{aysep={},yysep={,}}
\newtheorem{theorem}{Theorem}[section]

\def\phi{\varphi}        
\def\A{\forall}           
\def\E{\exists} 
\def\mand{\, \wedge \, }

\def\proves{\vdash}

\usepackage{hyperref}
\usepackage[usenames,dvipsnames]{color}
\hypersetup{colorlinks=true, urlcolor=blue,citecolor=Purple}

\title{On consistency and existence in mathematics}

\author{Walter Dean\\
\footnotesize{Department of Philosophy}\\
\footnotesize{University of Warwick}\\
\footnotesize{Coventry, CV4 7AL UK}\\
\footnotesize{\href{maitlto:W.H.Dean@warwick.ac.uk}{\texttt{W.H.Dean@warwick.ac.uk}}}\thanks{This is an extended version of a paper delivered to the Aristotelian Society on 15 June 2020 and which has been published in their Proceedings (Volume CXX, Part 3, pp. 349-393, \href{http://doi.org/10.1093/arisoc/aoaa017}{\texttt{DOI:10.1093/arisoc/aoaa017}}.).  Thanks are owed to Andrew Arana, Patricia Blanchette, Michael Detlefsen, Alberto Naibo, and Marco Panza for comments and discussion and to audiences in South Bend, Saint Andrews, Paris, and London where this material was presented.}}
\date{}

\begin{document}

\maketitle

\begin{abstract}
\noindent This paper engages the question \textsl{Does the consistency of a set of axioms entail the existence of a model in which they are satisfied?} within the frame of the Frege-Hilbert controversy.   The question is related historically to the formulation, proof, and reception of G\"odel's Completeness Theorem.   Tools from mathematical logic are then used to argue that there are precise senses in which Frege was correct to maintain that demonstrating consistency is \textsl{as difficult as it can be} but also in which Hilbert was correct to maintain that demonstrating existence given consistency is \textsl{as easy as it can be}.
\end{abstract}
\section{Introduction}

This paper aims to provide a contemporary reappraisal of a well-worn debate in philosophy of mathematics: does the consistency of a set of  axioms or other mathematical statements entail the existence of a model in which they are satisfied?  As most readers will be aware, this forms part of a complex of issues which has come to be known as the \textsl{Frege-Hilbert controversy}.  The controversy arose as a consequence of the publication in 1899 of  David Hilbert's \textsl{Grundlagen der Geometrie} \citeyearpar{Hilbert1899}.  This prompted Gottlob Frege to initiate an exchange of letters  wherein he and Hilbert debated a wide range of topics surrounding the axiomatic method, inclusive of the following:\footnote{Frege and Hilbert's letters are reprinted in translation in \citep{Frege1980a}.}
\begin{example}
\label{initques}
\begin{compactenum}[i)]
\item Should axioms be understood as basic truths about a domain or as implicit definitions of systems of arbitrary objects and relations in which they are satisfied?
\item What is the relation between the syntactic expression of axioms in a formal language and the propositions (or thoughts) which they express?
\item What is the correct method for demonstrating the consistency or logical independence of a set of axioms?
\end{compactenum}
\end{example}

Certain aspects of how Frege and Hilbert reacted to these questions can be attributed to the features of the logicist program which Frege had initiated in \citeyearpar{Frege1884} and the finitist consistency program which Hilbert would announce in \citeyearpar{Hilbert1904} and further develop in the 1910s-1930s.   More generally, however, Frege can be understood as an inheritor of the ``traditional'' view that geometry is the study of space while Hilbert was one of the original exponents of the ``modern'' view that mathematics in general -- and geometry in particular -- is concerned with the study of abstract structures and the statements which they satisfy.  The ultimate grounds of their disagreement were thus deep-seated.   But their exchange was also instrumental in bringing into focus a number of questions which are still with us today.  

Most readers will also be aware that the Frege-Hilbert controversy has itself inspired a substantial literature in contemporary philosophy of mathematics.\footnote{\citep{Resnik1974a}, \citep{Hallett1990,Hallett2010}, and \citep{Blanchette1996,Blanchette2012} provide paradigmatic waypoints for the present paper.}  Much of this has focused on how Frege's and Hilbert's views should be understood in the historical context of their exchange.  This preceded by at least 15 years the isolation of many of the features which we now take to be characteristic of mathematical logic -- e.g. a clear distinction between syntax and semantics, deductive and semantic consequence, or the order of variables and formulas.  One concern which pervades the contemporary discussion of the Frege-Hilbert controversy is thus that of making clear how Frege's and Hilbert's understanding of  consistency and existence (as well as attendant  notions such as logical consequence, validity, and independence) differed from currently accepted definitions.

There are indeed risks involved with separating questions like (\ref{initques}i-iii) from the historical milieu in which Frege and Hilbert debated them.   Nonetheless, the perspective which I will adopt in this paper is that there are also insights to be gained by viewing their controversy in light of subsequent work in mathematical logic.   One reason for adopting this approach is that the sequence of steps separating Hilbert's understanding of the notions just mentioned around the time of \textsl{Grundlagen der Geometrie} and that which he had adopted by the time of  \textsl{Grundz\"uge der theoretischen Logik} \citeyearpar{Hilbert1928} -- which is arguably the first ``modern'' textbook in mathematical logic -- is relatively short.   Another reason is to draw attention to a sequence of mathematical results obtained from the late 1920s onward -- some well-known, some less so -- which were not only inspired by Frege and Hilbert's exchange but can also be used to sharpen our understanding of the grounds of their disagreement and its contemporary significance (or at least so I will suggest).   


Within this frame I will argue for three basic theses:
\begin{example}
\label{points}
\begin{compactenum}[i)]
\item The Frege-Hilbert controversy informed G\"odel's formulation and proof of the Completeness Theorem for first-order logic as well as shaping its reception.
\item There is a precise sense in which Frege was correct to maintain that the general problem of demonstrating the consistency of a set of axioms is \textsl{as difficult as it can be}.
\item There is a precise sense in which Hilbert was correct to maintain that the general problem of demonstrating the existence of a structure satisfying a set of axioms is \textsl{as easy as it can be} conditional on their proof-theoretic consistency.
\end{compactenum}
\end{example}

The first of these claims is historical and will be developed in \S 2 and \S 3 alongside an overview of Hilbert's project in  \textsl{Grundlagen der Geometrie} and the developments in metamathematics which grew out of it.  Therein I will highlight in particular the role of Paul Bernays in locating the concerns which animated Hilbert's exchange with Frege within the context of what we now call \textsl{Hilbert's program}.  This includes the task of distilling from Hilbert's prior work a general technique for investigating the consistency of arbitrary axiom systems which Bernays referred to as \textsl{the method of arithmetization}.   The isolation of this method anticipates subsequent developments in computability theory -- in particular the proof of the unsolvability of the \textsl{Entscheidungsproblem} by Church and Turing, the definition of a \textsl{complete problem} by Post, and the introduction of the \textsl{arithmetical hierarchy} for classifying problems by Kleene -- which I will suggest can be naturally brought to bear on Frege's remarks about consistency in \S4 and Hilbert's remarks about existence in \S5. Finally in \S 6 and \S 7, I will apply these observations to explore some vestiges of the Frege-Hilbert controversy which linger in contemporary philosophy of mathematics -- e.g. in regard to the relationship between Hilbert's views and structuralism as well as the \textsl{gap} which I will suggest remains between consistency and existence in light of the metamathematical interaction of the Completeness Theorem with both arithmetic and geometry.   

\section{From the \textsl{Grundlagen der Geometrie} to the Completeness Theorem}

The first edition of Hilbert's \textsl{Grundlagen der Geometrie} was based on lecture courses delivered in G\"ottingen between 1894 and 1898 and was originally published in 1899 as part of a \textsl{Festschrift} for Gauss and Weber.  To fix some of its contributions which will be relevant here, Hilbert's own introduction can be quoted directly:\footnote{For more on the context and influence of  \textsl{Grundlagen der Geometrie} see the introduction to \citep{Hilbert1899} in \citep{Hilbert2004} and for mathematical reconstructions of the results mentioned below see \citep{Hartshorne2000}.  Page references and numbering below refer to the English translation of the 10th German edition \citep{Hilbert1971}.  Although substantial changes were made between editions, those relevant here will be flagged individually.}
\begin{quote}
{\footnotesize
Geometry, like arithmetic, requires only a few and simple principles for its logical development. These principles are called the axioms of geometry $\ldots \P \ldots$ This present investigation is a new attempt to establish for geometry a \textsl{complete}, and as \textsl{simple as possible}, set of axioms and to deduce from them the most important geometric theorems in such a way that the meaning of the various groups of axioms, as well as the significance of the conclusions that can be drawn from the individual axioms, come to light.  \\ \hspace*{1ex} \hfill \citep[p. 2]{Hilbert1971}}
\end{quote}  

Although Hilbert's axiomatization is in many ways similar to that of Euclid, his choice of primitive notions and his organization of theorems differ somewhat from that of the \textsl{Elements}.   The language in which he formalizes plane geometry includes the sorts \textsl{point} and \textsl{line} and the binary relation of \textsl{incidence} (relating a point and a line), the ternary relation of \textsl{betweenness} (relating triples of points), the 4-ary relation of \textsl{congruence of line segments} (relating two pairs of points), and the 6-ary relation of \textsl{congruence of angles} (relating two triples of points).  Hilbert's division of axioms is designed to facilitate not only the reconstruction of some of Euclid's theorems but also  exploration of the consistency and independence of various subsets of principles.   This gives rise to his organization of axioms into the categories \textsl{Incidence} (I), \textsl{Order} (II), \textsl{Congruence} (III), \textsl{Parallels} (IV), and \textsl{Continuity} (V).   The first three of these respectively pertain to the properties of the incidence,  betweenness, and the two congruence relations.  Hilbert took Playfair's form  [PP] of Euclid's Fifth Postulate as his parallel axiom.  The choice of continuity axioms was subject to the most revision across the various editions of the \textsl{Grundlagen}.   For present purposes, however, I will take this category to include only the so-called Circle-Circle Intersection Principle [CCP].\footnote{CCP states `If from distinct points $A$ and $B$, circles with radius $AC$ and $BD$ are drawn such that one circle contains points both in the interior and in the exterior of the other, then they intersect in two points, on opposite sides of $AB$'.   The explicit inclusion of CCP rather than a stronger principle such as Hilbert's Axiom of Completeness or the continuity axioms of Dedekind or Cantor (which appear in later editions of the \textsl{Grundlagen}) among Hilbert's axioms is an anachronism (as is the omission of the Archimedean axiom) .  But it is useful for understanding the structure of Hilbert's independence proofs and (per, e.g., \citealp{Baldwin2018}) it is also arguably more in keeping with the aims of the original text.}    

Although Hilbert originally expressed his axioms in natural language, the first-order formalization of these and related axioms was in fact a frequent source of examples in his later work (as I will discuss below).  And in fact it is straightforward to see how his axioms can be formalized in a two-sorted first-order language $\mathcal{L}_H$ with sorts $A,B,C\ldots$ for points, $\ell,\ell',\ell'',\ldots$ for lines and predicates $\mathit{In}(A,\ell)$, $\mathit{Btw}(A,B,C)$, $\mathit{Cong}_1(A,B,C,D)$, and $\mathit{Cong}_2(A,B,C,D,E,F)$ for incidence, betweenness, and congruence.  I will borrow from  \citet{Baldwin2018} the name $\mathsf{HP}$ (for `Hilbert planes') for the theory containing the axioms in groups (I), (II), and (III) formulated in $\mathcal{L}_H$.

In the first chapter of \citeyearpar{Hilbert1899}, Hilbert illustrates that these axioms are sufficient for developing  plane geometry by proving theorems from Euclid such as side-angle-side congruence for triangles.  He then introduces one of the main goals of the \textsl{Grundlagen} at the beginning of chapter II:
\begin{quote}
{\footnotesize
The axioms formulated in the five groups in Chapter I are not contradictory to each other, i.e., it is impossible to deduce from them by logical inference a result that contradicts one of them. In order to realize this a set of objects will be constructed from the real numbers in which all axioms of the five groups are satisfied. \hfill \citep[p. 29]{Hilbert1971}}
\end{quote}
Hilbert then describes how this can be accomplished by constructing a set of real numbers $\mathbb{P}$ (denoted $\Omega$ in the original text) which he describes as follows 
\begin{quote}
{\footnotesize Consider the field $[\mathbb{P}]$ arising from the number $1$ and the application of a finite number of times of the four arithmetic operations of addition, subtraction, multiplication, division and the fifth operation  $|\sqrt{1+x^2}|$, where $x$ denotes a number that results from these five operations.  \hfill \citeyearpar[p. 29]{Hilbert1971}}
\end{quote}
The set $\mathbb{P}$ is an example of what is now called an \textsl{ordered Pythagorean field} -- i.e. an ordered field which is also closed under the operation $\sqrt{x^2 + y^2}$.   In fact it is not difficult to see that $\mathbb{P}$ is the \textsl{minimal} such subfield of the real numbers $\mathbb{R}$ containing the rational numbers $\mathbb{Q}$ with this property.  A similar structure is the minimal \textsl{Euclidean field} $\mathbb{E}$ formed by closing $\mathbb{Q}$ under the operations of addition, subtraction, multiplication, division and arbitrary square roots of non-negative elements. 

Hilbert next described how the domains of points and lines as well the other predicates of $\mathcal{L}_H$ can be interpreted over $\mathbb{P}$.  These descriptions collectively determine a structure of the form $\mathcal{P} = \langle P^{\mathcal{P}}, L^{\mathcal{P}},\mathit{In}^{\mathcal{P}},\mathit{Btw}^{\mathcal{P}},\mathit{Cong}_1^{\mathcal{P}},\mathit{Cong}_2^{\mathcal{P}} \rangle$ via definitions such as the following:
\begin{example}
\label{infdefns}
\begin{compactenum}[i)]
\item $P^{\mathcal{P}} =$ the set of ordered pairs $A = (x,y) \in \mathbb{P}^2$
\item $L^{\mathcal{P}} =$ the set of triples $\ell = (u,v,w) \in \mathbb{P}^3$ such that $u$ and $v$ are not both $= 0$ 
\item $\mathit{In}^{\mathcal{P}} =$ the set of pairs of points $A = (x,y)$ and lines $\ell = (u,v,w)$ such that $u x + v y + w = 0$
\item $\mathit{Btw}^{\mathcal{P}} =$ the set of triples of distinct points $A = (x_1,y_1), B= (x_2,y_2),C = (x_3,y_3) $ incident on the same line $\ell = (u,v,w)$ and such that $x_1 < x_2 < x_3$ or $x_1 > x_2 > x_3$ if $u \neq 0$ and $v \neq 0$ (i.e. the line is not vertical) or if $y_1 < y _2 < y_3$ or $y_1 > y_2 > y_3$ otherwise
\end{compactenum}
\end{example}

As Hilbert gave these definitions himself, it is not a major anachronism to understand him as describing what would now be referred to as an \textsl{analytical model} of $\mathcal{L}_H$ -- i.e. a first-order structure whose sorts and relational symbols are defined over the real numbers $\mathbb{R}$.    But it would be a bigger stretch to understand him as having precisely anticipated at this stage the inductive definition of truth in such a structure.  What Hilbert did do, however, is to prove informally that various sets of the axioms he considers are satisfied in $\mathcal{P}$ and similar models.

The title of the second chapter of the \textsl{Grundlagen} is `The consistency and the mutual independence of the axioms'.   There can be little doubt that Hilbert understood  consistency deductively at this time -- i.e. a set of axioms is defined to be consistent if `it is impossible to deduce  from them by logical inference a result that contradicts one of them'.   But as there is no indication of a formal proof system in the \textsl{Grundlagen}, there is no apparent way in which he could have expressed the \textsl{soundness} of `logical inference' as a mathematical theorem at this time.   It is clear, however, that Hilbert also understood that were any of the sets of geometrical axioms he considers deductively inconsistent, an inconsistency would be transferred by his model constructions into the analytical domain -- i.e. 
\begin{quote}
{\footnotesize
Every contradiction in the consequences of the line and plane axioms $\ldots$ would therefore have to be detectable in the arithmetic of the field $[\mathbb{P}$]. \hfill \citeyearpar[p. 30]{Hilbert1971}}
\end{quote}

Hilbert's constructions are thus prototypical of what we would now call \textsl{model theoretic consistency proofs}. In fact if we adopt the familiar notation $\mathcal{M} \models \Gamma$ to express that the sentences $\Gamma$\footnote{Unless otherwise noted, the symbols $\Gamma, \Delta, \ldots$ should be understood in the sequel as ranging over \textsl{finite} sets of first-order sentences.   Since if $\Gamma = \{\gamma_1, \ldots, \gamma_n\}$ then $\Gamma$ is equivalent to a single sentence $\wedge \Gamma = \gamma_1 \mand \ldots \mand \gamma_n$, results which were originally formulated for single sentences thus also apply trivially to finitely axiomatizable theories.  This  covers almost all of the geometric systems originally considered by Hilbert as well as systems such as G\"odel-Bernays set theory [$\mathsf{GB}$] which are (more than) sufficient to formalize the fragments of analysis involved in Hilbert's model constructions.  But as will be evident to contemporary readers, many of the results discussed below also hold \textsl{mutatis mutandis} for computably axiomatizable theories such as $\mathsf{ZF}$ as well.}  are satisfied in the structure $\mathcal{M}$, then all of the consistency proofs in the \textsl{Grundlagen} can be understood as conforming to the following \textsl{model-theoretic consistency test}:
\begin{example}[(MCT)] 
If $\Gamma$ is a set of axioms over a language $\mathcal{L}$, then a sufficient condition for the \textsl{consistency} of $\Gamma$ is that there exist an $\mathcal{L}$-model $\mathcal{M}$ such $\mathcal{M} \models \Gamma$.   
\end{example} 

It is also clear that Hilbert understood the concept of independence deductively -- i.e. as corresponding to the fact `that no essential part of anyone of these groups of axioms can be deduced from the others by logical inference.' \citeyearpar[p. 32]{Hilbert1971}.   His independence proofs thus took the familiar form of showing that there exist models in which certain axioms while others are not.  Hilbert's proofs thus conformed to the following \textsl{model-theoretic independence test}:
\begin{example}[(MIT)]
If $\Gamma$ is a set of $\mathcal{L}$-sentences and $\phi$ is an $\mathcal{L}$-sentence, then a sufficient condition for the \textsl{independence of $\phi$ from $\Gamma$} is that there exist an $\mathcal{L}$-model $\mathcal{M}$ such that $\mathcal{M} \models \Gamma$ but $\mathcal{M} \not\models \phi$.
\end{example} 

If we assume that negation is in our language, then the problem of checking the independence of $\phi$ from $\Gamma$ via (MIT) can clearly be reduced to checking the consistency of $\Gamma \cup \{\neg \phi\}$ via (MCT).   Some of the specific results which Hilbert obtained in \citeyearpar{Hilbert1899} can be presented in this way as follows:
\begin{example}
\label{mcr}
\begin{compactenum}[i)]
\item CCP is independent of $\mathsf{HP} + \mathrm{PP}$.   This is shown in \S37 wherein Hilbert observes that the model $\mathcal{P}$ described in (\ref{infdefns}) based on the minimal Pythagorean field $\mathbb{P}$ is such that $\mathcal{P} \models \mathsf{HP} + \mathrm{PP} + \neg \mathrm{CCP}$.  On the other hand, he also shows that the model $\mathcal{E}$ which is constructed in the same manner as $\mathcal{P}$ starting from the minimal Euclidean field $\mathbb{E}$ satisfies $\mathsf{HP} + \mathrm{PP} + \mathrm{CCP}$.\footnote{These results are implicit in Hilbert's formulation of a general criterion for the possibility of ruler and compass constructions such as that required to demonstrate Theorem I.1 of Euclid's \textsl{Elements} -- i.e. `An equilateral triangle can be constructed on any given segment'.     In \S 37 this is shown to fail in structures such as $\mathcal{P}$ which do not satisfy CCP.   See \citep[\S 11]{Hartshorne2000}, \citep[pp. 204-206]{Hilbert2004}, \citep[\S 8.4.3]{Hallett2008} for discussion of Hilbert's engagement with various weakenings of the continuity axioms he considers which (like CCP) are sufficient to derive this statement.\label{CCPnote}}
\item PP is independent of $\mathsf{HP} + \mathrm{CCP}$.   This is demonstrated in \S 10 wherein it is shown that within the model $\mathcal{E}$ it is possible to define another structure $\mathcal{O} \models \mathsf{HP} + \mathrm{CCP} + \neg \mathrm{PP}$.\footnote{This construction is described only briefly in the text.   However in the lectures \citep[pp. 347-359]{Hilbert2004} on which \citeyearpar{Hilbert1899} is based, Hilbert had presented in more detail the so-called  \textsl{Poincar\'e disk model} $\mathcal{O}$ of hyperbolic geometry which has as its points $P^{\mathcal{O}}$ those of $\mathbb{E}$ inside a fixed circle $C$ centered at point $O$ and as its lines $L^{\mathcal{O}}$ those sets of points in $P^{\mathcal{O}}$ lying on circles which are orthogonal to $C$ or which lie on diameters of $C$.  See also \citep[\S 39]{Hartshorne2000}.}   
\end{compactenum}
\end{example}

The foregoing by no means exhaust the mathematical contributions of the \textsl{Grundlagen}.  But it is precisely these results and their demonstrations which originally caught Frege's eye in his correspondence with Hilbert. It is difficult to summarize their exchange in a manner which does justice to the broad range of issues in the background.   But the following passages from Frege's first letter [IV/3] to Hilbert about geometry introduce some of the issues which will be of concern below:\footnote{Both the labeling of letters and the page numbering below refer to the reprinting of Frege and Hilbert's correspondence in \citep{Frege1980a}.   I have also adapted the page and section of Frege's references to \textsl{Grundlagen der Geometrie} to match those in \citep{Hilbert1971}.}

\begin{quote}
\footnotesize{I was interested to get to know your \textsl{Festschrift} on the foundations of geometry, the more so as I myself had earlier been concerned with them, but without publishing anything. As was to be expected, there are many points of contact between my earlier unrealized attempts and your account, but also many divergencies. In particular, I thought I could make do with fewer primitive terms.  \hfill (p. 34)

The explanations of [\S 1, p. 3 and \S 3, pp. 5-6] are apparently of a very different kind, for here the meanings of the words `point', `line', `between' are not given, but are assumed to be known in advance.  At least it seems so.  But it is also left unclear what you call a point. One first thinks of points in the sense of Euclidean geometry, a thought reinforced by the proposition that the axioms express fundamental facts of our intuition.  But afterwards [\S 9, p. 29] you think of a pair of numbers as a point. \hspace*{1ex} \hfill (p. 35)

I call axioms propositions that are true but are not proved because our knowledge of them flows from a source very different from the logical source, a source which might be called spatial intuition. From the truth of the axioms it follows that they do not contradict one another. There is therefore no need for a further proof. \hfill (p. 37)} 
\end{quote}

The first two of these passages respectively introduce the fact that Frege had previously considered axiomatizing geometry on the basis of different non-logical terms than those employed by Hilbert and also that he assigned a more significant role to intuition in characterizing geometric axioms.  I will return to these aspects of the controversy below.  But it is the third passage which introduces the theme which will be initially relevant here -- i.e. since for Frege axioms have to correspond to \textsl{true} propositions about a given domain, once a given set of propositions are accepted as axioms, it follows without the need for additional proof that they are consistent.   

Hilbert [IV/4] replied to Frege as follows about this point:
\begin{quote}
\footnotesize{You write: `I call axioms propositions $\ldots$ From the truth of the axioms it follows that they do not contradict one another.' I found it very interesting to read this very sentence in your letter, for as long as I have been thinking, writing and lecturing on these things, I have been saying the exact reverse: if the arbitrarily given axioms do not contradict one another with all their consequences, then they are true and the things defined by the axioms exist. This is for me the criterion of truth and existence. \\ \hspace*{1ex} \hfill (p. 39)}
\end{quote}
This is among the most commonly quoted passages in the correspondence and is also largely responsible for the slogan `consistency implies existence' which has come to be associated with Hilbert.   But although such a view is often understood as a cornerstone of Hilbert's foundational standpoint, the early date of the correspondence should also be kept in mind.  For as I will discuss further in \S 3, his views on mathematical existence appear to have evolved substantially in light of the finitist consistency program he would later announce.\footnote{Another complication in basing a case for Hilbert's acceptance of the slogan on this passage is that he suggests here that the consistency of axioms entails not just the existence of a model in which they are satisfied but also their \textsl{truth}.   On the other hand in a contemporaneous formulation of a similar point Hilbert fails to mention truth in conjunction with existence :  `[I]f it can be proved that the attributes assigned to the concept can never lead to a contradiction by the application of a finite number of logical processes, I say that the mathematical existence of the concept $\ldots$ is thereby proved' \citeyearpar[p. 10]{Hilbert1900}.  Some evidence that Hilbert had in mind something closer to the contemporary notion of `truth in a structure' (rather than `truth simpliciter') is also provided by how he continued in his reply to Frege:  `The proposition ``Every equation has a root'' is true, and the existence of a root is proven, as soon as the axiom ``Every equation has a root'' can be added to the other arithmetical axioms, without raising the possibility of contradiction, no matter what conclusions are drawn. This conception is indeed the key to an understanding of $\ldots$ my \textsl{Festschrift} $\ldots$' (IV/4, p.40)   (Note that while the statement in question is false in the real field $\mathbb{R}$, its consistency is attested to by the fact it becomes true in the complex field $\mathbb{C}$.)}

This point will bear directly on the case I will give for theses (\ref{points}ii,iii) below.  It is also related to another of Frege's concerns which emerges in his second letter to Hilbert: 

{\footnotesize 
\begin{quote}
[W]e must ask, What means have we of demonstrating that certain $\ldots$ requirements $\ldots$ do not contradict one another? The only means I know is this: to point to an object that has all those properties, to give a case where all those requirements are satisfied. It does not seem possible to demonstrate the lack of contradiction in any other way. \hfill (p. 43)
\end{quote}
}
\noindent Frege is here using `requirement'  in something like the postulational sense in which he takes Hilbert to be using the term `axiom'.  What he is saying in this passage  (and others to be considered below) can thus be understood as suggesting that the existence of a satisfying structure is in fact a \textsl{necessary} condition to show the consistency of a set of propositions.  It hence seems reasonable to understand him as holding that any suitable criterion for consistency must impose a condition \textsl{at least as stringent} as (MCT) and similarly for independence with respect to (IT).

Hilbert's consistency and independence proofs in \citeyearpar{Hilbert1899} satisfy these criteria in an apparently paradigmatic way.   The fact that Frege was still dissatisfied with them attests to the fact that he wished to impose a yet stricter criteria.   The task of precisely articulating the additional conditions which Frege wished to enforce is delicate.  However it is generally agreed that Frege's reservations derived in part from his understanding of axioms not as \textsl{sentences} (i.e. linguistic items) but rather what he called \textsl{thoughts} (i.e. non-linguistic entities closer to propositions).  He also held that thoughts have a determinate subject matter and thus do not to admit to reinterpretation of the terms which we employ to express them.  Thus although Frege can be understood as acknowledging that results like (\ref{mcr}i-iii) present models satisfying various axiom systems in Hilbert's sense (i.e. sets of sentences), he still maintained that these models are not of the right sort to demonstrate consistency of the systems in his own sense (i.e. sets of thoughts).  In particular, Hilbert's models are not composed from geometrical points and lines which Frege took to be expressed by the terms `point' and `line'  but rather from pairs and triples of real numbers over fields like $\mathbb{P}$ or $\mathbb{E}$.\footnote{For more on these points see \citep{Blanchette1996,Blanchette2012}.}

I will return to this (familiar) component of the ideological background which separated Frege from Hilbert in \S 4.  But it is equally notable that already at the time of his exchange with Frege, Hilbert had himself become dissatisfied with (MCT) and (IT) as ultimate criteria for consistency and independence.   This concern comes almost to the surface in one of Frege's subsequent letters [IV/7] in which he thanks Hilbert for sending him a copy of his address `Mathematische Probleme' \citeyearpar{Hilbert1900} and then presses him for further details on another technique which Hilbert alludes to therein.\footnote{`It seems to me that you believe yourself to be in possession of a principle for proving lack of contradiction which is essentially different from the one I formulated in my last letter and which, if I remember right, is the only one you apply in your [\citeyear{Hilbert1899}].  If you were right in this, it could be of immense importance, though I do not believe in it as yet, but suspect that such a principle can be reduced to the one I formulated and that it cannot therefore have a wider scope than mine. It would help to clear up matters if in your reply to my last letter $\ldots$ you could formulate such a principle precisely and perhaps elucidate its application by an example.' [IV/7], pp. 49-50 \label{fnote1}}  The passage about which Frege was inquiring is most likely the following:

{\footnotesize
\begin{quote}
In geometry, the proof of the compatibility of the axioms can be effected by constructing a suitable field of numbers, such that analogous relations between the numbers of this field correspond to the geometrical axioms. Any contradiction in the deductions from the geometrical axioms must thereupon be recognizable in the arithmetic of this field of numbers. In this way the desired proof for the compatibility of the geometrical axioms is made to depend upon the theorem of the compatibility of the arithmetical axioms. 

On the other hand a direct method is needed for the proof of the compatibility of the arithmetical axioms. The axioms of arithmetic are essentially nothing else than the known rules of calculation, with the addition of the axiom of continuity $\ldots$ I am convinced that it must be possible to find a direct proof for the compatibility of the arithmetical axioms, by means of a careful study and suitable modification of the known methods of reasoning in the theory of irrational numbers. \hfill \citeyearpar[p. 1104]{Hilbert1900}
\end{quote}
}

Hilbert here takes five steps in quick succession: i) he summarizes the technique of his geometric consistency proofs in a manner which is  reminiscent of the modern notion of \textsl{model theoretic interpretability} (see \S 3);  ii) he observes that this method provides a means of converting any potential proof of a contradiction from the geometric axioms into a contradictory proposition about the relevant numerical fields -- i.e that the model constructions described above provide only \textsl{relative consistency proofs} (as we would now put it); iii) he alludes to the fact that he has recently introduced an axiomatic theory for reasoning about real numbers relative to which the prior step can be made precise;  iv) he states that unlike in the case of geometry, the consistency of this other axiomatization must be proven `directly' rather than via interpretation; v) he indicates a direction by which an appropriately `direct' proof of the consistency of analysis might be found.  

As I will discuss further in \S 3, it is fair to say that points i) and ii) reflect the consensus view of the contribution of \textsl{Grundlagen der Geometrie}.  With respect to point iii), it appears likely that Hilbert had also sent Frege a copy of his paper `\"Uber den Zahlbegriff' \citeyearpar{Hilbert1900a} in which an axiom system similar to what we would now call \textsl{analysis} (or \textsl{second-order arithmetic}) is indeed presented.  The situation with respect to iv) and v) is more complicated.  For Hilbert failed at this time to provide either a characterization of what he meant by a `direct consistency proof' as distinct from (MCT) or anything more than a promissory note that such a proof would be forthcoming.  It was thus quite reasonable of Frege to press Hilbert for further details.  

Despite the confidence he projected, the prior passage forms part of the announcement of what later became known as \textsl{Hilbert's second problem} (cf., e.g., \citealp{Kreisel1976}).   While the exchange just recorded effectively marked the end of Hilbert's correspondence with Frege, we now know that it would take 25 years of work before it would be possible to provide a satisfactory characterization of what form such a proof might take.   On the other hand, Hilbert's announcement of this goal also anticipated the development of what he would come to call \textsl{metamathematics}  within a stream of developments which would come to include G\"odel's proof of the Completeness Theorem in 1929.   As many of the intervening events are well-known (and aptly described in \citealp[I.3]{Sieg2013} and in the introduction to \citealp{Hilbert2013}), it will suffice to summarize a few points from the intervening years.  

With the exception of his \citeyearpar{Hilbert1904} address `\"Uber die Grundlagen der Logik und der Arithmetik' -- which anticipates several aspects of his finitist consistency program  -- Hilbert worked mostly on mathematical physics until the late 1910s, at which time he began offering lecture courses on the foundations of mathematics.   The intervening years had seen Poincare's critique of Hilbert's use of induction in the consistency proof which is sketched in \citeyearpar{Hilbert1904}, the publication of \textsl{Principia Mathematica}, and Hilbert's own engagement with the paradoxes, logicism, and type theory in light of \textsl{Principia}.\footnote{On which see \citep{Mancosu2003}.}   Hilbert reacted to these developments starting in his 1917-1918 lectures  `Prinzipien der Mathematik' by proposing that the axiomatic method he had developed for geometry should be extended to all of mathematics.  His lectures also take significant steps towards modern mathematical logic in terms of precision and a clear delineation between propositional, first-order, and higher-order logic.  

Hilbert was aided in these steps by Paul Bernays whom he had hired as his assistant for the foundations of mathematics in 1917.   Amongst Bernays's duties was to prepare protocols of Hilbert's lectures.  This informed his (second) \textsl{Habilitationschrift} \citeyearpar{Bernays1918a} on the propositional calculus in which Bernays formulated and proved soundness and completeness in the contemporary manner -- i.e. `Every provable formula is universally valid, and conversely' (p. 236).\footnote{Page references are to Bernays's \textsl{Habilitationschrift} as reprinted \citep{Hilbert2013} (which also contains Hilbert's lectures).   Hilbert had previously originally shown that the propositional calculus has the property now known as \textsl{Post completeness} -- i.e. if any unprovable formula is added to the axioms as a schema, then the calculus becomes inconsistent.   It is now natural to view this as a corollary of \textsl{deductive completeness} as formulated by Bernays -- i.e. every valid formula is provable.  However Post completeness bears an obvious resemblance to the property of \textsl{descriptive completeness}  which Hilbert had hoped to capture by his so-called \textsl{Vollst\"andigkeit} axiom in the second edition of the \textsl{Grundlagen der Geometrie} -- i.e. `To a system of points, straight lines, and planes, it is impossible to add other elements in such a manner that the system thus generalized shall form a new geometry obeying all of the five groups of axioms' (p. 25).   See \citep{Zach1999}, \citep{Detlefsen2014a}, and \citep[\S 11]{Baldwin2018} for more on the relationship between these and several other notions of completeness which are relevant to this historical context.} On this basis Bernays additionally observed that soundness implies that the calculus is consistent -- i.e. `it is not possible to derive both a formula and its negation' -- and \textsl{decidable} -- i.e. that there exists a `uniform procedure by which can one decide for every formula of the calculus after a finite number of steps whether or not the expression is provable' (p. 240).  

The major development separating  Bernays's propositional proof  from G\"odel's first-order one was the publication of Hilbert and Ackermann's textbook \textsl{Grundz\"uge der theoretischen Logik} in 1928.  This was G\"odel's source for his dissertation \citeyearpar{Godel1929a} and consolidates many of the advances made by Hilbert and his collaborators during  the 1920s.  These include an axiomatization of first- and second-order logic extending that which Bernays had used for propositional logic and the formulations of the formal definitions of provability ($\Gamma \proves \phi$) and consistency ($\Gamma \not\proves \bot$) for these systems.\footnote{In fact \citet[\S I.10, \S III.5]{Hilbert1928} provided  the original definition for what we now call the `Hilbert system' for first-order logic.} These formalisms are distinguished from Russelian type theory (e.g. with respect to the dispensability of the Axiom of Reducibility) and shown to be sufficient for the axiomatization of several portions of mathematics -- inclusive of fragments of geometry (\S III.10) and analysis (\S IV.8).

Hilbert and Ackermann also extended the notion of \textsl{universal validity} [\textsl{Allgemeing\"ultigkeit}] which had been employed by Bernays -- i.e. a statement `which always yields a true statement for any interpretation of the variables'  -- to first-order logic.   This is used to demonstrate the consistency of the first-order axioms and to illustrate the validity and non-validity of a number of specific first-order formulas via the construction of arithmetical models and counter-models (\S III.9-10).    No formal definition of truth in a model or of logical consequence is stated.  But the approach adopted for demonstrating first-order validity and consequence are all compatible with the contemporary definitions of $\models \phi$ and $\Gamma \models \phi$ now associated with Tarski -- i.e. $\phi$ is valid if it is true in all models and $\phi$ is a logical consequence of $\Gamma$ just in case it is true in models in which $\Gamma$ is satisfied.\footnote{See \citep[pp. 44-45]{Hilbert2013} for further discussion of how logical validity is understood in in these sources.}

The most enduring contribution of the \textsl{Grundz\"uge} may be Hilbert and Ackermann's explicit formulation of the completeness and decidability of their axioms for first-order logic as open problems.  These are respectively stated in the text as follows:

\begin{quote}
\footnotesize{Whether the axiom system is complete at least in the sense that really all logical formulas that are correct for all domains of individuals can be derived from it is an unsolved question. We can only say purely empirically that this axiom system has always sufficed for any application.  \hfill (p. 68)

The  decision problem [\textsl{Entscheidungsproblem}] is solved if one knows a procedure, which permits the decision of the universal validity or satisfiability of a given logical expression by finitely many operations. The solution of the problem of decision is of fundamental importance to the theory of all domains whose propositions can be logically described using finitely many axioms.\footnote{Appearing immediately after this passage in the first edition of the \textsl{Grundz\"uge} (pp. 74-76) is a example which Hilbert and Ackermann use to illustrate the importance of the \textsl{Entscheidungsproblem}  -- i.e. the question of whether what they call \textsl{Pascal's Theorem} (a projective form of Pappus's Theorem) is provable from Hilbert's Order, Incidence, and Parallels axioms  alone.   This is formulated in detail by first showing how Pascal's Theorem ($\phi$) and the axioms ($\Gamma$) and can be formulated in a first-order language similar to that described above.  In \citeyearpar{Hilbert1899} Hilbert had shown that $\phi$ is not  entailed by $\Gamma$ (i.e. essentially $\Gamma \not\models \phi$) in the course of showing the equivalence of Pascal's Theorem with the commutativity of multiplication in the so-called \textsl{skew field} constructed from a geometry satisfying the Axiom groups I, II, IV via the segment arithmetic described in \S14-\S15.    Hilbert and Ackermann go on to illustrate the significance of the \textsl{Entscheidungsproblem} by observing that if the problem were solved in the \textsl{positive}, then this rather involved construction could be replaced by an algorithmic demonstration that $\not\proves \bigwedge \Gamma \rightarrow \phi$.\label{epng}} \hfill (pp. 73-74)}
\end{quote}

As we have seen, Bernays had answered the analogous versions of these questions in the positive for the propositional calculus.   And while Church and Turing would later provide a negative answer to the latter (to which I will return in \S 4), G\"odel promptly answered the former positively in his dissertation.  As this was juxtaposed between the publication of the \textsl{Grundz\"uge} in 1928 and his announcement of his incompleteness theorems in late 1930, the context which G\"odel navigates in its introductory section is quite complex.   But as the following passage makes clear,  it is evident that he understood the significance of his result relative to the issues which had been brought into the open by Frege's exchange with Hilbert:\footnote{See the introduction to \citep{Godel1929a}, \citep{Kennedy2011a}, and \citep[\S 4.2]{Baldwin2018} for more on the context of G\"odel's dissertation.  It will also be useful to note here that the sort of \textsl{failure} of completeness for second-order theories which G\"odel anticipated in the K\"onigsberg lecture in which he first announced Theorem \ref{ct} will not be relevant here until \S 6.}

{\footnotesize
\begin{quote}
Here `completeness' is to mean that every valid [\textsl{allgemein giltige}] formula expressible in the restricted functional calculus can be derived from the axioms by means of a finite sequence of formal inferences. This assertion can easily be seen to be equivalent to the following: Every consistent axiom system consisting of only [first-order formulas] has a realization.   (Here `consistent' means that no contradiction can be derived by means of finitely many formal inferences.) The latter formulation seems also to be of some interest in itself, since the solution of this question represents in a certain sense a theoretical completion of the usual method for proving consistency $\ldots$ for it would give us a guarantee that in every case this method leads to its goal, that is, that one must either be able to produce a contradiction or prove the consistency by means of a model.  L. E. Brouwer, in particular, has emphatically stressed that from the consistency of an axiom system we cannot conclude without further ado that a model can be constructed. But one might perhaps think that the existence of the notions introduced through an axiom system is to be defined outright by the consistency of the axioms and that, therefore, a proof has to be rejected out of hand. \hfill \citep[p. 61]{Godel1929a}
\end{quote}
}

Although G\"odel mentions neither Frege nor Hilbert by name, there can be little doubt that by  `the usual method for proving consistency' G\"odel had in mind what \citet{Hilbert1934} would later called `the method of exhibition' (and has here been called MCT).  It is thus also significant that he appears to acknowledge the potential cogency of each of the following positions:
\begin{compactenum}[i)]
\item Consistency \textsl{automatically} entails the existence of a satisfying model -- i.e. if $\Gamma$ is consistent, then the existence of $\mathcal{M} \models \Gamma$ follows `\textsl{without further ado}'.
\item Consistency is \textsl{insufficient on its own} to entail the existence of a satisfying model -- i.e. even if $\Gamma \not\proves \bot$, then the existence of $\mathcal{M} \models \Gamma$ only follows `\textsl{with further ado}' (and thus potentially not at all).   
\end{compactenum}

We have seen that Hilbert at one time endorsed a version of i) whereas Frege himself emphatically stressed ii).   But independently of their philosophical debate, it is also important to keep in mind that the result which G\"odel described in the foregoing passage and then proceeded to \textsl{prove} can be understood to have the following form:\footnote{Indeed it was in this contraposed form that G\"odel stated and proved the completeness of first-order logic: `The completeness theorem that we must now prove states $\ldots$: \textsl{Every valid logical expression is provable}. Clearly, this can also be expressed thus: \textsl{Every logical expression is either satisfiable or refutable}, and we shall prove it in this form.' \citeyearpar[pp. 74-75]{Godel1929a}.   In the address in which G\"odel originally announced the Completeness Theorem he would also go on to briefly indicate an argument for why it is does not apply to higher-order logic \citeyearpar[p. 29]{Godel1930c}.    His argument goes via the now well-known observation that his \textsl{incompleteness} theorems can be combined with the provable to categoricity of the axioms of second-order arithmetic $[\mathsf{PA}_2$] to show that there can be standard second-order model of (presumably) consistent theories like $\mathsf{PA}_2 + \neg \mathrm{Con}(\mathsf{PA_2})$.   I will return to discuss the significance of this observation to contemporary considerations in \S 6 (see also note \ref{nsol} in regard its bearing on geometry).}

\begin{theorem} \label{ct}
If a set of axioms $\Gamma$ is consistent $(\Gamma \not\proves \bot)$, then there exists a model $\mathcal{M}$ in which $\Gamma$ is satisfied $(\mathcal{M} \models \Gamma)$.   \label{compthm}
\end{theorem}
As G\"odel evidently realized, such a result can naturally be taken to provide mathematical confirmation of the slogan `consistency implies existence'.   But whereas in \citeyearpar{Hilbert1900} Hilbert introduced the slogan almost glibly, we have also seen that by the mid-1920s he had come to an understanding on which this statement both admits to a precise mathematical formulation and also that relative to this formulation it requires proof.  On the other hand, G\"odel himself does not appear to  share the reservations with accepting the principle of the sort he attributes to Brouwer but we will see also came to be shared by Hilbert (at least in some sense).   But at the same time, it is also easy to appreciate why G\"odel would have wished to stress that the correctness of the slogan notwithstanding, there is indeed some `ado' involved with actually showing that every consistent set of first-order sentences has a model.    This point will be further borne out in light of the further metamathemtical analysis of the Completeness Theorem to be considered below.

\section{Bernays and method of arithmetization}

We have seen that Hilbert's consistency proofs in  \textsl{Grundlagen der Geometrie} took the form of model constructions whereby a set of axioms $\Gamma$ is shown to be consistent by describing a structure $\mathcal{M} \models \Gamma$.  In the case where $\Gamma$ includes at least the axioms of Incidence and Order, Hilbert observed that the domain of $\mathcal{M}$ must already be infinite \citeyearpar[p. 8]{Hilbert1971}.   However he also stressed that the models employed in demonstrations (\ref{mcr}i-ii) are all \textsl{countable}  (p. 41) and also obtained by closing the natural numbers under what he described as the `calculating operations' [\textsl{Rechnungsoperationen}] of addition, subtraction, multiplication, and division together with extraction of square roots (p. 29, p. 102).   Thus despite the fact that these models are naturally termed \textsl{analytical} (in the sense that their domain and relations are comprised of objects defined over the real numbers $\mathbb{R}$), it is also understandable why Hilbert's models are considered as  steppingstones to what are now known as \textsl{arithmetical models} (whose domain and relations are defined over the natural numbers $\mathbb{N}$).\footnote{Although it will be useful to retain the distinction between analytical and arithmetical models here, \citet{Bernays1967} used the latter term to describe the models constructed in the \textsl{Grundlagen der Geometrie} itself (as have other commentators).   It should be noted in this regard that models constructed in service of proving results like (\ref{mcr}i,ii) are only required to satisfy weak continuity principles like CCP.   Once stronger principles such as Hilbert's \textsl{Vollst\"andigkeitsaxiom} or the completeness axioms of Cantor or Dedekind are added, then the systems cease to have countable models presuming these principles are interpreted as second-order axioms (rather than first-order schema) and also that the so-called \textsl{standard semantics} for second-order logic is employed.   But Hilbert did not make further use of these principles in the \textsl{Grundlagen} and they are not needed for the results discussed here.    It will also become clear below it is also possible to view such stronger principles as being satisfied in (countable Henkin) models of second-order arithmetic of the sort now studied in Reverse Mathematics. \label{nsol}}

Hilbert's work in geometry was taken into account in a series of expository papers by Bernays in the 1920s-mid-1930s (several of which are reprinted in \citeyear{Bernays1976}). Taken together with the introductory chapters of the first volume of \textsl{Grundlagen der Mathematik} \citeyearpar{Hilbert1934}, these accounts locate this work relative to the mature form of what is now referred to as \textsl{Hilbert's program} in the foundations of mathematics.  I will return to some of the specific characteristics of the program in regard to geometry in \S 5.  But for the moment it will suffice to recall that its goal is typically described as that of proving the consistency of `infinitary mathematics' from the so-called `finitary standpoint' of `contentual [\textsl{inhaltliche}] mathematics'.  Infinitary mathematics can here be taken to correspond to systems of analysis or set theory of the sort in which Hilbert's work in geometry or mathematical physics can be formalized.   On the other hand, Hilbert paradigmatically characterized finitary mathematics in terms of concrete operations performed on `number signs' typified by  calculations involving the `construction and deconstruction' of numerals in stroke notation.  Such operations were often described as `intuitive' rather than `logical'.  But in mature formulations of the program, an effort was also made to characterize finitary reasoning axiomatically via systems resembling (but also occasionally extending) what is now known as \textsl{Primitive Recursive Arithmetic} [$\mathsf{PRA}$].\footnote{Recall that $\mathsf{PRA}$ is formulated in a first-order language containing terms for all primitive recursive functions and is axiomatized by their defining equations together with the induction schema for quantifier-free formulas.  See \citep[\S 1-\S 2]{Tait2005} and \citep[\S 2]{Dean2017b} for more on the relationship between Hilbert's finitism and $\mathsf{PRA}$.}

It was against this backdrop in which Bernays introduced the expression `the method of arithmetization' in his (1930) paper `Die Philosophie der Mathematik und die Hilbertsche Beweistheorie':\footnote{The expression `the method of arithmetization' is now often used to refer in a narrow sense to G\"odel's \citeyearpar{Godel1931a} so-called `arithmetization of syntax'.   But as is evident from this and related passages, Bernays understood the method in a broader sense which derives part of its motivation from Kronecker's general views about the fundamental role played by natural numbers in mathematics but also subsumes the specific techniques developed by Dedekind, Weierstrass, and Cantor for the arithmetizaiton of analysis.   He would later suggest in the first volume of \textsl{Grundlagen der Mathematik}  \citeyearpar[pp. 2-3, pp. 17-18]{Hilbert1934} -- wherein   `Methode der Arithmetisierung' appears as an \textsl{indexed term} -- that its sense can be precisified in terms of something akin to the modern notion of interpretability in a number theoretic theory (possibly of the second-order -- cf. \citealp[Sup. IV]{Hilbert1939}).   Bernays would go on to invoke the method repeatedly in his mathematical work -- e.g. when he presents G\"odel's  technique for coding syntax in arithmetic in the second volume \citeyearpar[\S 4]{Hilbert1939} he thus treats it as an \textsl{application} of the method of arithmetization to metamathematics rather than coextensive with the method itself.  But it should also be kept in mind that Bernays also drew attention to the limitations of the method -- e.g. in \citeyearpar[p. 152]{Bernays1941d} he writes that `the arithmetization of geometry in analysis and set theory is not complete' and in \citeyearpar[p. 65]{Bernays1970} that `concepts such as those of a continuous curve and of a surface, as developed especially in topology, can probably not be reduced to the idea of number'.  See, however, note \ref{proj} below.}
\begin{quote}
\footnotesize{We thus come to a differentiation between the elementary mathematical standpoint and a systematic standpoint that goes beyond it. This differentiation is not drawn artificially or merely ad hoc, but rather it corresponds to the duality of the points of departure that lead to arithmetic, namely, on the one hand the combinatorial activity with ratios in discrete quantities [\textsl{mit Verh\"altnissen im Diskreten}] and on the other hand the theoretical demand that is placed on mathematics from geometry and physics. The system of arithmetic does not emerge only from a constructive and intuitively contemplating activity [\textsl{konstruierenden und anschaulich betrachtenden T\"atigkeit}], but rather mostly from the task to grasp exactly and master theoretically the geometric and physical ideas of set, area, tangent, velocity, etc. The method of arithmetization is a means to this end.  In order to serve this purpose, however, arithmetic must \textsl{extend} its methodical point of view from the original elementary standpoint of number theory to a \textsl{systematic} view in the sense of the aforementioned postulates.  \hfill (p. 253)}
\end{quote}

Bernays here makes a number of points which can be understood not only to refine those of \citet{Hilbert1900} but which also bear on the legacy of the Frege-Hilbert controversy.  First, the central role of arithmetical models in providing consistency proofs is reaffirmed and generalized to include not just geometry but also applications of infinitary mathematics in physics.   Second, the importance of arithmetic itself comes into sharper focus as it is now assigned the dual role of both helping to characterize finitary mathematics -- i.e. `the combinatorial activity with ratios in discrete quantities' -- and also of `mastering theoretically' concepts from other domains such as geometry and physics.    Finally, it is acknowledged that in order to achieve this purpose  the `elementary standpoint' of number theory needs to be extended in order to encompass `the aforementioned postulates' -- a reference to \textsl{analysis} which Bernays takes to include  `the idea of an infinite totality [that] cannot be verified by intuition [but is graspable] only in the sense of an idea-formation.' \citeyearpar[p. 252]{Bernays1930}

By the late 1920s Hilbert and his collaborators had indeed developed additional proof-theoretic techniques for proving consistency to which I will return in \S4.  But within the realm of model-theoretic proofs, they also introduced a systematic distinction between two kinds of axioms sets: 
\begin{example}
\label{fininf}
\begin{compactenum}[i)]
\item Axioms $\Gamma$ which can be  be shown to be consistent in the manner of (MCT) by presenting a $\mathcal{M} \models \Gamma$ with a \textsl{finite domain}.
\item  Axioms $\Gamma$ for which it can be shown that any model $\mathcal{M} \models \Gamma$ must have an \textsl{infinite domain} (presuming that $\Gamma$ is satisfiable at all).
\end{compactenum}
\end{example}
\citet{Hilbert1928} observed that it is possible to extend the method of arithmetical interpretations which \citet{Bernays1918a} had employed to demonstrate the consistency of the propositional calculus to also demonstrate the consistency of the axioms for first-order logic in the manner of (\ref{fininf}i) by observing that the axioms are satisfied in a one element model and the rules preserve truth.   However they also observed that the significance of their consistency proof should not be overestimated as in the case where $\Gamma$ contains mathematical axioms satisfying (\ref{fininf}ii) it provides  `absolutely no assurance' of consistency.  About this circumstance they go on to observe
\begin{quote}
{\footnotesize This problem, whose solution is of fundamental importance for mathematics, is incomparably more difficult than the question dealt with here. The mathematical axioms actually assume an infinite domain of individuals, and there are connected with the concept of infinity the difficulties and paradoxes which play a role in the discussion of the foundations of mathematics.  \hfill \citeyearpar[p. 65-66]{Hilbert1928}}
\end{quote}

As \textsl{Grundz\"uge der theoretischen Logik} was the source for G\"odel's dissertation, such observations presumably informed what he meant by the `usual method of proving consistency' \citeyearpar{Godel1929a}.  But consideration of the contrast between (\ref{fininf}i) and (\ref{fininf}ii) also figured prominently in how the general consistency problem was framed by Hilbert and Bernays at the beginning of the first volume of the \textsl{Grundlagen der Mathematik}.\footnote{Although officially joint with Hilbert, Bernays is generally credited as being the primary author of \textsl{Grundlagen der mathematik}.  On the other hand, the views expressed in the introductory chapters of the first volume (\S 1 and \S 2) largely conform to Hilbert's published addresses and lectures from the around this time.  I will thus  adopt the following convention below: i) general philosophical views expressed in this portion of the first volume will be attributed jointly to Hilbert and Bernays; ii) the details of the development of metamathematics in the second volume -- inclusive of the proof of Theorem \ref{act} and attendant reflections -- will be attributed individually to Bernays.} Therein, the process of demonstrating consistency in the manner of i) is referred to as the `method of exhibition' according to which `the finite domain of individuals together with the graphs chosen for the predicates $\ldots$ constitutes a model in which we can concretely point out that the axioms are satisfied' \citeyearpar[p. 12]{Hilbert1934}.   Proofs of this sort are contrasted with those in which $\Gamma$ satisfies (\ref{fininf}ii) which they refer to as instances of \textsl{formal axiomatics} (the characterization of which I will discuss further in \S 5).   In this case, they remark
\begin{quote}
{\footnotesize We are therefore forced to investigate the consistency of theoretical systems without considering actuality, and thus we find ourselves already at the standpoint of formal axiomatics. $\P$
Now, one usually treats this problem -- both in geometry and the disciplines of physics -- with the \textsl{method of arithmetization}.  The objects of a theory are represented by numbers or systems of numbers and the basic relations by equations or inequations, such that, on the basis of this translation, the axioms of the theory turn out either as arithmetical identities or provable sentences $\ldots$ or as a system of conditions whose joint satisfiability can be demonstrated via arithmetical existence sentences $\ldots$ This approach presupposes the validity of arithmetic, i.e. the theory of real numbers (analysis). And so we come to ask ourselves what kind of validity this is. \hfill \citeyearpar[p. 3]{Hilbert1934}} 
\end{quote}
The last sentence of this passage again foretells of the developments in proof theory which are presented in the second volume of  \textsl{Grundlagen der Mathematik}.   But it also provides further details both on how the method of arithmetization operates and how it might be used to mediate between the infinitary structures -- consideration of which Hilbert \& Bernays took be necessary for the development of geometry and physics -- and their representation in forms accessible to intuition.  

Contemporary readers are also likely to see in such passages a yet more explicit precedent for the concept of \textsl{interpretability} as we have already seen is employed in Hilbert's geometric consistency proofs.  For not only are the domains of points and line in his models countable, but they also form subsets of the real numbers which are \textsl{definable} by formulas in the language $\mathcal{L}^2_{\mathsf{Z}}$ of second-order arithmetic (or `analysis' as it is also called) -- i.e. the language $\mathcal{L}^1_{\mathsf{Z}} = \{0,+,\times,<\}$ of first-order arithmetic supplemented with second-order quantifiers intended to range over the powerset of $\mathbb{N}$ -- as are the extensions of the non-logical expressions of $\mathcal{L}_H$.  For instance relative to a particular means of formalizing real numbers in this language, it is possible to view the model $\mathcal{P} \models \mathsf{HP}$ as arising from a particular set of $\mathcal{L}^2_{\mathsf{Z}}$-formulas $\pi(x,y),\lambda(x,y,z),\iota(A,\ell),\beta(A,B,C),\kappa_1(A,B,C,D),\kappa_2(A,B,C,D,E,F)$ which respectively define the sorts of points and lines over the minimal Pythagorean field $\mathbb{P} \subseteq \mathbb{R}$ and the corresponding relations of incidence, betweenness, and congruence.   

There are, however, a number of distinct notions which go by the name `interpretability' in mathematical logic. Some of these are syntactic and relate pairs of theories  and others are semantic and relate pairs of models.  Paradigmatic of these classes are the following definitions:\footnote{The first of these definitions is the notion now called \textsl{relative interpretability}.  This definition is traditionally credited to \citet[\S I.V]{Tarski1953}.  However a specific case of this definition was employed by  \citet{Ackermann1937} to show the consistency of $\mathsf{ZF}$ set theory with the Axiom of Infinity negated relative to first-order Peano arithmetic. Variants of the second definition can be found under different names in many sources.  But the consensus seems to be that model theoretic interpretability is a folklore notion with  historical antecedents including Descartes' identification of points in the Euclidean plane and pairs of real numbers on which Hilbert's constructions build.}
\begin{example}
\label{interp}
\begin{compactenum}[i)]
\item \textsl{A theory} $\mathsf{S}$ \textsl{is proof-theoretically interpretable in a theory} $\mathsf{T}$ just in case there is a map $(\cdot)^*: \mathcal{L}_{\mathsf{S}} \rightarrow \mathcal{L}_{\mathsf{T}}$ which associates the non-logical symbols of $\mathcal{L}_{\mathsf{S}}$ with formulas of $\mathcal{L}_{\mathsf{T}}$ defining objects of the same types, together with a \textsl{domain predicate} $\delta(x)$ such that if $\phi^*$ denotes the result of replacing the $\mathcal{L}_{\mathsf{S}}$-symbols in $\phi$ with their images under $(\cdot)^*$ in $\mathcal{L}_{\mathsf{T}}$ and restricting quantifiers by $\delta(x)$ then the following holds:
\begin{center} For all $\phi \in \mathcal{L}_{\mathsf{S}}$, if $\mathsf{S} \proves \phi$, then $\mathsf{T} \proves \phi^*$.
\end{center}
\item \textsl{A structure $\mathcal{A}$ is model-theoretically interpretable in a structure $\mathcal{B}$} with domain $B$ just in case there exists a map $(\cdot)^*: \mathcal{L}_{\mathcal{A}} \rightarrow \mathcal{L}_{\mathcal{B}}$ associating the primitive expressions of $\mathcal{L}_{\mathcal{A}}$ with those of $\mathcal{L}_{\mathcal{B}}$ and a domain predicate $\delta(x)$ as in (\ref{interp}i) such that the $\mathcal{L}_{\mathcal{A}}$-structure $\mathcal{A}^*$ with domain $A^* = \{a \in B : \mathcal{B} \models \delta(a)\} \subseteq B$ and with non-logical symbols similarly interpreted in $\mathcal{B}$ by their images under $(\cdot)^*$ -- e.g. $P^{\mathcal{A}^*} = \{\vec{a} \in B^k : \mathcal{B} \models P^*(\vec{a})\} \subseteq B^k$ for a $k$-ary predicate $P$ --  is such that $\mathcal{A}^*$ is isomorphic to $\mathcal{A}$.
\end{compactenum}
\end{example}

A number of authors have suggested that Hilbert's geometric consistency and independence proofs may be reconstructed in terms of one or the other of these notions.\footnote{See, e.g., \citep[\S 3.3]{Hallett2010} for a reconstruction using definition (\ref{interp}i) and \citep[\S 4.3]{Eder2018} for a reconstruction using (\ref{interp}ii).}  But the pattern of the proofs in \citeyearpar{Hilbert1899} does not exactly conform to either of these templates.  For while we have seen that Hilbert describes the systems of geometry which are to be \textsl{interpreted} with sufficient precision to treat them as first-order theories extending $\mathsf{HP}$, he did not concern himself with showing that their images under the translation scheme described above are \textsl{provable} in an interpreting theory of analysis (as is required for reconstruction via \ref{interp}i).   And while we have seen that Hilbert's model constructions are also specified with adequate precision to treat them as defining \textsl{interpreting} analytical models, he says little to suggest that he viewed geometric theories as coming along with synthetic models which are \textsl{interpreted} in such structures (as is required for reconstruction via \ref{interp}ii).   

It is perhaps possible to read into Hilbert's failure to substantially discuss synthetic interpretations of geometry a number of morals relevant to the Frege-Hilbert controversy.\footnote{Although in the introduction to \citeyearpar{Hilbert1899} Hilbert describes geometry as `equivalent to the logical analysis of our perception of space' (p. 3), the term `spatial geometry' [\textsl{r\"aumlichen Geometrie}] is then immediately repurposed as an expression to denote arbitrary models of Hilbert's axioms for planes as well as points and lines.   Notably absent from \textsl{Grundlagen der Geometrie} are thus descriptions of `Euclidean space', `projective space', etc. which Frege presumably thought should have accompanied his axiomatic systems to describe their intended interpretations.   On the other hand a rationale for abandoning spatial intuition in the axiomatic development of geometry is provided in the introduction to the 1898-1899 lecture course \citep[pp. 221-223]{Hilbert2004} on which \citeyearpar{Hilbert1899} is most closely based.} But already in \citeyearpar[\S 13]{Hilbert1899} Hilbert had given a version of the axioms for analysis he would later present in \citeyearpar{Hilbert1900a}.   Additional evidence that he ultimately came to understand his geometric results in terms of something like proof-theoretic interpretability is provided by the  prior passage from \citeyearpar{Hilbert1934} wherein it is explicitly noted that geometrical axioms may be translated into \textsl{provable }arithmetical statements. This is exactly what is required for a proof-theoretic interpretation but is left out in \citeyearpar{Hilbert1899}.    

Separating \textsl{Grundlagen der Geometrie} from \textsl{Grundlagen der Mathematik} is not only the work in metamathematics surveyed in \S 2 but also Hilbert's engagement with the foundational problems specifically posed by infinitary mathematics.    The \textsl{locus classicus} for his discussion of such matters is his address `\"Uber der Unendliche' \citeyearpar{Hilbert1925}.   This was also based on a lecture course with the same name which Hilbert had given in G\"ottingen in 1924-1925 in which a broad range of themes are developed.    But what is most relevant here is how the contrast between cases (\ref{fininf}i) and (\ref{fininf}ii) is repeatedly used to motivate Hilbert's finitist consistency program as described above.   The following passage provides a summary of Hilbert's conclusions in this regard:\footnote{A related passage from the first volume of \textsl{Grundlagen der Mathematik} \citeyearpar[pp. 15-16]{Hilbert1934} is as follows: `In view of this difficulty in proving consistency, we could now try to use some other infinite domain of individuals which is not a mere product of thought (such as the number series), but is taken from the realm of sense perception or physical reality. If we take a closer look, however, we realize that wherever we believe that we encounter infinite manifolds in the realm of qualia or in physical reality, there can be no actual detection of such a manifold. The conviction of the existence of such a manifold actually rests on a mental extrapolation, which requires an examination of its justification at least, as necessarily as the conception of the totality of the number series.'}
\begin{quote}
{\footnotesize The final result then is: nowhere is the infinite realized; it is neither present in nature nor admissible as a foundation in our rational thinking -- a remarkable harmony between being and thought. We gain a conviction that runs counter to the earlier endeavors of Frege and Dedekind, the conviction that, if scientific knowledge is to be possible, certain intuitive conceptions and insights are indispensable; logic alone does not suffice. The right to operate with the infinite can be secured only by means of the finite. \hspace*{1ex} \hfill \citeyearpar[p. 392]{Hilbert1925}} 
\end{quote}
Hilbert reaches this conclusion after recording a long list of examples which testify to the central role played by infinitary concepts and structures in contemporary mathematics.   In the case of geometry he observes both that the Euclidean axioms lead to the assumption that space is infinite and also that introduction of `ideal elements' (e.g. points at infinity) makes the axiomatization of geometry `as simple and perspicuous as possible' (p. 373).   But he also observes that current physical theory leaves open the possibility that physical space is bounded in extent before concluding that  `Euclidean geometry, as a structure and a system of notions, is consistent in itself, but this does not imply that it applies to reality' (pp. 371-372). Such considerations are at least suggestive of why Hilbert may have originally demurred from offering a synthetic characterization of axiom systems like $\mathsf{HP}$ which do not have finite models.    But they also give rise to the question of whether by the mid-1920s he could have continued to endorse the slogan `consistency implies existence' while also allowing that a theory could be `consistent in itself' while failing to possess a model `in reality'?

Bernays's awareness of this issue is also vividly attested in the following note he made at around this time:\footnote{This is part of a note found in Hilbert's \textsl{Nachlass} in G\"ottingen [Cod. 685:9, 2] with  the title `\textsl{Existenz und Widerspruchsfreiheit}'.   It has been reproduced in the original and in translation by \citet[p. 479]{Sieg2002a} who dates it to between 1925 and 1928.}
\begin{quote}
{\footnotesize  The claim: ``Existence = consistency'' can only refer to a system as a whole. Within an axiomatic system the axioms decide about the existence of objects. $\P$ If, for a system as a whole, consistency is to be synonymous with existence, then the proof of consistency must consist in an exhibition.  $\P$ (All consistency proofs up to now have been either direct exhibitions or indirect ones by reduction; in the latter case a certain other system is already taken as existent.  -- \textsl{Frege} has defended with particular emphasis the view that any proof of consistency has to be given by the actual presentation of a system of objects). $\P$ In proof theory, laying a new foundation of arithmetic, consistency proofs are not given by exhibition. From this foundational standpoint it does not hold any longer that existence equals consistency. Indeed, it is not the opinion that the possibility of an infinite system is to be proved, rather it is only to be shown that operating with such a system does not lead to contradictions in mathematical reasoning.}
\end{quote} 

As should now be evident, the tenability of the restrictive view about mathematical existence described here is at least called into question by G\"odel's Completeness Theorem. But as the foregoing passages also make clear, Hilbert and Bernays were at this time operating from within a framework which took as central the interaction between epistemological concerns about how we can come to know that various theories are consistent and ontological ones about what it would mean for a satisfying model to exist in the case that they are.  It is striking how far removed these considerations are from the concerns and presuppositions of contemporary model theory.  For in contemporary practice doubts about the consistency of theories or the ontological status of the models largely given way to the goal of applying specific model theoretic techniques  -- e.g. quantifier elimination, omitting types, infinitary logic -- to study questions internal to specific branches of infinitary mathematics.\footnote{The extent of this transformation is aptly illustrated by \citet{Baldwin2018} via the applications of work in stability and classification theory to algebraic geometry and combinatorics.}

But as we shall see in \S 4 and \S 5, the legacy of Hilbert and Bernays's work in mathematical logic is more directly related to the development of \textsl{computability theory} than it is to \textsl{model theory} in the manner in which these subjects are now understood.   In addition to G\"odel's incompleteness theorems, some of the early work in this area is already taken into account in the second volume of \textsl{Grundlagen der Mathematik}.  This contains a novel result of Bernays which has come to be called the \textsl{Arithmetized Completeness Theorem}.    A simple formulation of this result is as follows:
\begin{theorem}
\label{act}
Suppose that $\Gamma$ is a finite set of sentences over an arbitrary first order language $\mathcal{L}_{\Gamma}$ and that $\Gamma$ is consistent -- i.e. $\Gamma \not\proves \bot$.   Then $\Gamma$ is model-theoretically interpretable in the standard model of first-order arithmetic -- i.e. the $\mathcal{L}^1_{\mathsf{Z}}$-structure $\mathcal{N} = \langle \mathbb{N},+,\times,<,0\rangle$.   
\end{theorem}
\noindent Suppose we now officially define an \textsl{arithmetical model} for a first-order language $\mathcal{L}_{\Gamma}$ -- which, for ease of illustration, we assume contains only predicate letters $P_1,\ldots,P_k$ -- to be a structure $\mathcal{M} = \langle A,R_1,\ldots,R_n \rangle$ such that $A \subseteq \mathbb{N}, R_i \subseteq \mathbb{N}^{r_i}$ (where $r_i$ is the arity of $P_i$) and also that each of these sets is \textsl{arithmetically definable} -- i.e. there exist $\mathcal{L}^1_{\mathsf{Z}}$-formulas $\delta(x),\psi_1(\vec{x}), \ldots, \psi_n(x)$ such that $A = \{n \in \mathbb{N} : \mathcal{N} \models \delta(n)\}$ and $R_i = \{\vec{n} \in \mathbb{N}^{r_i} : \mathcal{N} \models \psi_i(\vec{n})\}$.  Theorem \ref{act} is then equivalent to the statement: \textsl{Every finite consistent set of first-order sentences $\Gamma$ has an arithmetical model}.  

Theorem \ref{act} is, of course, a version of the L\"owenheim-Skolem Theorem -- i.e. every consistent set of sentences $\Gamma$ over a countable language has a countable model.  G\"odel's \citeyearpar{Godel1929a} original proof of completeness established this directly by constructing a model $\mathcal{M} \models \Gamma$ which has either domain $\mathbb{N}$ (in the case that $\Gamma$ has at least one infinite model) or otherwise has a domain consisting of finitely many equivalence classes over $\mathbb{N}$ (in the case that $\Gamma$ has only finite models).   Bernays's result strengthens this observation by showing how G\"odel's method can be formalized in a manner which allows for the construction of arithmetical formulas defining the non-logical symbols in $\mathcal{L}_{\Gamma}$ so that the resulting collection of arithmetical sets forms a model of $\Gamma$.   In \S 5 I will return to the question of whether this construction can additionally be viewed as yielding a finitary reduction of the problem posed by case \ref{fininf}ii) -- e.g. in that of proving the consistency of theories which can be shown to possess no finite models.

\section{On the difficulty of consistency}

In the prior two sections I have laid out a case for the first of my initial theses (\ref{points}i) -- i.e. that G\"odel's Completeness Theorem was not only obtained in light of the concerns brought to the fore by the Frege-Hilbert controversy but that these considerations also contributed to its subsequent reception.   In this section, I will present a case for thesis (\ref{points}ii) -- i.e. that Frege was correct to anticipate that there is a precise sense in which the problem of determining the consistency of a set of sentences (so that, e.g., completeness might be applied) is \textsl{as difficult as it can be}.   Laying this out will require indicating an analysis of what it means for a mathematical problem to be \textsl{difficult} -- an account I will also employ in \S 5 to provide an account of the relative \textsl{ease} of demonstrating existence given consistency. Although the proposal I will present is general,  it will still be useful to approach it via the historical frame of the preceding sections.\footnote{Although informal attributions of `ease' and `difficulty' are common in mathematical practice, philosophers have thus far paid little attention to their basis or the formal judgement or their relationship to notions such as \textsl{hardness}, \textsl{completeness}, or \textsl{degree of difficulty} (or \textsl{unsolvability}) studied in subjects such as computability theory and complexity theory.    The account developed here should be understood as exemplifying how the notion of \textsl{problem difficulty} from \citep{Dean2019} -- i.e. that of deciding an infinite class of yes-no questions about set membership -- relates to that of \textsl{propositional difficulty} -- i.e. that of determining the truth or falsity of a single mathematical statement.} 

A first consideration emerges from Frege's original misgivings about Hilbert's goal of obtaining a `direct' consistency proof for arithmetic or analysis.  We have seen that Hilbert stated this goal in \citeyearpar{Hilbert1900} but provided few details as to the form which such a proof might take.   Although this announcement itself played a role at the end of his correspondence with Hilbert, Frege also commented extensively on the nature and role of consistency proofs in both his mathematical and philosophical work.  I will now suggest that these accounts collectively suggest that he identified \textsl{two dimensions of difficulty} involved in such proofs.  

The first dimension can be introduced by recalling that a central part in Frege's own philosophy of mathematics was played by what we would now call \textsl{conceptual analysis}.   In fact a key feature of his logicism was the \textsl{denial} of the view that the surface grammar of textbook statements of mathematical theorems and other principles provide a reliable guide to the non-logical primitives which should appear in a fundamental axiomatization of their subject matter.   He rather held that the notions expressed by grammatically simple expressions employed in practice are often amenable to analysis which reveals that they admit to decomposition into more elementary concepts.    Such a process was to be repeated, potentially uncovering multiple conceptual strata which -- to suit Frege's ends -- ought to be of an increasingly logical character.  

Frege famously applied this method to arithmetic in the program he announced in his  \textsl{Begriffsschrift} (1879) and then developed more fully in his \textsl{Grundlagen} (1884) and the two volumes of his \textsl{Grundgesetze} (1893, 1903).  The details and exigencies of Frege's proposed analysis of number theoretic statements in terms of equicardinality and the analysis of the latter in terms of the extensions of concepts are sufficiently well known that there is no need to enter into them here.\footnote{For present purposes an equally illustrative example is provided by Frege's \citeyearpar[pp. 27-32]{Frege1880} proposal to demonstrate the statement `the sum of two multiples of a number is in its turn a multiple of that number' -- i.e. $\A w \A x \A y \E z(w \times x + w \times y = z \cdot w)$ -- from axioms stating the associativity of addition -- i.e. $\A x \A y \A z((x+ y) + z = x + (y +z))$ -- and that $0$ is the additive identity -- i.e. $\A x(x = x + 0)$.   Frege argues that it is possible to derive this statement without `presupposing any multiplication theorem' or even `the concept of multiplication' by analyzing the concept expressed by `multiple of' in terms of that expressed by `$x$ is a multiple of $y$ iff $x$ appears in the sequence $0, y, y+y, (y + y) + y, \ldots$'  This analysis can potentially be further decomposed using the notion of `hereditary in the $f$-sequence' which Frege attempts to characterize logically in \S III of his \textsl{Begriffsschrift}.  See \citep[\S 1]{Blanchette2012} for further reconstruction of this example.}  As Frege's general views about content evolved over time, it is not easy to give a concise account of the circumstances under which he would have regarded a proposed analysis as preserving the proposition (or \textsl{thought}) expressed by a sentence.   But what is most significant here is that he maintained the centrality of conceptual analysis throughout the evolution of his program and also that he proposed this method should be applied not just to arithmetical discourse but that of other sciences as well.  Both points are evident from the following passage from his late paper `Logic in mathematics':
\begin{quote}
{\footnotesize In the development of a science it can indeed happen that one has used a word, a sign, an expression, over a long period under the impression that its sense is simple until one succeeds in analysing it into simpler logical constituents. By means of such an analysis, we may hope to reduce the number of axioms; for it may not be possible to prove a truth containing a complex constituent so long as that constituent remains unanalysed; but it may be possible, given an analysis, to prove it from truths in which the elements of the analysis occur.  \hfill \citeyearpar[p. 209]{Frege1914a}}
\end{quote}

Frege's remark in his initial letter to Hilbert  that he had earlier considered an axiomatization of geometry in which he could `make do with fewer primitive terms' [IV/3, p. 34] comes into sharper focus in this context.  For it highlights the possibility that Hilbert's consistency and independence proofs might be invalidated in virtue of a given $\mathcal{L}_{H}$-statement becoming provable or refutable once the concepts expressed by its non-logical terms are subjected to proper analysis.\footnote{In fact \citet[p. 247]{Frege1914a} provides the expected case in point by considering the possibility that the Parallel Postulate becomes derivable if additional relationships between the concepts expressed by `straight line', `parallel', and `intersect' are taken into account by the geometric axioms.} In voicing such concerns, Frege isolates the following question illustrating the first of the two dimensions of difficulty alluded to above:\footnote{This formulation presupposes that consistency is a notion which is appropriately applied to sets of propositions in addition to their expressions as sentences in a fixed language.   Although Frege most often speaks of consistency of \textsl{concepts} specified by sets of defining predicates, it seems likely that he would not have objected to these other formulations.  But of course it is also possible to formulate a question similar to (D1) where we ask not after the consistency of the set of propositions $\Gamma^*$ but rather after that of a set of sentences $\Gamma^+$ which is obtained from $\Gamma$ by systematically replacing the non-logical terms of $\mathcal{L}_{\Gamma}$ with expressions in $\mathcal{L}_{\Gamma^+}$ which provide linguistic formulations of appropriate conceptual analyses.   Such a characterization takes a step towards reformulating (D1) in a manner which could potentially be assimilated to the second dimension of difficulty (D2) discussed below -- e.g. by asking after how hard it is to determine whether there exists a proof-theoretic interpretation of $\Gamma$ in a theory $\Gamma^+$ satisfying various adequacy conditions on successful conceptual analyses.  A related proposal is given by \citet{Eder2015} in his reconstruction of Frege's unpublished (1906) attempt to understand Hilbert's independence proofs on his own terms.}
\begin{itemize}
\item[(D1)]  Suppose we are given a set of sentences $\Gamma$ formulated over a language $\mathcal{L}_{\Gamma}$ and a definition of formal derivability from axioms $\Gamma \proves \phi$.  Let $\Gamma^*$ denote the thoughts expressed by the sentences in $\Gamma$.   How can we determine if the $\mathcal{L}_{\Gamma}$-formulation of the thoughts in $\Gamma^*$ provides an adequate analysis of the concepts (or \textsl{senses}) which they contain so that the formal consistency of $\Gamma$ (i.e. $\Gamma \not\proves \bot$) provides a sufficient condition for the consistency of the thoughts $\Gamma^*$?
\end{itemize}

It is easy to adduce examples which illustrate that this is not a mathematically idle concern.\footnote{In fact a paradigm example is provided by Tarski's \citeyearpar{Tarski1959a} later axiomatization of geometry in a language which contains only a sort for points and equidistance and betweenness predicates.   Not only might such an axiomatization be regarded as simpler than Hilbert's but it also has different metatheoretic properties -- e.g. Tarski's theory is unlike $\mathsf{HP}$ in that it is deductively complete and hence decidable (as I will discuss in \S 6).   But as has been stressed by \citet{Kreisel1967}, it is also possible to look at subsequent debates about set theoretic independence through the lens of (D1) -- e.g. can further application of conceptual analysis (or `informal rigour') to the concepts \textsl{set} and \textsl{membership} be used to justify axioms which decide formally independent statements such as the Continuum Hypothesis?} But as Frege himself admitted,  it is at best unclear whether (D1) admits to a definitive answer even in specific cases.  For not only are we left with the question of when a proposed conceptual analysis of an expression is \textsl{correct}, but there is also the question of whether there is a means of determining whether we have reached an `ultimate' strata of concepts beneath which no further analysis is possible.  And in both cases, it seems that Frege simply did not provide an answer.\footnote{See, e.g., \citeyearpar[p. 209 ff.]{Frege1914a}.  But of course this issue is yet more complex since it interfaces in various ways with Frege's views about other difficult topics -- e.g. the composition of thoughts as sense complexes, the identity of sense, and the role of definitions in mathematics.  See, e.g., \citep[\S 2]{Blanchette2012}.}

The foregoing observations about Frege's outlook on consistency proofs are well-known.  But it is evident that Frege additionally anticipated another central difficulty about the notion of consistency which -- while itself well-known -- admits to a precise mathematical analysis while also engaging with his misgivings about the possibility of non-model-theoretic consistency proofs.   This dimension can be introduced by observing that in the course of highlighting (D1), Frege also called attention to the fact that even at a fixed level of analysis the problem of checking consistency is often non-trivial.  One of his formulations of this point is as follows:
\begin{quote}
{\footnotesize It is completely wrongheaded to imagine that every contradiction is immediately recognizable; frequently the contradiction lies deeply buried and is only discovered by a lengthy chain of inference.\\ \hspace*{1ex} \hfill \citeyearpar[p. 194/179]{Frege1906}}
\end{quote}

The use of the metaphorical expression `deeply buried' is typical of Frege's remarks about the sort of effort which is often required to provide revealing conceptual analyses.  But he also acknowledged that the process of analyzing concepts is itself distinct from that of deductively deriving consequences from their formulations within a fixed theory.\footnote{E.g. `The effect of the logical analysis of which we spoke will then be precisely this -- to articulate the sense clearly. Work of this kind is very useful; it does not, however, form part of the construction of the system, but must take place beforehand.'  \citeyearpar[p. 211]{Frege1914a}}  In using expressions such as `lengthy chain of inference' it would thus appear that Frege is pointing to another aspect of what often makes it difficult to determine if a given set of sentences is consistent.   One way of formulating the question underlying this second dimension is as follows:
\begin{itemize}
\item[(D2)]  Let $\Gamma$ be a fixed set of sentences and $\Gamma \proves \phi$ a notion of deductive consequence as above.   Given that formal derivations can be of unbounded  length -- and thus there is no \textsl{a priori} bound which can be placed on the length of the derivation of a contradiction -- how can we decide if $\Gamma \not\proves \bot$?   
\end{itemize}

It is again easy to adduce examples which illustrate that this is not a mathematically idle concern.\footnote{Frege made the prior remark as part of a reaction to Shoenflies's discussion of Russell's paradox -- a contradiction which had famously afflicted the system of his \textsl{Grundgesetze}.   But although we now do not think of this contradiction as `lying deeply buried' itself, the history of mathematical logic provides other examples where axiom systems have been found to be inconsistent only in virtue of lengthy (and sometimes non-obvious) derivations.   For instance, \citet{Curry1941} observed that the complete derivation of a contradiction discovered by Kleene and Rosser in one of Church's original formulations of the lambda calculus runs to 162 published pages.  (In fact it was in an effort to simplify this result by which \citet{Curry1942} was led to what we now call `Curry's paradox'.)  \citet{Rosser1942} similarly showed that the system in the first edition of Quine's \textsl{Mathemtaical Logic} was inconsistent by presenting a dense 15-page derivation of the Burali-Forti paradox.   But of course the true force of (D2) is best illustrated by the existence of axiom systems whose consistency is still taken to be an open problem -- e.g. $\mathsf{ZF}$ plus the existence of a Reinhardt cardinal or Quine's New Foundations (Holmes's recently claimed consistency proof remaining as yet unpublished).}  But note that in the case where $\Gamma$ is a finite set of first-order sentences, it is also straightforward to see that (D2) is equivalent to Hilbert and Ackermann's \textsl{Entscheidungsproblem} -- i.e. to the question of determining in the general case whether $\phi$ is derivable from $\Gamma$.  For in this case we then have following sequence of biequivalences: 
\begin{example}
\label{red1}
 $\Gamma \proves \phi$ iff $\Gamma \proves \neg \neg \phi$ iff $\Gamma \proves \neg \phi \rightarrow \bot$ iff $\Gamma \cup \{\neg \phi\} \proves \bot$
\end{example}
Thus if we possessed a general method for determining whether an arbitrary finite set of first-order sentences is consistent, then we could also use this method to decide if $\Gamma \proves \phi$ by determining whether it yields a negative answer when applied to $\Gamma \cup \{\neg \phi\}$.

Although Frege's reservations about Hilbert's consistency proofs are traditionally explained in terms of his preoccupation with (D1), the previous passage suggests that he was also well aware of the challenge posed by (D2).  Additional evidence to this effect is provided by his more extensive discussion of consistency proofs in the second volume of his \textsl{Grundgesetze}.   Therein we find the following highly germane passage:

\begin{quote}
{\footnotesize 
How is it to be recognised that properties do not contradict each other? There seems to be no other criterion than to find the properties in question in one and the same object $\ldots \P \ldots$ Or is there perhaps a different way to prove the freedom from contradiction? If there were, this would be of the highest significance for all mathematicians who ascribe the power of creation to themselves. And yet hardly anyone seems concerned to find such a method of proof. Why not? Probably because of the view that it is superfluous to prove freedom from contradiction since any contradiction would surely be noticed immediately. How nice if it were so! How easy all proofs would then be!
The proof of the Pythagorean theorem would then go as follows:
\begin{quote}
``Assume the square of the hypotenuse is not of equal area with the squares of the two other sides taken together; then there would be a contradiction between this assumption and the familiar axioms of geometry. Therefore, our assumption is false, and the square of the hypotenuse is of an area exactly equal to the squares of the two other sides taken together.''
\end{quote}}

{\footnotesize
$\ldots$ Absolutely any proof could be conducted following this pattern. Unfortunately, the method is too easy to be acceptable. Surely, we see that not every contradiction lies in plain view. Moreover, we lack a sure criterion for the cases where from the non-obviousness of a contradiction we may infer its absence.  \hfill \cite[\S143-\S144]{Frege1903}
}
\end{quote}

Like \citet{Hilbert1900}, Frege begins here by calling attention to the fact that consistency proofs have traditionally proceeded by exhibiting models.  But rather than taking this as an incentive to develop a new technique for proving consistency which avoid this exigency, Frege then appears to offer a \textsl{reductio} of the supposition that such a method might exist.\footnote{\citet[p. 237]{Frege1897} also writes that `In my conceptual notation inference is conducted like a calculation $\ldots$ in the sense that there is an algorithm $\ldots$ which govern[s] the transition from one sentence or from two sentences to a new one in such a way that nothing happens except in conformity with the rules'.   The issue he is addressing in the prior passage is thus not whether the step-by-step process of derivation is effective, but rather that of whether there is a general algorithmic method for determining whether a statement is derivable from axioms by a proof of arbitrary length.} One way of reconstructing his argument is as follows:
\begin{example}
\begin{compactenum}[i)]
\item Suppose there were a mathematical method which allowed us to determine whether an arbitrary set of sentences $\Gamma$ is consistent without exhibiting a satisfying model -- i.e. an effective procedure $\alpha$ which when applied to a specification of $\Gamma$ allowed us to determine if $\Gamma \not\proves \bot$ in a finite number of steps, taking into account both that no general bound can be placed on the length of a proof of a contradiction and also that $\Gamma$ may not possess a finite model which would rule out the existence of such a proof.
\item In virtue of (\ref{red1}), the general task of deciding if an arbitrary mathematical sentence $\phi$ follows $\Gamma$ could then be replaced by using $\alpha$ to decide if $\Gamma \cup \{\neg \phi\}$ is consistent.
\item In practice, we find that the task of deciding whether a statement follows from a set of axioms is of considerable difficulty -- e.g. while we know that the Pythagorean Theorem is derivable from the axioms $\mathsf{HP}$, this is only so in virtue of the fact that we have expended effort in finding a proof, not as a result of applying a general algorithmic decision procedure.\footnote{In the passage cited above, Frege offers both the Pythagorean Theorem and the Law of Quadratic Reciprocity from number theory as examples of this phenomena.   Although there are many well-known proofs of both statements, they are each typically taken to be at least somewhat involved or non-obvious in their respective domains.   They both can thus both be taken to illustrate instances in which what \citet[p. 376]{Detlefsen1990} calls the \textsl{inventional complexity} of a statement --  i.e. that `encountered in coming up with a proof in the first place' -- would be trivialized if there were a general method for proving consistency of the sort which Frege considers.}
\item The existence of $\alpha$ as described in i) would thus trivialize a problem which we know to be difficult in practice.  We can thus conclude that no such method can exist.  
\end{compactenum} 
\end{example}

An obvious rejoinder to this argument is that our ignorance of a general method for deciding consistency -- however well confirmed in mathematical practice -- is not itself sufficient to conclude that no such method can exist.    But in juxtaposing these concerns, Frege can also be seen as anticipating by more than 30 years another well-known turn of events in mathematical logic.   For as we now know, Church and Turing answered the \textsl{Entscheidungsproblem} in the negative in 1936.   Relative to the widely accepted analysis of `effective procedure' as `Turing computable function' (to which I will return in a moment), Frege was thus indeed prescient in suspecting that no general method for checking the consistency of finite sets of fist-order sentences can exist.

This is often presented as a \textsl{negative} (or ``limitative'') result.   But it may also be contrasted with several \textsl{positive} results pertaining to special cases of the \textsl{Entscheidungsproblem} and to `direct' consistency proofs which had previously been obtained by Hilbert and his collaborators.  By the 1920s, Hilbert had come to realize that a first step towards a positive characterization of this notion already follows from the introduction of a formal definition of derivability -- i.e. what has here been denoted here by $\Gamma \proves \phi$.  For in the course of defining such a relation, one also defines a class of derivations $\mathfrak{Der}$ whose members  $\mathfrak{D}_0,\mathfrak{D}_1,\ldots$ are what we typically call `proofs' --  e.g. in the case of `Hilbert proofs' for first-order logic as defined by \citet{Hilbert1928} a finite sequence of formulas whose elements are either members of $\Gamma$, logical axioms, or follow from earlier statements by rules such as \textsl{modus ponens}.    

Once such a definition of derivability is in place, the definition of the consistency of $\Gamma$ then takes the form
\begin{example}
\label{derdefn}
$\Gamma \not\proves \bot$ if and only if for all $\mathfrak{D} \in \mathfrak{Der}$, it is not the case that $\mathfrak{D}$ has hypotheses contained in $\Gamma$ and conclusion $\bot$.
\end{example}
Hilbert repeatedly called attention to the fact that the structure of this definition makes clear that to prove the consistency of $\Gamma$ it suffices to demonstrate a \textsl{universal} assertion about finite combinatorial objects (i.e. proofs) rather than an \textsl{existential} one about the existence of a (possibly infinite) satisfying model.\footnote{E.g. `To prove consistency we therefore need only show that $0 \neq 0$ cannot be obtained from our axioms by the rules in force as the end formula of a proof, hence that $0 \neq 0$ is not a provable formula. And this is a task that fundamentally lies within the province of intuition just as much as does in contentual number theory the task, say, of proving the irrationality of $\sqrt{2}$, that is, of proving that it is impossible to find two numerals $a$ and $b$ satisfying the relation $a^2 = 2b^2$, a problem in which it must be shown that it is impossible to exhibit two numerals having a certain property.  Correspondingly, the point for us is to show that it is impossible to exhibit a proof of a certain kind.'  \citeyearpar[p. 471]{Hilbert1927} }   And on this understanding, he and his collaborators did in fact make a considerable amount of progress in the 1920s and early 1930s towards the development of what can  be regarded as direct methods of proving consistency.

The methods in question are typified by the technique of \textsl{cut elimination} which pertains to an alternative proof system for first-order logic known as the \textsl{sequent calculus} originally due to Gerhard Gentzen.\footnote{As a means for proving consistency, cut elimination is related to (but arguably more general than) the slightly earlier methods of $\varepsilon$-\textsl{substitution}  -- i.e. Hilbert's \textsl{Ansatz} \citeyearpar{Hilbert1923} -- and also that of \textsl{expansion} or \textsl{r\'eduite} -- i.e. Herbrand's Theorem \citeyearpar{Herbrand1930a}.   These techniques are connected in such a way that most of the general points made below about `direct consistency proofs' could be formulated in terms of each of them.    See, e.g., \citep{Rathjen2018} for an overview of their technical and historical relationship.}  Gentzen showed in his dissertation (written under Bernays) that there existed an effective procedure $\gamma$ whereby applications of the so-called \textsl{cut rule} (a generalized form of \textsl{modus ponens}) can be successively eliminated such that a derivation $\mathfrak{D}$ of a sequent $\Gamma \Rightarrow \Delta$ can be transformed into a so-called \textsl{cut-free} derivation $\gamma(\mathfrak{D}) = \mathfrak{D}^*$ of the same sequent.  Cut-free derivations $\mathfrak{D}^*$ can also be shown to possess the so-called \textsl{subformula property} -- i.e. if $\mathfrak{D}^*$ is a cut-free derivation of $\Gamma \Rightarrow \Delta$, then all of the formulas it contains are subformulas of formulas appearing in either $\Gamma$ or $\Delta$ (a feature which is potentially violated by the cut rule).  A consequence is that the sequent $\emptyset \Rightarrow \emptyset$ is not derivable in Gentzen's system which, on the intended interpretation, corresponds to the non-derivability of a contradiction from no premises.  In this way cut elimination yields another proof of the consistency of pure first-order logic (i.e. $\not\proves \bot$)  which can be contrasted with Hilbert \& Ackermann's inductive consistency proof using arithmetical models. 
 
Such a proof can be carried out either by directly reasoning about derivations as finite combinatorial objects or via a well-known formalization within in a fragment of $\mathsf{PRA}$.\footnote{See, e.g., \citep[\S 3c, \S 5d]{Hajek1998}.}   In both cases the sort of reasoning involved is paradigmatically finitary in the sense described by \citet{Hilbert1934} -- e.g. not only can the cut elimination operation $\gamma$  be understood as acting on `concrete' (or `intuitively accessible') proofs, proving that it has the requisite properties does not require reasoning involving unbounded quantifiers over the natural numbers.    In this way Gentzen's method is prototypical of the sorts of positive results which Hilbert hoped would ultimately be capable of supplying a consistency proof for analysis.\footnote{Hilbert's remarks on Ackermann's attempted consistency proof \citeyearpar[pp. 477-479]{Hilbert1927} via the $\varepsilon$-substitution method are prototypical here.   But the context is complicated by the fact that once Ackermann's subsequent modifications to the proof are take into account, the scope of his result is more limited than Hilbert appears to have assumed.  See, e.g., the introduction to Bernays's commentary on \citeyearpar{Hilbert1927} in \citep[p. 485-486]{Heijenoort1967}.} 
 
As mathematical axioms typically disturb the symmetric nature of the logical rules on which Gentzen's method depends, the foregoing method on its own is not sufficient to demonstrate the consistency of substantial fragments of analysis or even arithmetic.\footnote{This obstacle can be overcome in certain cases by the techniques of \textsl{partial cut elimination} or \textsl{ordinal analysis} (as introduced by Gentzen in his 1936 consistency proof for first-order arithmetic).  However these methods are most directly applicable to theories in which arithmetic can be interpreted and are thus only indirectly relevant here.}   These difficulties notwithstanding, I will describe below how traditional proof-theoretic techniques can be combined with more recent methods from Reverse Mathematics to prove the consistency of the specific geometric theories with which Hilbert was concerned in \citeyearpar{Hilbert1899} in a manner which is arguably `direct'.\footnote{G\"odel's second  incompleteness theorem also places well-known constraints on the relationship between $\Gamma$ and a theory $\mathsf{T}$ in which a potential consistency proof can be formalized.  But note that this is only a concern when both $\Gamma$ is `stronger' than $\mathsf{T}$ (in the sense of interpretability) and the consistency of $\mathsf{T}$ is itself in doubt.  And as I will discuss further in \S 5, these conditions are not satisfied in the specific case when $\mathsf{T}$ can be taken to be $\mathsf{PRA}$ and $\Gamma$ is one of Hilbert's geometric theories.}  Although Hilbert could not have cited such results at the time of his correspondence with Frege, I will also suggest that this does indeed bear retrospectively on how we should understand some aspects of their disagreement with respect to geometry.   But we have also seen how their exchange led them to individually recognize the dimension of difficulty embodied by (D2) as well as its relevance to geometric consistency proofs -- in Frege's case by using it to argue that no non-model-theoretic method of showing consistency could exist and in Hilbert's case by explicitly formulating the decidability of first-order logic as an open mathematical question (see also note \ref{epng}).

Partial solutions to special cases of the \textsl{Entscheidungsproblem} for a fragment of first-order logic were in fact obtained by \citet{Bernays1928c}.   However two other well-known developments separate the context in which the problem was originally framed -- wherein Hilbert clearly appears to have expected a positive solution -- from that in which \citet{Church1936a} and \citet{Turing1936} showed that no general finitary decision procedure for first-order logic can exist.  First, Church and Turing were able to rely on their independently motivated analyses of the notion `finite decision procedure'  -- respectively in the form of their definitions of $\lambda$-definability and computability by a Turing machine.  And second, in order to provide a uniform analysis of decidability for different domains, they were also able to make use of G\"odel's \citeyearpar{Godel1931a} \textsl{arithmetization of syntax} so as to formulate the \textsl{Entscheidungsproblem} as a decision problem about numbers rather than formulas.

These developments are also sufficiently familiar to require little comment.  But to fix the aspects which will be relevant here, recall that G\"odel's method allows us to define a mapping $\ulcorner \cdot \urcorner : \mathfrak{Form} \rightarrow \mathbb{N}$ which effectively assigns a numerical code $\ulcorner \phi \urcorner$ (or \textsl{G\"odel number}) to every sentence $\phi$ in the class $\mathfrak{Form}$ of well-formed first-order formulas.  This in turn allows us to define the set $\textsc{Val} = \{\ulcorner \phi \urcorner : \phi \text{ is a valid first-order formula}\} \subseteq \mathbb{N}$ and also its characteristic function $\mathit{val}(n)$ which returns the value $1$ just in case $n = \ulcorner \phi \urcorner$ if $\phi \in \textsc{Val}$ and $0$ otherwise.   The relevant form of the result proven by Church and Turing can now be formulated as follows:\footnote{See, e.g., \citep{Soare2016} or \citep{Dean2020b} for a review of  the notions from computability theory employed below and \citep{Borger2001} for both a modern presentation of results (\ref{undeca})-(\ref{pi01}) and historical discussion of the decision problem in a broader context.}   
\begin{theorem} \label{undec} $\mathit{val}(x)$ is neither a $\lambda$-definable function nor a Turing computable function. 
\end{theorem}
\noindent Church and Turing also supplied arguments for the adequacy of their models by providing what can reasonably be called conceptual analyses of the notion `effectively computable function'.  Taken together with results demonstrating the extensional coincidence of the Turing computable and $\lambda$-definable functions (as well as those determined by several other models), these are now generally taken to confirm \textsl{Church's Thesis} -- i.e. the claim that the Turing computable (and thus also $\lambda$-definable, etc.) functions coincide with those which are effectively computable by a finite procedure.

But Turing in fact showed a bit more than this which is relevant to understanding (D2).  For in the course of proving Theorem \ref{undec} he demonstrated the following: 
\begin{example}
\label{undeca}
For every Turing machine $T$ and input $n \in \mathbb{N}$, it is possible to effectively construct a first-order formula $\phi_{T,n}$ such that $T(n)\downarrow$ -- i.e. the computation of $T$ on $n$ halts -- if and only if $\proves \phi_{T,n}$ -- i.e. $\phi_{T,n}$ is provable in Hilbert \& Ackermann's axiomatization of first-order logic.  
\end{example}
Suppose we now let $T_0,T_1,\ldots$ be an effective enumeration of Turing machines -- say in increasing order of their G\"odel numbers -- and also define the sets $H = \{\langle i,x\rangle : T_i(x) \downarrow\}$ -- i.e. the \textsl{Halting Problem} consisting of pairs of indices and inputs to Halting Turing machine computations -- and $\textsc{Prov} = \{\ulcorner \phi \urcorner : \ \proves \phi\ \ \& \ \phi \in \mathfrak{Form}\}$ -- i.e. the set of (codes of) provable first-order formulas.   Note also that it follows from the Soundness and Completeness Theorems for first-order logic that $\textsc{Prov} = \textsc{Val}$.     Thus in proving (\ref{undeca}) \citet[pp. 259-263]{Turing1936} also demonstrated the following:
\begin{example} \label{red2} There exists a Turing computable function $f$ such that for all Turing machines $T_i$ and inputs $x$, $\langle i,x \rangle \in H$ if and only if $f(i,x) \in \textsc{Val}$.
\end{example}
What (\ref{red2}) reports is that the Halting Problem for Turing machines is reducible to the problem of deciding first-order validity -- i.e. if there were an effective method for solving the \textsl{Entscheidungsproblem}, then this method could also be employed to solve Turing's Halting Problem.

It is, of course, also a commonplace that $H$ is itself a \textsl{difficult} problem in the sense that Turing showed that it is not solvable by an algorithm (assuming Church's Thesis).    But although I have thus far employed the term `difficult'  informally, Turing's analysis was subsequently employed by \citet{Post1944} to define a general notion of a  \textsl{degree of difficulty} (or \textsl{unsolvability}) which has subsequently been studied in great detail within computability theory.   For instance one way of understanding the difficulty of $H$ is to observe that it is a so-called \textsl{Turing complete} set -- i.e. were we to have access to a (notional) `oracle' for deciding membership in $H$, then we could effectively decide membership in every other computably enumerable set.  But in order to relate the foregoing observations to (D2) more directly, it will be useful to employ not this notion of a so-called \textsl{Turing reduction} but that of a \textsl{many-one reduction} as exemplified by (\ref{red2}).

To this end, first observe that in addition to being Turing complete, the problem $H$ can also be shown to be a so-called \textsl{$\Sigma^0_1$-complete set}.  This means two things:
\begin{example}
\label{sigma01}
\begin{compactenum}[i)]
\item $H$ can be defined by a $\Sigma^0_1$-\textsl{formula} of  $\mathcal{L}^1_{\mathsf{Z}}$ -- i.e. there is an open formula in the language of first-order arithmetic of the form $\psi(x) = \E y \chi(x,y)$ where $\chi(x,y)$ contains only bounded quantifiers such that  $H = \{n : \mathcal{N} \models \psi(\overline{n})\}$.
\item Every set $A \subseteq \mathbb{N}$ which is similarly definable by a $\Sigma^0_1$-formula is many-one reducible to $H$ (notation: $A \leq_m H$) in the sense exemplified by (\ref{red2}) -- i.e. there is a Turing computable function $f(x)$ such that $x \in A$ if and only if $f(x) \in H$ for all $x \in \mathbb{N}$.
\end{compactenum}
\end{example}
Next note that although the validity of $\phi$ is defined in terms of truth with respect to \textsl{all} models, the extensional equivalence of the set $\textsc{Val}$ with $\textsc{Prov}$ also shows that the former set possesses an \textsl{existential} (i.e. $\Sigma^0_1$) definition -- i.e. $\textsc{Val} = \{\ulcorner \phi \urcorner : \mathcal{N} \models \E y \mathrm{Proof}_{\mathsf{FOL}}(\ulcorner \phi \urcorner,y)\}$ where $\mathrm{Proof}_{\mathsf{FOL}}(x,y)$ is a variant of G\"odel's well-known proof predicate restricted to derivability in pure first-order logic.   Putting this together with the $\Sigma^0_1$-completeness of $H$ and the transitivity of the relation $\leq_m$, this shows that $\textsc{Val}$ is also a $\Sigma^0_1$-complete set.

As a final step towards relating these observations to (D2), consider the set  $\textsc{Con} = \{\ulcorner \phi \urcorner \in \mathfrak{Form} : \ \not\proves \neg \phi \}$ -- i.e. the set of sentences which are \textsl{consistent} with the axioms of first-order logic in the sense that their negations are not provable.  It then follows that $\ulcorner \phi \urcorner \in \textsc{Con}$ if and only if $\ulcorner \neg \phi \urcorner \not\in \textsc{Val}$ -- i.e. a first-order formula is consistent just in case its negation is not a valid formula.   It thus follows from this that $\textsc{Con}$ is what is known as a \textsl{$\Pi^0_1$-complete} set -- i.e.
\begin{example}
\label{pi01}
\begin{compactenum}[i)]
\item $\textsc{Con}$ is definable by a $\Pi^0_1$-\textsl{formula} of the form $\theta(x) = \A y \eta(x,y)$ where $\eta(x,y)$ contains only bounded quantifiers -- in fact $\textsc{Con} = \{\ulcorner \phi \urcorner : \mathcal{N} \models \forall y \neg \mathrm{Proof}_{\mathsf{FOL}}(\ulcorner \neg \phi \urcorner,y)\}$.
\item  Every set $B \subseteq \mathbb{N}$ which is similarly definable by a $\Pi^0_1$-formula is many-one reducible to $\textsc{Con}$ -- i.e. $B \leq_m \textsc{Con}$.   
\end{compactenum}
\end{example}

The problem of deciding membership in $\textsc{Con}$ is thus an arithmetized version of the original question of checking the consistency of a finite set of axioms $\Gamma$.   We have seen that Hilbert stressed that the logical form of an assertion of consistency has a universal structure (\ref{derdefn}) which, via arithmetization, is parallel to that of (\ref{pi01}i).   But in addition to this, the $\Pi^0_1$-\textsl{hardness} of $\textsc{Con}$ (as reported by \ref{pi01}ii)  shows that the problem of deciding $n \in A$ for any $\Pi^0_1$-definable set $A$ can be uniformly reduced to that of deciding whether a particular first-order formula $\phi_{n,A}$ is consistent with the axioms of first-order logic.   

The significance of this can be further reinforced by observing that the $\Pi^0_1$-definable sets include $\mathrm{Tr}_{\Pi^0_1} = \{\ulcorner \phi \urcorner : \phi \text{ is a $\Pi^0_1$-formula of $\mathcal{L}^1_{\mathsf{Z}}$ and } \mathcal{N} \models \phi\}$ -- i.e. set of all \textsl{true} $\Pi^0_1$-statements about the natural numbers.  It thus follows that deciding not just the provability from specified axioms but also the \textsl{truth} of any $\Pi^0_1$-sentence about the natural numbers can be uniformly reduced to checking the consistency of a single first-order formula.  Many famous open problems in number theory fall into this class -- e.g. the Goldbach Conjecture and the Riemann Hypothesis.  Given the extensive efforts which have been mounted to prove these statements over the course of many years, it does indeed seem implausible that there could exist a uniform procedure which  allows us to algorithmically determine their truth values.  Thus even when Frege's concerns about (D1) are set aside, his expectations about the intrinsic difficulty of determining consistency are indeed confirmed by the computability-theoretic analysis of the \textsl{Entscheidungsproblem}.  For once the logical form of consistency statements -- as both he and Hilbert understood them -- is taken into account, problem (D2) is indeed \textsl{as difficult as it can be} (per thesis \ref{points}ii).   

\section{On the ease of existence}

The foregoing argument for thesis (\ref{points}ii) makes use of concepts from computability theory to provide an analysis of the difficulty of determining consistency in the sense of (D2).   The notion of difficulty which is at issue here is an epistemological one -- i.e. it is evident that both Frege and Hilbert had amongst their concerns the practical question of determining how we can \textsl{come to know} that certain sets of axioms $\Gamma$ are consistent.   On the other hand thesis (\ref{points}iii) -- i.e. that Hilbert was correct to maintain that demonstrating the existence of a model satisfying $\Gamma$ is as \textsl{easy} as it can be conditional on its consistency -- has both epistemological and ontological dimensions.   For assuming $\Gamma$'s consistency, there is the question of gauging the difficulty of constructing a model of $\mathcal{M} \models \Gamma$ and determining its properties.   But 
there is also the question of what it means for such a structure to exist in the first place.

A useful waypoint in navigating this relationship is Bernays's \citeyearpar{Bernays1950} paper `Mathematische Existenz und Widerspruchsfreit'.  Therein he makes the following remark about its titular concern:
\begin{quote}
{\footnotesize The common acceptance of the explanation of mathematical existence in terms of consistency is no doubt due in considerable part to the circumstance that on the basis of the simple cases one has in mind, one forms an unduly simplistic idea of what consistency (compatibility) of conditions is. One thinks of the compatibility of conditions as something the complex of conditions wears on its sleeve, as it were, such that one need only sort out the content of the conditions clearly in order to see whether they are in agreement or not. In fact, however, the role of the conditions is that they affect each other in functional use and by combination. The result obtained in this way is not contained as a constitutent part of what is given through the conditions. It is probably the erroneous idea of such inherence that gave rise to the view of the tautological character of mathematical propositions. \hfill (p. 98)}
\end{quote}
This passage makes clear that Bernays was aware not only of the conceptual problem posed by (D2) but also how the technical details pertaining to the \textsl{Entscheidungsproblem} bear on the Frege-Hilbert controversy.\footnote{By this point he had in fact offered a state-of-the-art account of the status of the decision problem for a number of logical systems -- inclusive of the results of Church and Turing -- in \citep[Sup. II]{Hilbert1939}.}  But he suggests that although there may be a sense in which the consistency of $\Gamma$ \textsl{does} entail the existence of a satisfying model, he also observes that such a structure need not be `contained as a constituent' or `inherent' within $\Gamma$ itself.  Much as Frege's discussion of consistency proofs anticipates the unsolvability of the \textsl{Entscheidungsproblem}, I will suggest below that Bernays's remarks foreshadow several specific results which grew out of his own work on the Arithmetized Completeness Theorem.  But in order to appreciate this it will be useful to return again to the problem of demonstrating the consistency of mathematical theories which can be shown to possess no finite models.    

As we have seen, one example of such a theory is provided by the conjunction of Hilbert's Incidence and Order axioms for geometry -- i.e. I.1-8 and II.1-4 in \citep{Hilbert1971}.  Theorem 7 of this edition states the following: `Between any two points on a line there exist an infinite number of points.'   Since axiom II.1 requires the distinctness of points $A,B,C$ standing in the betweenness relation, this result can be obtained from Theorem 3 which states `For any two points $A$ and $C$ there always exists at least one point $D$ on the line $AC$ that lies between $A$ and $C$.'  This can in turn be obtained by an argument which uses Axiom II.4 -- a version of Pasch's axiom which is informally glossed as `If a line enters the interior of a triangle, it also leaves it'  -- to construct $D$ as the intersection of $AC$ and a line $EG$ determined relative to another point $F$ not incident on $AC$ whose existence is inferred by axiom I.3 -- i.e. `There exist three points which do not lie on a line'.  

In the first volume of \textsl{Grundlagen der Mathematik} \citeyearpar[pp. 5-6]{Hilbert1934} Hilbert and Bernays illustrate how it is possible to formalize this subset of Hilbert's axioms in a one-sorted first-order language $\mathcal{L}_G$ containing predicates $Gr(x,y,z)$ (for \textsl{geraden}) intended to express that $x,y,z$ are collinear points and $Zw(x,y,z)$ (for \textsl{zwischen}) intended to express that $x$ is between $y$ and $z$.  They then present nine sentences in this language $\{\theta_1,\ldots,\theta_9\} = \Theta$ by which the axioms mentioned above can be deduced and thus from which the formalization of Theorem 3 can also be derived in first-order logic.   It is then easy to see that the domain of any model $\mathcal{M} \models \Theta$ must be infinite.    

The foregoing situation thus presents a characteristic example of case \ref{fininf}ii) -- i.e. a situation in which we can prove mathematically that any model satisfying a given set of axioms must be infinite (presuming one exists at all).   As we have seen, \citet{Hilbert1928} described the problem of proving consistency in this case as `incomparably more difficult' than that in which a finite model can be exhibited.   \citet{Hilbert1934} concretely illustrated this contrast by first observing that a subset of Hilbert's axioms which weakens II.4 is satisfied in a 5-element model (p. 13).  But they also note that it is easy to see that a theory of an irreflexive, transitive, and serial relation can only have infinite models.  And upon observing that that this theory is satisfied in the structure $\langle \mathbb{N}, < \rangle$, they conclude by remarking that this `exhibition does not conclusively settle the issue of consistency' but rather (again) `reduces it to the consistency of number theory' (p. 15).  

But as we are now in a better position to appreciate, this is at least potentially an \textsl{overstatement}.   For is at least possible that the axioms under consideration in this case are sufficiently weak that their consistency can be demonstrated by a `direct' consistency proof which be itself be formalized in a suitably finitary theory of arithmetic.   Bernays returned to this point immediately after his proof of the Arithmetized Completeness Theorem \citeyearpar[p. 253]{Hilbert1939}.   He observes there that the question of the consistency of a finitely axiomatizable theory $\Gamma$ is reduced via arithmetization to the truth of a particular  $\Pi^0_1$-statement $\mathrm{Con}(\Gamma)$ which can be defined as $\A x \neg \mathrm{Proof}_{\mathsf{FOL}}(x,\ulcorner \bigwedge \Gamma \urcorner)$.\footnote{Bernays constructed a $\Pi^0_1$-formula expressing that $\phi$ is irrefutable by showing that there exists a primitive recursive function $\mathfrak{q}(x)$ such that the truth of $\A x (\mathfrak{q}(x) = 0)$ is equivalent to the non-derivability of $\neg \phi$.   But the structure of $\mathfrak{q}(x)$ is unlike the more familiar proof predicate employed by \citet{Godel1931a} (and also here) in that it is based not irrefutability in a Hilbert system but rather what we now call \textsl{Herbrand consistency} -- i.e. the fact that the so-called \textsl{Herbrand expansions} $\phi_n$ of $\phi$ are truth funcitonally satisfiable for all $n$ (see, e.g., \citealp[\S III.3c]{Hajek1998}). Bernays's presentation of Herbrand's Theorem itself \citeyearpar[\S III.3]{Hilbert1939},  suggests he understood this notion as an intermediary between a syntactic definition of irrefutability and a semantic definition of satisfiability (on which see also \citealp{Franks2009}).  Although the details cannot be considered further here, his proof of Theorem \ref{act2}i \citeyearpar[\S 4.2]{Hilbert1939} thus appears to be connected in a deeper way to how Bernays viewed the relationship between consistency and existence. \label{herbrandnote}} But in this case the truth of $\mathrm{Con}(\Gamma)$ implies the existence of an arithmetical model of $\Gamma$ in virtue of Theorem \ref{act}. 

It is easy to see that the truth of $\mathrm{Con}(\Gamma)$ is equivalent to the membership of $\ulcorner \bigwedge \Gamma \urcorner$ in the set $\textsc{Con}$ defined above.   As we have seen that $\textsc{Con}$ is a $\Pi^0_1$-complete set, there can be no means of effectively determining if $\mathrm{Con}(\Gamma)$ is a true arithmetical statement \textsl{in the general case}.   On the other hand, this does not preclude that in certain instances $\mathrm{Con}(\Gamma)$ can be proven in an appropriate theory.   And in this regard we can now record the following result which was promised in \S 4:\footnote{A proof-theoretic demonstration of a special case of this result using Hebrand's Theorem is given by \citet{Beeson2015}.  But as noted by \citet[10.3.4]{Baldwin2018}, Theorem \ref{congeo} can be approached more uniformly from the standpoint of Reverse Mathematics.  A central observation in this regard is that systems like $\mathsf{HP}$ or $\mathsf{HP} + \mathrm{CCP}$ whose models are coordinatized by Pythagorean or Euclidean fields are also satisfied in so-called \textsl{real closed fields} (see \S 7).  The existence of the real closure $\widetilde{\mathbb{Q}}$ of $\mathbb{Q}$ -- i.e. the ordered field of the real algebraic numbers -- can be proven in the subsystem of second-order arithmetic known as $\mathsf{RCA}_0$ consisting of $\mathrm{I}\Sigma_1$ plus comprehension for $\Delta^0_1$-definable sets \citep[II.9.7]{Simpson2009}.   Reasoning in $\mathsf{RCA}_0$ it is then possible to computably construct models $\mathcal{M} \models \Delta$ over $\widetilde{\mathbb{Q}}$.  Via a formalized cut elimination argument, it is also possible to prove the formalization of the Soundness Theorem for first-order logic in $\mathsf{RCA}_0$ \citep[II.8.8]{Simpson2009} from which it thus follows that $\mathsf{RCA}_0 \proves \mathrm{Con}(\Delta)$.  But since $\mathsf{RCA}_0$ is conservative over $\mathsf{PRA}$ for $\Pi^0_2$-statements \citep[IX.1.11]{Simpson2009} -- and as $\mathrm{Con}(\Delta)$ is a $\Pi^0_1$-statement -- $\mathsf{PRA}$ proves the consistency of $\Delta$ as well.  \label{congeonote}} 
\begin{theorem} \label{congeo}  Let $\Delta$ denote any of the sets of Hilbert's geometric axioms considered in \textnormal{(\ref{mcr}i,ii)}.    Then $\mathsf{PRA} \proves \mathrm{Con}(\Delta)$.
\end{theorem}

It is reasonable to assume that Theorem \ref{congeo} provides a demonstration that the systems considered in (\ref{congeo}i,ii) are consistent which is `direct' in the manner Hilbert came to understand this term.   Theorem \ref{act} can thus be invoked to obtain arithmetical models of these theories in exactly the manner which Bernays envisioned.  For instance when applied to the theory $\Theta$ described above, such a model will take the form of a structure $\mathcal{T} = \langle \mathbb{N},G,Z\rangle \models \Theta$ where $\mathbb{N}$ can be understood as the domain of `points' and $G,Z \subseteq \mathbb{N}^3$ are the extensions of the collinearity and betweenness predicates $Gr(x,y,z)$ and $Zw(x,y,z)$ as defined by the $\mathcal{L}^1_{\mathsf{Z}}$-formulas  $\gamma(x,y,z)$ and $\zeta(x,y,z)$ whose existence is yielded by Theorem \ref{act}.\footnote{Hilbert and Bernays do not explicitly construct such a model of geometry themselves.   On the other hand, at the beginning of the second volume of \textsl{Grundlagen der Mathematik} they did anticipate Theorem \ref{congeo} by showing how the consistency of a geometrical theory similar to $\Theta$ can be proven using the First Epsilon Theorem (pp. 38-48).  \citet[pp. 88-94]{Bernays1954} would also go on to describe in detail how such a model could be constructed for the set theoretic system $\mathsf{GB}^-$ describe below.   He would later go on to characterize such constructions as follows: `By making the deductive structure of a formalized theory one's object of study, as suggested by Hilbert, that theory is, as it were, projected into number theory. The resulting number-theoretic structure is, in general, essentially different from the structure intended by the theory; nevertheless, it can serve to recognize the consistency of the theory, from a viewpoint that is more elementary than the assumption of the intended structure.'  \citeyearpar[p. 63]{Bernays1970}. For more on the nature and significance of such `projections', see \citep{Dean2017c} and \citep[\S 3]{Sieg2020}. \label{proj}}

With an example of such a structure finally in hand, many of the original battle lines in the Frege-Hilbert controversy may be revisited.  For despite Frege's specific misgivings, it is still possible to view models like $\mathcal{P}$ or $\mathcal{E}$ described in \S 2 from the perspective of analytic geometry as consisting of points in Euclidean space coordinatized by pairs of real numbers.  But the fact that the `points' of $\mathcal{T}$ are \textsl{natural numbers} is illustrative of how models obtained by Theorem \ref{act} are even further removed from spatial intuition.   This is reinforced by observing that the extension of the predicate $\zeta(x,y,z)$ which interprets the betweenness relation $Zw$ in $\mathcal{T}$ may hold between a given triple of natural numbers -- say $\mathcal{T} \models Zw(\overline{42},\overline{2},\overline{17})$ and thus also $\mathcal{N} \models \zeta(\overline{42},\overline{2},\overline{17})$  -- despite the fact that the numbers themselves are not between one another with respect to the \textsl{less than} relation on $\mathbb{N}$ -- i.e. $\mathcal{N} \not\models \overline{2} < \overline{42} < \overline{17}$.\footnote{This reflects the fact that the specific correlation between natural numbers and the points which one might take to form a synthetic model $\mathcal{S} \models \Theta$ is determined by largely unconstrained decisions about how the arithmetization of syntax underlying Theorem \ref{act} is carried out.   Thus even when $\mathcal{S}$ is countable -- in which case its points may be put into one-to-one correspondence with $\mathbb{N}$ -- $\mathcal{T}$ will be unlike $\mathcal{P}$ or $\mathcal{E}$ in the sense that it cannot (at least without additional artifice) be regarded as arising from the coordinatization of a synthetic structure.}  One can thus envision Frege doubling down on his basic critique that since Hilbert's models reinterpret the geometric primitives in a manner which divorces them from their intended meanings, Hilbert had no right to describe the reinterpreted statements as axioms and thus also no justification for asserting that their deductive consequences express truths of geometry.

But one of Hilbert's best-known responses to Frege still seems apt:
\begin{quote}
{\small 
But it is surely obvious that every theory is only a scaffolding or schema of concepts together with their necessary relations to one another, and that the basic elements can be thought of in any way one likes. If in speaking of my points I think of some system of things, e.g. the system: love, law, chimney sweep $\ldots$ and then assume all my axioms as relations between these things, then my propositions, e.g. Pythagoras' theorem, are also valid for these things. In other words: any theory can always be applied to infinitely many systems of basic elements. \hfill (IV/4, p. 42)
}
\end{quote}

This passage has traditionally been taken to support the view that Hilbert regarded systems of mathematical axioms as \textsl{implicit definitions}.   But by the time of \textsl{Grundlagen der Mathematik} he and Bernays had also come to systematically distinguish between what they called `formal axiomatics' and `contentual axiomatics'. The former applies to geometrical theories like $\Theta$ whose consistency cannot be secured by exhibition and whose axioms thus not only lack `a special epistemic relation to [a] specific subject matter' (p. 2) but also  `transcend the realm of experience and intuitive self-evidence' (p. 3).   In such cases, Hilbert and Bernays do explicitly endorse an implicit interpretation:\footnote{See also \citep[p. 497]{Bernays1967}.}
\begin{quote}
{\small 
In formal axiomatics $\ldots$ the basic relations are not conceived to be contentually determined from the outset; rather, they receive their determination only \textsl{implicitly} through the axioms.   And any considerations within an axiomatic theory may make use only of those aspects of the basic relations that are explicitly formulated in the axioms.  $\P$ Thus, in axiomatic geometry, whenever we use names that correspond to intuitive geometry -- such as ``lie on'' or ``lie between'' -- this is just a concession to custom, and a means of simplifying the connection of the theory with intuitive facts. Actually, however, in formal axiomatics, the basic relations play the role of \textsl{variable} predicates. \hfill \citeyearpar[p. 7]{Hilbert1934}
}
\end{quote}

By contrast, Hilbert and Bernays describe contentual axiomatics as being grounded `in common experience [which] presents its first principles either as self-evident facts or formulates them as extracts from experience-complexes' (p. 2).    And it is precisely this sort of experience which they took to ground elementary number theory about which they remark
\begin{quote}
{\small 
In the field of arithmetic, we are not concerned with problematic issues connected with the question of the specific character of geometrical knowledge; and it is indeed also here, in the disciplines of elementary number theory and algebra, that we find the purest manifestation of the standpoint of direct contentual thought that has evolved without axiomatic assumptions. $\P$ This methodological standpoint is characterized by \textsl{thought experiments} with things that are assumed to be \textsl{concretely present}, such as numbers in number theory $\ldots$ \hfill \citeyearpar[p. 20]{Hilbert1934}
}
\end{quote}
After this they go on to describe in detail the contentual basis of primitive recursion and quantifier-free induction in virtue of operations on numbers represented in stroke notation.   This exposition has traditionally been taken (e.g. by \citealp[\S 1-\S 2]{Tait2005}) to present the canonical justification for the theory $\mathsf{PRA}$ from the finitary standpoint.   

At least by the time of \textsl{Grundlagen der Mathematik},  Hilbert and Bernays may thus be understood as having made a systematic distinction between `contentual axiomatics' -- as exemplified by $\mathsf{PRA}$ -- and `formal axiomatics' -- as exemplified by the theories of geometry, physics, analysis, and set theory whose consistency they proposed to investigate using contentual theories.   They also held that this distinction was to be mirrored in how the theories are to be interpreted.   For on the one hand they regarded contentual theories as coming along with specific interpretations grounded in intuition.   But on the other hand, they thought of formal theories as lacking a specific subject matter but potentially satisfiable by variable systems of objects and relations subject, of course, to their consistency.  

This in turn suggests the possibility of drawing a similar distinction about the nature of mathematical existence in the two cases.  Bernays would go on to explore such a possibility in his later philosophical writings.  For instance in \citeyearpar{Bernays1950} he suggests that it is still possible to regard the existential consequences of formal theories to justify assertions of what he refers to as `relative existence' -- e.g. if the existence of a model of Euclidean geometry is assumed, then we can understand the consequences of theories like $\mathsf{HP}$ such as  `there is a line connecting two given points' to express a true assertion of relative existence with respect to this model (p. 99).\footnote{\citet{Parsons2014} -- who describes `Mathematische Existenz und Widerspruchsfreit' as `one of the most important contributions to mathematical ontology of its time' -- also takes these remarks as evidence that Bernays held a structuralist view about the nature of mathematical objects more generally.  I will return to this possibility in \S 6.}   But when we ask after the `independent existence' of such infinitary structures themselves, he suggests that we are led via a familiar route to consider the theory of the mathematical continuum (i.e. analysis), and from there to the number series via arithmetization.   He finally asks
\begin{quote}
\footnotesize{But where do all these reductions lead? We finally reach the point at which we make reference to a theoretical framework [\textsl{ideellen Rahmen}]. It is a thought-system that involves a kind of methodological attitude; in the final analysis, the mathematical existence posits [\textsl{Existenz-Setzungen}] relate to this thought system $\ldots$ But we notice here again that we cannot simply identify existence with consistency, for consistency applies to the framework as a whole, not to the individual thing being posited as existent.  \hspace*{1ex} \hfill \citeyearpar[p. 100]{Bernays1950}}
\end{quote}

Bernays does not explicitly identify what he takes such a framework to be.   But he does remark that mathematicians work within a basic set of assumptions about which there can be no `\textsl{de facto} doubt'  and whose consistency serves as `the precondition for the validity of the existence posits made within the theoretical framework' (p. 441).  This suggests that what he had in mind is indeed the standpoint of contentual mathematics within which the proof-theoretic investigations undertaken in \textsl{Grundlagen der Mathematik} are carried out.  On the basis of this interpretation, I will now suggest that the continued study of the Completeness Theorem allows us to further refine our understanding of mathematical existence in regard to formal axiomatics.   

To this end it will be useful to formulate the following extensions to Theorem \ref{act}:\footnote{The proof of part i) of this theorem can be extracted from Bernays's original proof of Theorem \ref{act} for $\mathsf{Z} = \mathsf{PA}$ but is stated more explicitly by \citet{Feferman1960}.   Part ii) is stated by \citet[p. 394]{Kleene1952} and builds on prior work of Kreisel, Mostowski, and Hasenjager.   Part iii) is implicit in \citep{Feferman1960} but explicitly stated (e.g.) in \citep{Kaye1991}.  See \citep{Dean2020a} for further reconstruction.}
\begin{theorem}
\label{act2}
Let $\mathsf{T}$ be an arbitrary computably axiomatizable theory, $\mathrm{Con}(\mathsf{T})$ its canonical consistency statement formulated in $\mathcal{L}^1_{\mathsf{Z}}$ and $\mathsf{Z}$ an $\mathcal{L}^1_{\mathsf{Z}}$-theory extending $\mathrm{I}\Sigma_1$.  Then:
\begin{compactenum}[i)]
\item $\mathsf{T}$ is proof-theoretically interpretable in $\mathsf{Z}+\mathrm{Con}(\mathsf{T})$.
\item If $\mathsf{T}$ is consistent, then there exists an arithmetical model $\mathcal{M} \models \mathsf{T}$  interpretable in $\mathcal{N}$ whose atomic diagram is $\Delta^0_2$-definable.
\item If $\mathsf{T}$ is consistent, then there exists an arithmetical model $\mathcal{M} \models \mathsf{T}$ interpretable in $\mathcal{N}$ whose elementary diagram is $\Delta^0_2$-definable.
\end{compactenum}
\end{theorem}
\noindent 

Although this formulation takes advantage of several notions which were not available in the 1930s,  the core parts of Theorem \ref{act2}i were shown by Bernays in \citeyearpar{Hilbert1939} for finitely axiomatizable theories.  What is meant here by $\mathrm{Con}(\mathsf{T})$ is thus simply the extension of Bernays's consistency predicate $\mathrm{Con}(\Gamma)$ to the case of provability from a decidable (but not necessarily finite set) of mathematical axioms $\mathsf{T}$, a special case of which had already been employed by \citet{Godel1931a}.    $\mathrm{I}\Sigma_1$ is a weak subsystem of first-order Peano arithmetic $[\mathsf{PA}]$ which limits the scope of the induction axiom to $\Sigma^0_1$-formulas.   What Theorem \ref{act2} thus reports is that there is an interpretation $(\cdot)^*$ of the language $\mathcal{L}_{\mathsf{T}}$ into the arithmetical language $\mathcal{L}^1_{\mathsf{Z}}$ such that if $\mathsf{T} \proves \phi$, then $\phi^*$ is derivable in $\mathsf{Z}$ together with the assumption of $\mathsf{T}$'s consistency expressed formally as $\mathrm{Con}(\mathsf{T})$.   Bernays demonstrated this by showing how it is possible to formalize G\"odel's original completeness proof in a fragment of $\mathsf{PA}$.   

Although it is easy to show that this fragment can be taken to be $\mathrm{I}\Sigma_2$, it is more involved to show that $\mathrm{I}\Sigma_1$ suffices.   But once this step is undertaken, it is also possible to take advantage of a well-known `proof theoretic reduction' of $\mathrm{I}\Sigma_1$ to $\mathsf{PRA}$ to argue that Theorem \ref{act2}i is accessible from the finitary standpoint.\footnote{On these steps see \citep{Dean2020a}, \citep[\S I.4b]{Hajek1998} and \citep{Feferman1988}.   But also note that Theorem \ref{act}i easily goes through in the case that $\mathsf{Z} = \mathsf{PA}$ which is a subtheory of the system $\mathsf{Z}_{\mu}$ in which Hilbert and Bernays conducted most of their metamathematical investigations in \citeyearpar{Hilbert1939}.  They also state that they regarded `the expression ``finitary'' not as a sharply delimited endpoint, but rather as a designation of a methodological guideline' \citeyearpar[p. 347]{Hilbert1939}.  Hence the technicalities involved with formalization in weak fragments of first-order arithmetic need not be taken as a central concern here.}  But now consider the result of applying Theorem \ref{act2}i to a theory such as $\Theta$.   Since $\mathsf{Z}$ extends $\mathsf{PRA}$ we have that $\mathsf{Z} \proves \mathrm{Con}(\Theta)$ by Theorem \ref{congeo} and thus also
\begin{example}  For all $\mathcal{L}_G$ sentences $\phi$, if $\Theta \proves \phi$, then $\mathsf{Z} \proves \phi^*$.
\end{example}
-- i.e. if a geometrical statement is derivable from $\Theta$, then its image under the sort of interpretation described above is provable in $\mathsf{Z}$.   This concretely illustrates what Hilbert and Bernays appear to have intended in their original description of the method of arithmetization in remarking that the geometrical axioms are translated into \textsl{provable} arithmetical statements \citeyearpar[p. 3]{Hilbert1934}.   And this in turn equips Hilbert with a reply to another of Frege's challenges (cf. letter IV/5 pp. 47-48) in that it demonstrates that geometrical \textsl{theorems} are translated by his methods into  \textsl{arithmetical truths}.\footnote{This assumes, of course, the \textsl{soundness} of arithmetical theories such as $\mathsf{PRA}$.   But this presumably would not have presented an issue for Hilbert and Bernays if they were indeed willing to regard them as contentual.} 

Bernays repeatedly described the Arithmetized Completeness Theorem as a `finite' or `proof theoretic sharpening of G\"odel's Completeness Theorem'.\footnote{E.g. \citeyearpar[p. VI, 191, 243]{Hilbert1939}.}     One way of understanding this is that in the form of Theorem \ref{act2}i what is yielded is not a model satisfying a given consistent set $\Gamma$, but rather arithmetical formulas $\psi_1(\vec{x}),\ldots,\psi_k(\vec{x})$ with arities matching those of the non-logical symbols of $\mathcal{L}_{\mathsf{T}}$.  If $\Gamma$ is indeed consistent, we have seen how the  $\psi_i(\vec{x})$s may be interpreted in the standard model of arithmetic $\mathcal{N}$ to determine extensions for the symbols of $\mathcal{L}_{\mathsf{T}}$ so as to then determine a model $\mathcal{M} \models \Gamma$.   In this way we are equipped with finitary \textsl{descriptions} of the sorts of sets which would exist were $\mathcal{N}$ to exist as a `completed totality' about which we can reason proof theoretically within the contentual theory $\mathsf{Z}$.   But note that Theorem \ref{act2}i itself does not entail the existence of $\mathcal{N}$ or any other infinitary structure.   And this in turn illustrates a means by which the Completeness Theorem may indeed be made accessible from the finitary standpoint.\footnote{See \citet[\S II.10]{Sieg2013} for a related reconstruction of the finitary standpoint in terms of what he calls `accessible domains'.   \citet[p. 53-56]{Quine1970a} also makes a similar point about how the Arithmetized Completeness Theorem can be used to `save on sets' in the course of advocating for a substitutional analysis of logical validity.} 

Some evidence that Bernays viewed matters in this way is provided by the discussion which follows his  presentation of G\"odel's original completeness proof in \citeyearpar{Hilbert1939}.   He observes that G\"odel's result leads to a reduction in the \textsl{Entscheidungsproblem} in the sense that it allows us to identify the valid formulas of first-order logic with the provable ones (p. 189). But on the other hand, he also notes that there is a step in G\"odel's proof which is non-constructive -- or as he puts it `requires the application of the ``tertium non datur'' for integers' (p. 188). In virtue of this he concludes that G\"odel's proof is not finitary as it stands and thus that completeness may only apply to `the set-theoretic treatment' of logic.\footnote{\citet{Bernays1950} formulates a similar point: `Of course -- from the standpoint of classical mathematics and logic -- this [identification of existence with consistency obtained from the Completeness Theorem]  is a valid equivalence. But using this equivalence to interpret existence statements is surely unsatisfactory: If the claim that there is an exception to a universal proposition is considered to be in need of a contentual explanation, since it is an existential statement, then the negation of that universal proposition certainly is not clearer as to its content.'  (p. 98)}   As a potential remedy he then introduces the notion of `effective satisfiability' [\textsl{effektive Erf\"ullbarkeit}] (p. 198) -- i.e. the satisfiability of a formula by sets which are not only arithmetically definable but also \textsl{computable} [\textsl{berechenbare}].  But after posing the question of whether completeness continues to hold in this case -- i.e. is every consistent first-order formula effectively satisfiable? -- he conjectures a negative answer (p. 199).\footnote{It is now well-known that G\"odel's result may not only be formalized using the principle known as \textsl{Weak K\"onig's Lemma} [$\mathsf{WKL}$] (i.e. `every infinite binary tree has an infinite path') but is in fact \textsl{equivalent} to this principle over the weak base theory $\mathsf{RCA}_0$ (see note \ref{congeonote}) in the case of computably axiomatizable theories.   But although a similar argument had been employed by \citet{Skolem1923b} in his attempt to avoid the use of the Axiom of Choice in the proof of L\"owenheim's Theorem, $\mathsf{WKL}$ is still \textsl{prima facie} non-constructive as selecting a path through a given tree may involve a choice between alternatives which cannot be effectively decided.  This was later made precise by \citet{Kleene1950a} to yield the result `there exists an infinite computable binary tree which does not contain a computable path'.  Bernays already anticipated this in observing that G\"odel's proof involves a non-finitary application of the law of the excluded middle to arithmetical formulas with unbounded quantifiers in order to obtain a model $\mathcal{M} \models \phi$ from a tree determined by a sequence of finite ``approximating'' models $\mathcal{M}_n$.  But on the other hand, it is now known that that the adjunction of $\mathsf{WKL}$ to $\mathsf{RCA}_0$ yields an extension which is conservative for $\Pi^0_2$-statements over $\mathsf{PRA}$.   This illustrates a sense in which the proof of the Completeness Theorem can be regarded as non-constructive but still finitary (see also note \ref{rmnote} and \citealp[p. 316]{Sieg2013} for a similar point).}

The passage quoted at the beginning of this section might be taken to suggest that Bernays had reservations about viewing the Completeness Theorem as sufficient demonstration that consistency implies existence.\footnote{See also \citep[p. 97]{Bernays1950}.}   As we have seen, however, Theorem \ref{act2}i yields the existence of a proof-theoretic interpretation of $\Theta$ in $\mathsf{Z}$ which supplies  $\mathcal{L}^1_{\mathsf{Z}}$-formulas $\gamma(x,y,z)$ and $\zeta(x,y,z)$ that define extensions for the predicates $Gr$ and $Zw$ in a model $\mathcal{T} \models \Theta$ which would exist were the existence of $\mathcal{N}$ itself to be granted.  But since  $\gamma(x,y,z)$ and $\zeta(x,y,z)$ are \textsl{arithmetical} formulas, it then becomes natural to ask whether they determine computable relations.   For if a positive answer were obtained then we would also be able to \textsl{decide} if the collinearity and betweenness relations hold between given triples of points in $\mathcal{T}$ by proving or refuting statements of the forms $\gamma(\overline{n},\overline{m},\overline{q})$ or $\zeta(\overline{n},\overline{m},\overline{q})$ within $\mathsf{Z}$.  Thus if $\Theta$ is effectively satisfiable, there would exist an effective means of making such determinations even if we do not wish to accord $\mathcal{T}$ (or $\mathcal{N}$) `independent existence'.\footnote{Although this may at first seem like a purely epistemological consideration, the definitional complexity of $\gamma(x,y,z)$ and $\zeta(x,y,z)$ will also determine whether these formulas are `absolute' with respect to arbitrary arithmetical models of $\Theta$ -- e.g. whether the fact that $\mathcal{T} \models Zw(\overline{n},\overline{m},\overline{q})$ is sufficient to ensure that this relationship continues to hold in another model $\mathcal{T}'$ defined relative to a potentially \textsl{nonstandard} interpretation of $\mathsf{Z}$.  For recall that the computable subsets of $\mathbb{N}^k$ are definable by $\Delta^0_1$-formulas of $\mathcal{L}^1_{\mathsf{Z}}$ which are in turn absolute with respect to nonstandard interpretations.  The question of whether a given theory $\mathsf{T}$ is effectively satisfiable thus also tracks the extent to which the predicates $\psi_i(\vec{x})$ provided by Theorem \ref{act2} rigidly determine their extensions.}

Of course other variants of the Completeness Theorem will often yield a great many models of $\Gamma$ -- e.g. one of every infinite cardinality if $\Gamma$ possesses an infinite model at all, one definable at every level of the arithmetical hierarchy extending $\Delta^0_2$ if $\Gamma$ is an essentially undecidable theory interpreting Robinson's $\mathsf{Q}$ etc.  According to the familiar perspective of \textsl{platonism} -- and absent various domain-specific constraints which might be thought to rule out various `unintended' interpretations (which I will discuss in the next section) -- these models will all exist on an equal ontological footing.   But recall that thesis (\ref{points}iii) concerns not just the ``bare'' existence of a model $\mathcal{M} \models \Gamma$ but rather the \textsl{difficulty} of demonstrating that such a structure exists.\footnote{It was, of course, Bernays  who introduced the term `platonism' (with the lowercase `p') into philosophy of mathematics in his paper `Sur le platonisme dans les math\'ematiques' \citeyearpar{Bernays1935}.  This term is now often used to describe a view which regards the structures studied in infintary mathematics as independently existing abstract objects.  But Bernays himself distinguished between different grades of platonism [\textsl{suppositions \guillemotleft platoniciennes\guillemotright}] the weakest of which accepts precisely the existence of the natural numbers and the applicability of \textsl{tertium non datur} to arithmetical formulas.  But he also considers a much stronger grade which postulates `the existence of a world of ideal objects containing all the objects and relations of mathematics' but concludes that such an `absolute platonism' has been shown to be inconsistent in light of the set-theoretic antinomies.  Much of his paper is thus devoted to characterizing various intermediate grades of platonism -- e.g. in the form of something akin to modern reconstructions of Weyl's predicativism.   As such, another way of approaching thesis (\ref{points}iii) is via the question: \textsl{Which of Bernays's grades of platonism allow for the passage from consistency to existence}?} 

Bernays's introduction of the notion of effective satisfiability can thus be seen as suggesting a further refinement into how we might understand the existence of infinitary structures.   For we can ask not only after the \textsl{cardinality} of a structure satisfying a given set of axioms $\Gamma$, but also the potential \textsl{complexity} of the formulas which determine an arithmetical model.   We can then take advantage of the well-known correspondence between arithmetical and computational complexity of the sort adverted to at the end of \S 4 to make distinctions not just between finite and infinite models but between ones which can and cannot be effectively constructed or decided.  

In order to make sense of this latter distinction, it is useful to recall the difference between the \textsl{atomic} and \textsl{elementary diagrams} of a model $\mathcal{M}$.  Supposing for simplicity that $\mathcal{M}$ has a purely relational language $\mathcal{L} = \{P_1,\ldots,P_k\}$, its atomic diagram $\mathrm{Diag}_{\mathrm{at}}(\mathcal{M})$ of $\mathcal{M}$ corresponds to the set of atomic sentences which it satisfies when new constant symbols $\texttt{a}_0,\texttt{a}_1,\ldots$ are added for each element of its domain $A$ -- i.e. the statements $P_i(\vec{\texttt{a}})$ such that $\mathcal{M} \models P_i(\vec{a})$ and $\vec{\texttt{a}} \in A^{r_i}$.   The elementary diagram $\mathrm{Diag}_{\mathrm{el}}(\mathcal{M})$ of $\mathcal{M}$ is similarly defined as the set of $\mathcal{L}$-formulas $\phi(\vec{\texttt{a}})$ such that $\mathcal{M} \models \phi(\vec{\texttt{a}})$.   Note also that if $\mathcal{M}$ is an arithmetical model, then $A \subseteq \mathbb{N}$ in which case the new constants can be taken to be numerals.   In this case both  $\mathrm{Diag}_{\mathrm{at}}(\mathcal{M})$ and  $\mathrm{Diag}_{\mathrm{el}}(\mathcal{M})$ will be countable and thus, relative to a suitable G\"odel numbering, we can also ask after their definitional (and hence computational) complexity.   Recall finally that a relation $R \subseteq \mathbb{N}^k$ is $\Delta^0_2$-definable if it is both $\Sigma^0_2$- and $\Pi^0_2$-definable -- i.e. there exist $\mathcal{L}^1_{\mathsf{Z}}$-formulas $\phi(\vec{x},y,z)$ and $\psi(\vec{x},y,z)$ such that $R = \{\vec{n} \in \mathbb{N}^k : \mathcal{N} \models \exists y \forall z \phi(\vec{n},y,z)\} =   \{\vec{n} \in \mathbb{N}^k : \mathcal{N} \models \forall y \exists z \psi(\vec{n},y,z)\}$.  

What Theorem \ref{act2}ii reports is  that any consistent computable axiomatizable theory $\mathsf{T}$ not only has an arithmetical model, but one whose atomic diagram is definable in this restricted manner.  In order to appreciate the significance of this, recall that it is a consequence of Kleene's \textsl{Hierarchy Theorem} that the $\Delta^0_2$-definable sets properly \textsl{extend} the $\Sigma^0_1$- and $\Pi^0_1$-definable sets but are properly \textsl{contained in} the $\Sigma^0_2$- and $\Pi^0_2$-definable sets -- i.e. if we equate a formula class with its name then $\Sigma^0_1 \varsubsetneq \Delta^0_2 \subsetneq \Sigma^0_2$ and $\Pi^0_1 \subsetneq \Delta^0_2 \subsetneq \Pi^0_2$.   This illustrates a sense in which deciding membership in a $\Delta^0_2$-definable set is \textsl{more difficult} than deciding membership in either a $\Sigma^0_1$-definable -- i.e. computably enumerable -- set or a $\Pi^0_1$-definable set.   On the other hand, this task is \textsl{easier} than deciding membership in either a $\Sigma^0_2$- or $\Pi^0_2$-definable set in the  sense measured by the arithmetical hierarchy.   

Note also that in the case that $\mathcal{M}$ is an arithmetical model, it is possible to go backwards from  the atomic diagram of $\mathcal{M}$ to effectively construct the extensions of the predicates in $\mathcal{L}$ as the sets $\{\vec{n} : P(\vec{n}) \in  \mathrm{Diag}_{\mathrm{at}}(\mathcal{M})\}$.   In this sense Theorem \ref{act2}ii may be understood as a \textsl{positive} result in the sense that one might \textsl{a priori} expect that certain mathematical theories contain predicates whose extensions must necessarily be highly complex in all of their models (as will be illustrated in \S 6).\footnote{For note that the $\Delta^0_2$ sets also correspond to those hich are \textsl{limit computable} in the sense of Putnam and Shoenfield (see, e.g., \citealp[\S 3.5]{Soare2016}).   In other words $A$ is $\Delta^0_2$ just in case there is a so-called \textsl{trial and error} procedure for deciding $n \in A$ -- i.e. one for which membership can be determined by following a guessing procedure which may make a finite number of initial errors before it eventually converges to the correct answer.}    But also recall that in light of \textsl{Post's Theorem}   the computable sets correspond to the $\Delta^0_1$-definable sets -- i.e. those which are definable by both $\Sigma^0_1$- and $\Pi^0_1$-formulas -- and that this class of sets is a proper subset of both the $\Sigma^0_1$- and the $\Pi^0_1$-definable sets. Thus in posing the question of whether completeness continues to hold if we demand not just satisfiability but \textsl{effective satisfiability}, Bernays was in effect asking whether Theorem \ref{act2}ii can be strengthened to assert that every consistent formula is satisfied in a model all of whose relations are $\Delta^0_1$-definable.

This question was answered in the \textsl{negative} via a series of results inspired by Bernays's conjecture leading to the following concisely statable theorem of \citet{Rabin1958}:
\begin{theorem} \label{rabin} Let $\mathsf{GB}^-$ be the axioms of G\"odel-Bernays set theory without the Axiom of Infinity formulated in the one-sorted language $\mathcal{L}_{S} = \{\in\}$.  Suppose that $\mathcal{M} = \langle \mathbb{N},E \rangle$ is an arithmetical model of $\mathsf{GB}^-$ where $E \subseteq \mathbb{N}^2$ interprets $\in$.  Then $E$ is not a $\Sigma^0_1$ $($i.e. computably enumerable$)$ relation.
\end{theorem}
\noindent  As \citet{Bernays1937} had shown that $\mathsf{GB}^-$ is finitely axiomatizable, this provides a precise negative answer to his question about effective satisfiability -- i.e. the conjunction of the axioms of $\mathsf{GB}^-$ is a single first-order sentence which is not effectively satisfiable.   On the other hand, $\mathsf{GB}^-$ is otherwise a much stronger theory than any of those considered above -- e.g. it is strong enough to prove the existence (as classes) of all arithmetically definable sets.   One might expect on this basis that the interpretation of $\in$ in any model of $\mathsf{GB}^-$ must inevitably be complex relative to the arithmetical hierarchy.   But although Theorem \ref{rabin} shows that membership can never be a computably enumerable relation (and thus also not a computable one) in any model of $\mathsf{GB}^-$, Theorem \ref{act2}ii can also be invoked to show that there is an arithmetical model of this theory in which $\in$ is $\Delta^0_2$-definable.   

This sequence of observations can be extended to show that there exist finite theories $\Gamma$  which are not satisfiable in any model consisting of sets definable by $\Sigma^0_1 \cup \Pi^0_1$-predicates. This result is due to \citet{Putnam1957} who described it as demonstrating that Theorem \ref{act2}ii is the `best possible' strengthening of the Completeness Theorem.    The foregoing results may thus be understood as confirming thesis (\ref{points}iii) in the following precise sense: while every consistent computably axiomatizable theory will have a $\Delta^0_2$-model, there exist single first-order sentences which do not have $\Sigma^0_1 \cup \Pi^0_1$ models.     If the difficulty of demonstrating the existence of $\mathcal{M} \models \Gamma$ is equated with the difficulty of deciding membership in its relations and the latter is equated with the position of the relevant sets in Kleene's arithmetical hierarchy,  then these results suggest that the former task is indeed \textsl{as easy as possible} conditional on having determined the consistency of $\Gamma$ (which we have seen has $\Pi^0_1$-complexity).    

The practical significance of these results can be illustrated relative to the specific theories we have been considering.  For on the one hand, after the 1930s Bernays turned towards set theory in part because he wanted to complete the project of axiomatizing analysis which he had begun in \citeyearpar[Sup. IV]{Hilbert1939}.   It was in this context where he introduced the theory $\mathsf{GB}^-$ in \citeyearpar{Bernays1942} and showed how it could be used to formalize real numbers as Dedekind cuts corresponding to infinite classes of rationals.  Relative to such a formalization, Bernays showed that results such as the Bolzano-Weierstrass Theorem were provable in $\mathsf{GB}^-$ from which it follows that this theory is also sufficiently strong to formalize Hilbert's original geometrical consistency proofs based on analytical models.\footnote{From the contemporary perspective this can be seen directly as $\mathsf{GB}^-$ is mutually interpretable with the theory of second-order arithmetic known $\mathsf{ACA}_0$ which extends $\mathsf{RCA}_0$ with comprehension or first-order arithmetical formulas with second-order parameters.   See note \ref{congeonote} and \citep[\S 2]{Dean2017b}.}  But on the other hand, we can see that $\mathsf{GB}^-$  is paradigmatic of a \textsl{formal} theory in Hilbert and Bernays's sense -- i.e. not only does it fail to have a finite model, but it also fails to be effectively satisfiable.  On the other hand, it follows from the work of Tarski discussed in the next section that while geometric theories like $\mathsf{HP}$ do not possess finite models, they are still effectively satisfiable.

These considerations point towards another means of precisifying thesis (\ref{points}iii) which builds on part iii) of Theorem \ref{act2}.   We have just seen that part ii) is optimal in the sense 
that while it implies that every computably axiomatizable theory has an arithmetical model whose elementary diagram is $\Delta^0_2$-definable, there are specific theories (like $\mathsf{GB}^-$) for which this is indeed the `simplest' model.   Part iii) entails the potentially more surprising fact that every such theory $\mathsf{T}$ has an arithmetical model $\mathcal{M}$ whose entire elementary diagram is $\Delta^0_2$-definable.   But note that $\mathrm{Diag}_{\mathrm{el}}(\mathcal{M})$ acts as a \textsl{satisfaction predicate} for $\mathcal{M}$ in the sense that the membership of $\phi(\vec{a})$ in this set determines if the formula $\phi(\vec{x})$ is satisfied by the vector objects $\vec{a}$ from its domain.  Thus for instance there are models $\mathcal{M}$ satisfying $\mathsf{GB}^-$ -- or even stronger theories which prove the existence of uncountable sets like $\mathsf{GB}$ or $\mathsf{ZF}$ -- in which not just the membership relation but the entire \textsl{truth relation} $\mathcal{M} \models \phi$ for sentences of arbitrary complexity is $\Delta^0_2$-definable.\footnote{Although this may seem like a substantial improvement upon Theorem \ref{act2}ii), the fundamental insight behind part iii) can already be read off from the familiar method of Henkin's completeness proof.  For instance a $\Sigma^0_2$-definition for $\mathcal{M} \models \phi(\vec{c})$ (where in this case $\vec{c}$ is an appropriate vector of `Henkin constants') can first be obtained by observing that a sufficient condition for the constructed model $\mathcal{M}$ to satisfy $\phi(\vec{c})$ is for there to \textsl{exist} a stage in the completion procedure at which this sentence can be added to the maximally complete set $\Gamma^* \supseteq \Gamma$  being constructed such that its adjunction preserves consistency -- a fact which we have seen is equivalent to a $\Pi^0_1$-statement.  But in virtue of the completeness of $\Gamma^*$, this is then equivalent to the $\Pi^0_2$-condition that there is no stage at which $\neg \phi(\vec{c})$ can be consistently adjoined.  Related observations about the Henkin procedure can also be used to provide an alternative  computability-theoretic analysis of the Completeness Theorem obtained by \citet{Jockusch1972}: every consistent computably axiomatized theory $\mathsf{T}$ has a completion $\mathsf{T}^* \supseteq \mathsf{T}$ which is is \textsl{low} in the Turing degrees -- i.e. $\mathrm{deg}(\mathsf{T}^*) =_T \emptyset'$ and thus also an arithmetical model $\mathcal{M}$ whose atomic diagram is of low degree.}

Such observations can in turn be used to illustrate another dimension of the difficulty of demonstrating the existence of a model $\mathcal{M}$ of a given consistent theory $\mathsf{T}$ -- i.e. rather than asking  how hard it is to decide $\mathrm{Diag}_{\mathrm{at}}(\mathcal{M})$, we can also ask what resources are required to prove the existence of the set $\mathrm{Diag}_{\mathrm{el}}(\mathcal{M})$ from which  $\mathcal{M}$ can be recovered as a structure.   The foregoing results illustrate that the existence of such a set can be demonstrated using a restriction of the comprehension scheme of a theory like $\mathsf{GB}$ or $\mathsf{ZF}$ to a small fragment of arithmetical formulas.   But using the more refined techniques of Reverse Mathematics (in the sense of \citealp{Simpson2009}), it can also be shown that a considerably weaker system of second-order arithmetic is already sufficient.  A paradigmatic result in this regard is that the theory $\mathsf{WKL}_0$ consisting of $\mathrm{I}\Sigma_1$ together with comprehension for $\Delta^0_1$-formulas and a form of K\"onig's Infinity Lemma for binary trees is sufficient to prove the existence of $\mathrm{Diag}_{\mathrm{el}}(\mathcal{M})$ for some model $\mathcal{M} \models \mathsf{T}$.   This theory can also be shown to be \textsl{weak} in other respects -- e.g. since $\mathsf{WKL}_0$ is $\Pi^0_2$-conservative over $\mathsf{PRA}$, it does not prove any consistency statements that are not already provable in $\mathsf{PRA}$.  This illustrates in yet another way in which the task of proving the existence of a model $\mathcal{M} \models \mathsf{T}$ is \textsl{easy} -- now measured in terms of the strength of the theory required to prove that $\mathrm{Diag}_{\mathrm{el}}(\mathcal{M})$ exists -- conditional on $\mathsf{T}$'s consistency.\footnote{See \citep[\S II.8, \S IV.3]{Simpson2009} for a presentation of the relevant results about $\mathsf{WKL}_0$ and \citep[\S 3]{Dean2017b} for further discussion of their historical context in regard to G\"odel's original completeness proof.  The framework of computable model theory and Reverse Mathematics also allows for the formulation of a number of more refined questions bearing on whether (or in what sense) consistency implies existence -- e.g. Is the difficulty of proving the existence of a model of $\mathsf{T}$ better understood in terms of constructing $\mathrm{Diag}_{\mathrm{at}}(\mathcal{M})$ or $\mathrm{Diag}_{\mathrm{el}}(\mathcal{M})$? To what extent can these sets differ in complexity?  How should theories be compared in strength in regard to their set (or model) existence consequences?  In order to demonstrate the existence of a set (or model) must we show that it is the \textsl{unique} set (or model) satisfying some expressible property? While germane to this paper, these questions are of a yet more technical nature than can be explored here see.  But, e.g., \citep{Ash2000} and \citep{Eastaugh2019}. \label{rmnote}}

\section{Consistency and existence, redux}

The foregoing arguments for theses (\ref{points}ii) and (\ref{points}iii) involve notions which were not available at the time of Frege's correspondence with Hilbert.   It would thus be an anachronism to ask whether they themselves would have accepted the precise analyses of the difficulty of determining consistency -- i.e. the $\Pi^0_1$-completeness of deciding membership in the set $\textsc{Con}$ of consistent first-order formulas -- or the ease of demonstrating existence -- i.e. the fact that every consistent computably axiomatizable theory has a $\Delta^0_2$-arithmetical model whose existence can be proven in a weak fragment of second-order arithmetic -- which I have proposed.   But it should now be evident that the definitions and results on which these analyses depend grew in a direct and conceptually motivated manner out of developments which their exchange helped to initiate.  It is also clear that they relate to the  problematic of whether (and in what sense) consistency implies existence within the historical frame of the Hilbert program.

It is, however, also reasonable to ask how the proposed analyses bear on these notions as they are employed in contemporary practice, both inside and outside of mathematics.   A paramount concern among contemporary readers is likely to be that even in its original form, the Completeness Theorem appears to allow existence to be demonstrated in too wide a range of cases.   An incipient version of this worry is already illustrated by another of Frege's famous challenges to Hilbert:
\begin{quote}
{\footnotesize Suppose we knew that the propositions
\begin{compactenum}[1)]
\item $\texttt{a}$ is an intelligent being 
\item $\texttt{a}$ is omnipresent
\item $\texttt{a}$ is omnipotent
\end{compactenum}
together with all their consequences did not contradict one another; could we infer from this that there was an omnipotent, omnipresent, intelligent being? This is not evident to me. \hfill (IV/5, p. 47)}
\end{quote}
Suppose we formalize 1) - 3) in the obvious manner by treating $\texttt{a}$ as a constant symbol and  `intelligent', `omnipresent', and `omnipotent' as primitive unary predicates to yield the set $\Omega = \{\mathit{Int}(\texttt{a}),\mathit{Omnipot}(\texttt{a}),\mathit{Omnipres}(\texttt{a})\}$.  In this case it is indeed easy to show that $\Omega$ is consistent, either by a technique such as cut elimination or by directly constructing a model $\mathcal{G}$ in which $\texttt{a}$ denotes an arbitrary object and then interpreting $\mathit{Int},\mathit{Omnipot},\mathit{Omnipres}$ to hold of it.   But since we could, e.g., take $\texttt{a}$ to denote the Empire State Building in such a model, in neither case do we take this to demonstrate that there is an omnipotent, omnipresent, and intelligent being.   

Although this example is trivial as it stands, one can readily foresee Frege extending it by subjecting the senses of these predicates to conceptual analysis so as to obtain another set of statements $\Omega^+$ which records additional relationships between them in an enriched language -- perhaps even in the manner of the ontological argument.  Suppose further that it is still possible to demonstrate the consistency of $\Omega^+$ by some suitably non-model theoretic means.    In this case, the Completeness Theorem could then be invoked to yield a model $\mathcal{G}^+ \models \Omega^+$ in which $\texttt{a}$ would denote some object $b$ in the domain of $\mathcal{G}^+$ which simultaneously satisfies the potentially complex formulas $\mathit{Int}^+,\mathit{Omnipot}^+,\mathit{Omnipres}^+$ which result from such an analysis.    But it requires only passing familiarity with a typical proof of completeness to realize that this will be far from sufficient to ensure the object $b$ will be endowed with all of the properties commonly associated with intelligent, omnipotent, omnipresent beings.  For instance in the case of G\"odel's original proof, $b$ will be a natural number -- say 17 -- and the case of Henkin's \citeyearpar{Henkin1949} proof $b$ will be one of countably many witness constants -- say $\texttt{c}_{17}$.   But  we do not take the existence of either sort of object as testament to the existence of God.   And if $\Omega^+$ were to similarly entail the existence of infinitely many objects satisfying  the descriptions of angels and archangels,  we would not conventionally take the existence of the overall structure $\mathcal{G}^+$ to confirm the existence of such a heavenly host.   

The whimsical quality of this example derives at least in part from its non-mathematical character.   For we can also imagine Hilbert attempting to accommodate Frege's point by observing that while he was presumably thinking of $\Omega$ as a  \textsl{contentual} theory (in virtue of being bound to a subject matter, however ethereal),  the slogan `consistency implies existence' was only intended to apply to \textsl{formal} theories.  But as we have seen, Hilbert and Bernays took theories of this latter sort to be divorced from particular subject matter and instead regarded them as implicit definitions of variable systems of objects and relations.   Several theorist have suggested on this basis that at the time of \citeyearpar{Hilbert1899} Hilbert can be regarded as embracing a form of \textsl{structuralism} -- e.g. rather than regarding the subject matter of geometry as being that of \textsl{space} and attendant spatial notions  of points and lines, Hilbert understood it as being the \textsl{relations} which hold between arbitrary systems of objects which satisfy axioms which we customarily express using the terms `point' and `line'.\footnote{For instance Shapiro \citeyearpar[p. 162]{Shapiro1997}, \citeyearpar[p. 205]{Shapiro2005} cites one of the passages from the correspondence quoted above in support of such an interpretation: `[I]t is surely obvious that every theory is only a scaffolding or schema of concepts together with their necessary relations to one another, and that the basic elements can be thought of in any way one likes $\ldots$ One only needs to apply a reversible one-one transformation and lay it down that the axioms shall be correspondingly the same for the transformed things' (IV/4, p. 42).  Here Hilbert does indeed appear to anticipate the fact that isomorphic structures are elementarily equivalent -- a fact which is further confirmed by his invocation of the `principle of duality' from projective geometry as an illustration (on which see \citealp{Eder2018}).  This appears to be compatible with the contemporary understanding of structures as isomorphism types.  And at least in retrospect, Hilbert's aim in articulating his \textsl{Vollstandigkeitsaxiom} can be understood as that of providing a categorical axiomatization of Euclidean geometry.    But on the other hand,  theories like $\mathsf{HP}$ which play the primary role in \textsl{Grundlagen der Geometrie} admit many non-isomorphic models (even in a given cardinality) -- a fact that Hilbert repeatedly \textsl{exploited} in his independence proofs.  This in turn suggests that what is `implicitly defined' by such theories are more akin to classes of models (some of which are \textsl{non-isomorphic}) than to `structures' in the sense of most contemporary expositions.}  \citet{Parsons2014} has additionally suggested that such a view is more systematically evident in the writings of Bernays from the 1920s onwards -- e.g. if we recall his notion of \textsl{relative existence} from  \citeyearpar{Bernays1950} it is not difficult to see him endorsing the view that `all there is to being' a particular point $A$ or line $\ell$ is to stand in the appropriate relation to other points and lines.

Such an interpretation is helpful for making sense of the prior example.   For if we maintain that `consistency implies existence' only applies when we are willing to adopt a structuralist view about the theory under consideration, then Frege's example can be explained away in virtue of the fact that we are not willing to regard the fact that an object $b$ stands in appropriate relations to a given set of objects comprising a model of $\Omega^+$ as a sufficient condition for its godliness.   On the other hand, it is more difficult to reconcile what Hilbert and Bernays say about implicit definition with a recognized form of structuralism once the status which they wished to assign arithmetical theories such as $\mathsf{PRA}$ is taken into account.   For since they regarded such theories as contentual, they presumably thought (\textsl{pace} \citealp[p. 165]{Shapiro1997}) that there is indeed more to an object being (say) the number 17 then standing in the 18th position in some $\omega$-sequence --  e.g. it must presumably also have the appropriate sort of contentual representation.\footnote{In at least the case of catgeorical theories, it may be possible to acknowledge the current point while also assimilating Hilbert's position to the familiar taxonomy of contemporary structuralisms.   Such an attempt is made by \citet{Doherty2019} who suggests that the general strategy of Hilbert's replies to Frege can be understood as anticipating a form of \textsl{non-eliminative} structuralism similar to that of \citet{Shapiro1997}.  (A related `conceptual' interpretation is also developed by \citealp{Isaacson2011a}.)  A more ambitious reconstruction is proposed by \citet{Sieg1990a,Sieg2014,Sieg2020} who seeks to synthesize the view of mathematical existence Hilbert held around 1900 based on his work in geometry with that which he came to develop in the 1920s based on his work in metamathematics.   Sieg suggests that Hilbert transitioned from an understanding of the axioms of (e.g.) analysis or number theory as characterizations of mathematical structures in a manner similar to \citet{Dedekind1888} during the prior period to one on which such theories need to be understood  `formally' in contradistinction to `contentual axiomatics' as described in \S 5 during the latter period.   But as Sieg also stresses, the proof theoretic strategies of finitary mathematics can be understood as a means of `reducing' or `projecting' infinitary structures to contentual representations given relative to inductively generated `accessible domains' as typified by the natural numbers.   The choice of underlying ontology aside, such an approach has an obvious affinity to various formulations of \textsl{eliminative} structuralism which identify `structures' with sets within the iterative hierarchy rather than taking them to be free-standing objects in their own right (see, e.g. \citealp{Reck2000}).  One could similarly look upon the interpretation developed in \S 5 as a means of understanding Hilbert and Bernays as `number theoretic structuralists' who seek to represent infinitary structures (inclusive of uncountable ones) via arithmetical models.}

This in turn raises another sort of concern about the sort of models whose existence is entailed by the Completeness Theorem.   For as was noted in passing above, G\"odel's original construction will always produce a model $\mathcal{M} \models \Gamma$ with domain $\mathbb{N}$ as long as $\Gamma$ possesses at least one infinite model.\footnote{It is also reasonable at this point to ask whether other  proofs of the Completeness Theorem fare any better in delivering what we would conventionally describe as `intended' models.   A canvas of other familiar alternatives -- e.g. the method of maximally complete sets employed by \citet{Henkin1949} or the interpretation in a Boolean algebra employed by \citet{Rasiowa1953} -- would suggest not.  But since a completeness proof must deliver a satisfying model in the case of \textsl{any} consistent set $\Gamma$, it is difficult to imagine how there could be a \textsl{general} method for constructing intended models which takes into account all of the different subject-specific associations we might have for different choices of $\mathcal{L}_{\Gamma}$.   And thus since all that will be common across different choices of $\Gamma$ will be the metatheoretic definitions of provability and satisfiability themselves, it seems difficult to escape the conclusion that there is indeed an affinity between the sort of mathematical existence which can be promised by a completeness proof and that which is common across most formulations of structuralism -- i.e. a conception of a model (or structure) $\mathcal{M} \models \Gamma$ as consisting of an arbitrary class of objects standing in the relations demanded by the definition of $\Gamma \proves \phi$ together with sufficiently many other relationships to yield a bivalent definition of $\mathcal{M} \models \phi(\vec{a})$.}  Together with Bernays's refinements, I have suggested in \S 5 that it is this feature which would provide a contentual representation of a satisfying structure as a set of  arithmetical predicates which would provide extensions for the non-logical symbols in $\mathcal{L}_{\Gamma}$ were the existence of the standard model of first-order arithmetic $\mathcal{N}$ itself to be granted.  But even if we are willing to take the latter step ourselves, we might still be resistant to accepting the induced arithmetical model as a sufficient condition for the existence of the sort of model we take to be described by $\Gamma$.   In fact it appears that in many branches of mathematics we have a prior practice of endowing the terms in their languages with meanings which are sufficiently precise to prejudice us in favor of some interpretations (and \textsl{against} others) but still which do not  completely determine the structure of a satisfying model in the manner which Hilbert and Bernays demanded of contentual theories.

A familiar case in point is provided by axiomatic set theories such as $\mathsf{ZF}$ (or $\mathsf{GB}$) and their extensions by large cardinal hypotheses.  For it is commonly said that we have a sufficient understanding of the membership relation in order to both justify the axioms of such theories and see that they are satisfied in the so-called cumulative hierarchy $\mathbb{V}$ defined by iterating the ranks $\mathbb{V}(0) = \emptyset, \mathbb{V}(\alpha + 1) = \mathcal{P}(\mathbb{V}(\alpha)), \mathbb{V}(\lambda) = \bigcup_{\alpha < \lambda} \mathbb{V}(\alpha)$.\footnote{See \citep{Kreisel1967} or \citep{Isaacson2011a} for a canonical exposition of this view.}   Such theories are, of course, famously far away from cases in which direct consistency proofs via methods like cut elimination are applicable.  But the informal picture is often taken as an adequate \textsl{bona fides} for the consistency of $\mathsf{ZF}$ even if it leaves open certain questions about the structure of $\mathbb{V}$.  Of course if $\mathrm{Con}(\mathsf{ZF})$ is accepted, Theorem \ref{act2}ii can be applied to yield (e.g.) an arithmetical model $\mathcal{U} = \langle \mathbb{N},E \rangle$ similar to the one considered in Theorem \ref{rabin} in which $E$ provides a $\Delta^0_2$-interpretation of $\in$.   But on the other hand, even if our grasp of $\mathbb{V}$ is partial in certain specific respects -- e.g. in regard to how far the $V(\alpha)$'s extend into the transfinite or about the cardinaility of $\mathcal{P}(\omega)$ -- the familiar picture is also conventionally taken to ensure us that it differs from $\mathcal{U}$ both in virtue of being uncountable and also in assigning $\in$  a highly complex interpretation.\footnote{For note that $\mathsf{ZF}$ proves the existence of every $\mathcal{L}^1_{\mathsf{Z}}$-definable subset of the natural numbers and also that $\mathbb{V}$ is an $\omega$-model and as such agrees with $\mathcal{N}$ on (the set-theoretic translations of) all arithmetical sentences.   It thus follows -- essentially by Tarksi's theorem on the undefinability of truth -- that the restriction of $\in^{\mathbb{V}}$ even to arithmetically definable members of $\mathcal{P}(\omega+1)$ cannot be arithmetically definable (and thus \textsl{ipso facto} not $\Delta^0_2$).} 

Of course Hilbert and Bernays may have welcomed this sort of example as it illustrates how it is possible to assimilate set theoretic reasoning to their finitary standpoint.\footnote{Were they to have done so explicitly, the contemporary tendency would be to brand them as `Skolemites' who perniciously deny the referential determinacy of set-theoretic notions.  But the application of such labels is anachronistic in several respects.  For despite Hilbert's early interest in the Continuum Hypothesis \citeyearpar{Hilbert1900,Hilbert1925},  Bernays original axiomatic development of set theory appears to have been more directly motivated by the goal of providing a smooth formalization of analysis in a manner which continued the treatment in higher-order arithmetic initiated in \citeyearpar[Sup. IV]{Hilbert1939}. Such a classification also fails to take into account the subtleties of Bernays's personal engagement  with the `Skolem paradox' in the course of which he both (e.g.) acknowledges that `the graduations of cardinalities are in a certain way unreal' while also stressing that `classical mathematics does not have reasons to transgress $\ldots$ the axiomatic and logistic analysis of mathematical theories' \citeyearpar[p. 117]{Bernays1957}.  An interesting question which cannot be explored here is whether his later work \citeyearpar{Bernays1976d} on the second-order reformulation of L\'evy's reflection scheme and large cardinals caused him to change his view.}    But if we do not adopt this perspective ourselves, then the foregoing example will be taken by many contemporary theorists to illustrate a case in which we believe we possess \textsl{prima facie} justification to accept the consistency of a theory while at the same time fearing that the sort of `existence' of a satisfying model which it is endowed by the Completeness Theorem is purchased too cheaply.   And of course this problem is compounded once  well-known formal independence and relative consistency proofs are taken into account.  For it is then easy to multiply examples in which the theorem appears to overgenerate models whose existence we might otherwise not wish to countenance.

A familiar but germane example is provided by the application of G\"odel's second incompleteness theorem to sufficiently strong theories of arithmetic such as $\mathsf{PA}$.   For in particular,  G\"odel's result shows that if $\mathsf{PA}$ is consistent, then $\mathrm{Con}(\mathsf{PA})$  is underivable in $\mathsf{PA}$ and thus that $\mathsf{PA} + \neg \mathrm{Con}(\mathsf{PA})$ is also consistent.   So if we are willing to accept that $\mathsf{PA}$ is consistent -- e.g. because we regard it as a contentual theory -- then the Completeness Theorem can be invoked to secure the existence of a model $\mathcal{M} \models \mathsf{PA} + \neg \mathrm{Con}(\mathsf{PA})$.   But of course such a model will not only be \textsl{nonstandard} -- i.e. not isomorphic to $\mathcal{N}$ -- but it will simultaneously affirm $\mathsf{Z}$ and deny its consistency in the relevant internal sense.  Shapiro has suggested that this situation in particular can be understood to present problems for both finitists like Hilbert and Bernays and contemporary structuralists.     

One way of understanding the problem which Shapiro \citeyearpar[pp. 134-135]{Shapiro1997},  \citeyearpar[pp. 71-74]{Shapiro2005} takes the existence of models like $\mathcal{M}$ to pose for finitists is that their existence requires them to acknowledge what might be called a `relativity of arithmetical notions'.   For once $\mathsf{PA}$ (or a weaker theory) has been accepted as contentual and the second incompleteness theorem has been demonstrated to hold for it, there is at least a temptation to regard a theory like $\mathsf{PA} + \neg \mathrm{Con}(\mathsf{PA})$ as an implicit definition of a structure whose existence must then be affirmed in light of the Completeness Theorem.  But as we have seen in \S 3, Hilbert and Bernays also anticipated G\"odel's observations about the entanglement of metamathematical and arithmetical notions as embodied, e.g., by the $\Pi^0_1$-definition of $\mathrm{Con}(\Gamma)$.   And thus since $\mathcal{M}$ is a model of the theory $\mathsf{PA}$, this might cause us to fear that the very notion of consistency is itself \textsl{model-relative} (and thus presumably also \textsl{non-contentual}) in the sense that $\mathrm{Con}(\mathsf{PA})$ can be true in one model but false in another.   But then if we persist in  equating consistency and existence, it might seem that the existence of a model of $\mathsf{PA}$ would itself be thrown into doubt.

But not only did Bernays anticipate this sort of scenario, it is also clear that -- at least in hindsight -- he would have been in a position to reply to Shapiro's concern.   For building on the work of \citet{Feferman1960}, the foregoing situation can be sharpened to show that there is a proof-theoretic interpretation $(\cdot)^*: \mathcal{L}^1_{\mathsf{Z}} \rightarrow \mathcal{L}^1_{\mathsf{Z}}$ of $\mathsf{PA} + \neg \mathrm{Con}(\mathsf{PA})$ in $\mathsf{PA}$.    But then reasoning within $\mathsf{PA}$ itself, it is possible to see that the model defined by $(\cdot)^*$ must be an \textsl{end extension} of what we take to be the standard (or `ground') model of $\mathsf{PA}$ and also that any element $b$ in this model which codes a proof of a contradiction must itself be nonstandard.   It is thus possible to see from within $\mathsf{PA}$ that $b$ is a sort of ideal object -- akin to a point at infinity in geometry -- and thus also that $\mathsf{PA} + \neg \mathrm{Con}(\mathsf{PA})$ is \textsl{not} a contentual theory.  But once this has been agreed, the existence of models like $\mathcal{M}$ need not be seen as threatening either the contentual status of $\mathsf{PA}$ or that of the notion of consistency itself.  For to repeat, Hilbert and Bernays took certain arithmetical theories to be contentual in virtue of their intrinsic relation to a specific subject matter -- e.g. calculations performed on numerals in stroke notation.  But since the argument just rehearsed can be conducted within such a theory itself, it shows in a contentual manner that the distinct theory $\mathsf{PA} + \neg \mathrm{Con}(\mathsf{PA})$ lacks this property.\footnote{With respect to Bernays's antiticipation of Shapiro's problematic see, e.g., \citeyearpar[p. 101]{Bernays1950} and \citeyearpar[p. 118]{Bernays1957}.  The relevant result of Feferman is Theorem 6.5 in \citeyearpar{Feferman1960} which he glosses as follows: `In other words we can construct a ``non-standard model'' of $\mathsf{PA}$ within $\mathsf{PA}$ which $\ldots$ we can verify, axiom by axiom, to be a model of $\mathsf{PA} + \neg \mathrm{Con}(\mathsf{PA})$' (p. 77).   See \citep[pp. 554]{Dean2020a} for the additional argument needed to verify within $\mathsf{PA}$ that this model must also be an end extension.}

In at least this sense, the existence of nonstandard (and possibly other sorts of ``unintended'') models which exist as consequences of the Completeness Theorem does not seem to block the argument which I have suggested in \S 5 allows for the acceptance of the slogan `consistency implies existence' from the finitary standpoint.  But at the same time Shapiro and other structuralists have also suggested for other reasons that consistency is still too liberal a standard for existence and ought to be replaced with a condition often referred to as \textsl{coherence}.   On Shapiro's \citeyearpar[pp. 135-136]{Shapiro1997} formulation this notion is not defined mathematically but rather treated as a primitive whose intended sense is likened to (but not identified with) that of satisfiability, potentially in a proper-class sized structure.\footnote{See also, e.g., Hellman \citeyearpar[p. 19]{Hellman1989}, \citeyearpar[p. 556]{Hellman2005}.}   With this notion in place, \citet[p. 95]{Shapiro1997} includes the following as an \textsl{axiom} of his structure theory: 
\begin{description}
\item[Coherence] If $\phi$ is a coherent formula in a second-order language, then there is a structure that satisfies $\phi$. 
\end{description}

By adopting such principles, contemporary structuralists attempt to obtain a portion of the ontological richness for mathematics which Hilbert boldly pronounced with the claim `if the given axioms do not contradict one another $\ldots$ then the things defined by the axioms exist'.   But on the other hand Shapiro explicitly \textsl{denies} that $\mathsf{PA} + \neg \mathrm{Con}(\mathsf{PA})$ is a coherent theory.   Part of his rationale for doing so is programmatic --  i.e. although $\mathsf{PA} + \neg \mathrm{Con}(\mathsf{PA})$ is consistent, it does not possess a model relative to the `standard' semantics for second-order logic which he wishes to employ for the interpretation of his axiomatization of structure theory.   It should be observed, however, that neither second-order logic nor recent debates about its proper interpretation have played a role in the foregoing.  But Shapiro also asserts his prior refusal to countenance  as structures nonstandard models of arithmetical theories \citeyearpar[p. 133]{Shapiro1997}.   As a result, he effectively cuts off his own development of structuralism from the hope of replacing his Coherence axiom with a mathematical result which might be \textsl{proven} in the spirit of the Completeness Theorem.

In order to see this, it suffices to observe that the set $\textsc{Val}^2$ of second-order validities is $\Pi_2$-complete in the \textsl{L\'evy hierarchy} \citep[Theorem 1]{Vaananen2001}.    $\textsc{Val}^2$  is thus not $\Sigma^m_n$-definable for any $m,n \in \mathbb{N}$ -- i.e.  it is ``above'' the generalized analytical hierarchy and hence of vastly greater complexity then the $\Pi^0_1$ set $\textsc{Val}$ of first-order validities.  It is also easy to see that even the set of first-order arithmetical sentences $\phi \in \mathcal{L}^1_{\mathsf{Z}}$ which are structure-theoretically satisfiable in Shapiro's sense corresponds to precisely those which are true in $\mathcal{N}$.\footnote{This follows since $\mathcal{N}$ is -- up to isomorphism -- the first-order reduct of the unique model of the axioms $\mathsf{PA}^2$ of second-order arithmetic relative to the standard semantics.}  But if we let  $\bigwedge \mathsf{PA}^2$ denote the conjunction of the axioms of second-order Peano arithmetic,  it then follows that the set of sentences $\phi \in \mathcal{L}^1_{\mathsf{Z}}$ such that $\bigwedge \mathsf{PA}^2 \rightarrow \phi$ is structure-theoretically valid is $\Delta^1_1$-complete and thus cannot correspond to the set of sentences provable relative to any computably enumerable deducibility relation for second-order logic.   More generally -- but still \textsl{apropos} of Shapiro's dialectic -- this illustrates the necessity of countenancing nonstandard models if we wish to vindicate a principle akin to `consistency implies existence' via a mathematical theorem rather than by fiat.\footnote{For instance given that $\mathsf{PA} \not\proves \mathrm{Con}(\mathsf{PA})$,  such structures are needed as counter-models to show that $\mathsf{PA} \not\models \mathrm{Con}(\mathsf{PA})$ in much the same way non-Archimedean models of geometry are needed to show that the Archimedean axiom is not a semantic consequence of $\mathsf{HP}$.  See \citep[pp. 27-30]{Isaacson2011a} for a related point. \label{solnote}}

\section{From arithmetic back to geometry}

This is not the place to assess whether the replacement of the mathematical definition of consistency with an extra-mathematical notion such as coherence is a necessary or even a salutary aspect of the development of structuralism in contemporary philosophy of mathematics.   But Shapiro's juxtaposition of G\"odel's completeness and incompleteness theorems highlights another aspect of the relationship between consistency and existence which can be used to provide a final perspective on the results discussed in \S 4 and \S 5.  For recall that another consequence of G\"odel's arithmetization of syntax is that the set $\mathrm{Thm}(\mathsf{T}) = \{\ulcorner \phi \urcorner : \mathsf{T} \proves \phi\}$ of theorems of a computably axiomatizable theory will always be at worst a $\Sigma^0_1$-definable set -- i.e. computably enumerable.  On the other hand, if $\mathsf{T}$ is a consistent theory which includes or interprets Robinson's $\mathsf{Q}$, a generalized form of  G\"odel's first incompleteness theorem entails that $\mathsf{T}$ is \textsl{essentially undecidable} -- i.e. $\mathrm{Thm}(\mathsf{T})$ cannot be a $\Delta^0_1$-definable set nor can $\mathrm{Thm}(\mathsf{T}^+)$ be computable for any consistent extension of $\mathsf{T}^+ \supseteq \mathsf{T}$ in the same language.   And thus in such cases the classification of $\mathrm{Thm}(\mathsf{T})$ as $\Sigma^0_1$ is exact.

Starting in the late 1940 Tarski and his collaborators demonstrated that a surprising variety of theories $\mathsf{T}$ which are not essentially undecidable are still simply undecidable -- i.e. although $\mathsf{T}$ may have some decidable extension, it is not possible to effectively decide whether $\mathsf{T} \proves \phi$ for are arbitrary $\mathcal{L}_{\mathsf{T}}$-sentence $\phi$.\footnote{See \citep{Tarski1953}, \citep{Greenberg2010}, and \citep{Makowsky2019} for more on the history of and references to the results cited in this section.}   This includes the first-order theories of the integers in the language with addition and multiplication $\mathrm{Th}(\mathbb{Z}) =\{\phi : \langle \mathbb{Z},+,\times\rangle \models \phi\}$, the theory of the rationals $\mathrm{Th}(\mathbb{Q}) =\{\phi : \langle \mathbb{Q},+,\times\rangle \models \phi\}$, and theory $\mathrm{Th}(\mathbb{F})$ of ordered fields.  On the other hand, the foregoing examples are all subtheories of a system now called $\mathsf{RCF}$ (for \textsl{real closed fields}) which Tarski also showed to be complete -- i.e. it proves or refutes every statement in its language -- and hence also decidable.\footnote{$\mathsf{RCF}$ can be characterized semantically as the first-order theory of the reals $\mathrm{Th}(\mathbb{R})$ or axiomatically via conjoining to the axioms of an ordered field the statements that every positive number has a square root, that no sum of squares is equal to $-1$, and every polynomial of odd degree in a single variable has a root.}

\citet{Tarski1959a} also presented a novel one-sorted axiomatization of geometry based on treating \textsl{equidistance} of points as a primitive notion together with betweenness.  In this language he presented a set of axioms together with a first-order \textsl{continuity axiom scheme} (which ensures the existence of definable Dedekind cuts) which form a theory $\mathsf{E}_2$ which Tarski called \textsl{elementary geometry}.  Tarski showed that $\mathsf{E}_2$ is sufficient to reconstruct Euclidean geometry in a manner similar to that in which Hilbert had done in \citeyearpar{Hilbert1899}.    But he also showed that $\mathsf{E}_2$ is satisfied in a model $\mathcal{R}$ formed by interpreting its language over the real field in a manner similar to Hilbert's construction of the models $\mathcal{P}$ and $\mathcal{E}$ over the minimal Pythagorean and Euclidean fields.  In this manner it is possible to show that $\mathsf{E}_2$  is mutually interpretable with $\mathsf{RCF}$ and thus also complete and decidable.   

Tarski additionally posed the question whether the theories of Pythagorean and Euclidean planes -- which we have seen in \S 2 are axiomatized by Hilbert's theories $\mathsf{HP}$ and $\mathsf{HP} + \mathrm{CCP}$ -- are decidable.   After much effort this question was answered in the \textsl{negative} by \citet{Ziegler1982} who showed that any finite subtheory of  $\mathrm{Th}(\mathbb{R}) =\{\phi : \langle \mathbb{R},+, \times\rangle \models \phi\}$  -- which can be be shown to include Hilbert's theories relative to a suitable coordinatization -- is undecidable.   Although the work that led to this result was again directly inspired by \textsl{Grundlagen der Geometrie} it was obtained more than 80 years after its original publication.   But as I will now suggest, it also provides a paradigmatic illustration of the complexity of the relationship between consistency and existence.   

Building first on the observation about arithmetical theories from above, suppose that $\mathsf{T}$ is a consistent computably axiomatizable theory interpreting $\mathsf{Q}$ -- a condition which includes both contentual theories like $\mathsf{PRA}$ and infinitary ones like $\mathsf{ZF}$.   As we have seen, in this case the set $\mathrm{Thm}(\mathsf{T})$ of theorems of $\mathsf{T}$ will be $\Sigma^0_1$-definable but not $\Delta^0_1$-definable.  On the other hand, it is a consequence of Theorem \ref{act2}iii that there is an arithmetical model $\mathcal{M} \models \mathsf{T}$ whose elementary diagram is $\mathrm{Diag}_{\mathrm{el}}(\mathcal{M})$ is $\Delta^0_2$-definable.   But since the restriction of $\mathrm{Diag}_{\mathrm{el}}(\mathcal{M})$ to $\mathcal{L}_{\mathsf{T}}$ is a complete theory extending $\mathrm{Thm}(\mathsf{T})$, this cannot be improved to $\Delta^0_1$-definable in virtue of the essential undecidability of $\mathsf{T}$.  

It thus follows that if we were to attempt to operationalize the application of the Completeness Theorem in order to determine if $\mathcal{M} \models \phi$ or $\mathcal{M} \models \neg\phi$ for each $\phi \in \mathcal{L}_{\mathsf{T}}$ we would reach the conclusion that infinitely many of these choices are not only unconstrained by the axioms and must be made non-effectively.    This concretely illustrates a sense in which the model $\mathcal{M}$ -- as \citet{Bernays1950} put it --  `is not contained as a constituent part of what is given through the conditions'.  In particular,   $\mathrm{Diag}_{\mathrm{el}}(\mathcal{M})$ must not only extend the axioms  $\mathsf{T}$ but it is also not even uniquely determined by their deductive consequences $\mathrm{Thm}(\mathsf{T})$.  

Since it is well-known that theories capable of interpreting arithmetic do not decide metamathematical assertions like their consistency statements,  this is in some sense an expected conclusion.   But Ziegler's result shows that a version of this phenomena also arises in the case of geometrical theories.   For although  theories like $\mathsf{HP}$ are not essentially undecidable, a central reason why Hilbert was originally interested in studying them is that they are far from being \textsl{complete} -- e.g. $\mathsf{HP}$ does not decide the Circle-Circle Principle (CCP), $\mathsf{HP} + \mathrm{CCP}$ does not decide the Parallel Postulate (PP), etc.  As we have seen, these independence results  were originally shown by specific model constructions -- e.g. the Beltrami, Klein and Poincar\'e models of non-Euclidean geometry, Hilbert's model of non-Desarguesian geometry, Dehn's models of non-Archimedean or non-Legendrian geometry.   But as Ziegler's result entails that there can be no procedure for determining if an arbitrary $\mathcal{L}_{H}$-sentence $\phi$ is provable or refutable in $\mathsf{HP}$, there is thus also no effective means of determining in which cases such constructions are possible.

This situation thus provides a concrete illustration of the significance of the unsolvability of the \textsl{Entscheidungsproblem} within the domain of applications which first attracted Hilbert towards metamathematics.     But it additionally illustrates that certain axiomatizations of geometry share with arithmetic and set theory the property of admitting a wide class of models which need not be effectively generated.   For on the one hand it is possible to use Tarski's decision algorithm for $\mathsf{RCF}$ to decide whether a given $\phi \in \mathcal{L}_{H}$ is true in $\mathcal{R}$.  This also allows us to see that $\mathsf{HP}$ is satisfied in a countable sub-model of $\mathcal{R}_0$ of $\mathcal{R}$ defined over the so-called \textsl{real algebraic numbers} (i.e. real roots of polynomials with rational coefficients) whose atomic diagram is decidable.   Thus while theories like $\mathsf{HP}$ do not possess \textsl{finite models}, they are still effectively satisfiable in the sense discussed in \S 5.   

On the other hand this is by no means the \textsl{only} means of completing $\mathsf{HP}$.   We have also seen that the sort of axiomatic extension $\mathsf{T}$ which Hilbert considered  possess consistency proofs which can be carried out in $\mathsf{PRA}$ and are thus arguably direct.  On this basis I have suggested in \S 5 that  Hilbert and Bernays were still in a position to acknowledge that there are indeed models for such theories which be obtained via the Completeness Theorem.   But in contradistinction to the situation with the theory $\mathsf{RCF}$ and the model $\mathcal{R}_0$, Ziegler's result demonstrates that we have no general guarantee that such a model $\mathcal{M} \models \mathsf{T}$ will be contained -- qua either atomic or elementary diagrams -- in the axioms themselves nor even that they can always be effectively generated by an auxiliary procedure.  

This illustrates a related phenomena which Bernays anticipated even before the incompleteness theorems were announced:
\begin{quote}
\footnotesize{We must only keep in mind the fact that the formalism of statements and proofs, with which we represent our idea-formation does not coincide with the formalism of that structure we intend in the
concept-formation.  The formalism is sufficient to formulate our ideas about infinite manifolds and to draw from these the logical consequences, but it is, in general, not capable of producing the manifold, as it were, combinatorially from within. \hfill \citeyearpar[p. 262]{Bernays1930}}
\end{quote}
Thus while it may be possible to maintain that consistency implies existence even from the finitary standpoint, the former proof-theoretic fact need not be understood as entailing the latter model-theoretic one in a manner which fixes completely the properties of the satisfying structure.

\bibliographystyle{plainnat}

\begin{thebibliography}{118}
\providecommand{\natexlab}[1]{#1}
\providecommand{\url}[1]{\texttt{#1}}
\expandafter\ifx\csname urlstyle\endcsname\relax
  \providecommand{\doi}[1]{doi: #1}\else
  \providecommand{\doi}{doi: \begingroup \urlstyle{rm}\Url}\fi

\bibitem[Ackermann(1937)]{Ackermann1937}
W.~Ackermann.
\newblock {Die Widerspruchsfreiheit der allgemeinen Mengenlehre}.
\newblock \emph{Mathematische Annalen}, 114\penalty0 (1):\penalty0 305--315,
  1937.

\bibitem[Ash and Knight(2000)]{Ash2000}
C.~Ash and H.~Knight.
\newblock \emph{Computable {S}tructures and the {H}yperarithmetical
  {H}ierarchy}, volume 144.
\newblock North-Holland, 2000.

\bibitem[Baldwin(2018)]{Baldwin2018}
J.~Baldwin.
\newblock \emph{Model Theory and the Philosophy of mathematical practice:
  Formalization without Foundationalism}.
\newblock Cambridge University Press, 2018.

\bibitem[Beeson et~al.(2015)Beeson, Boutry, and Narboux]{Beeson2015}
M.~Beeson, P.~Boutry, and J.~Narboux.
\newblock Herbrand's theorem and non-euclidean geometry.
\newblock \emph{Bulletin of Symbolic Logic}, 21\penalty0 (2):\penalty0
  111--122, 2015.

\bibitem[Bernays(1918)]{Bernays1918a}
P.~Bernays.
\newblock {Beitr\"{a}ge zur axiomatischen Behandlung des Logik-Kalk\"{u}ls
  (Habilitationsschrift, Universit\"{a}t G\"{o}ttingen)}, 1918.
\newblock Reprinted in \citep{Hilbert2013}, pp. 222-274.

\bibitem[Bernays(1930)]{Bernays1930}
P.~Bernays.
\newblock {Die Philosophie der Mathematik und die Hilbertsche Beweistheorie}.
\newblock \emph{{Bl\"{a}tter f\"{u}r deutsche Philosophie}}, 4:\penalty0
  326--367, 1930.
\newblock Reprinted in \cite{Bernays1976}, pp. 17--62 and translated in
  \cite{Mancosu1998}, pp. 234-265.

\bibitem[Bernays(1935)]{Bernays1935}
P.~Bernays.
\newblock Sur le platonisme dans les math{\'e}matiques.
\newblock \emph{L'enseignement mathematique}, 34:\penalty0 52-- 69, 1935.
\newblock {Reprinted in \citep{Bernays1976}, pp. 62-78}.

\bibitem[Bernays(1937)]{Bernays1937}
P.~Bernays.
\newblock {A system of axiomatic set theory: Part I}.
\newblock \emph{Journal of Symbolic Logic}, 2\penalty0 (1):\penalty0 65--77,
  1937.

\bibitem[Bernays(1942)]{Bernays1942}
P.~Bernays.
\newblock {A System of Axiomatic Set Theory: Part III. Infinity and
  Enumerability. Analysis}.
\newblock \emph{Journal of Symbolic Logic}, 7\penalty0 (2):\penalty0 65--89,
  1942.

\bibitem[Bernays(1950)]{Bernays1950}
P.~Bernays.
\newblock {Mathematische Existenz und Widerspruchsfreiheit}.
\newblock In \emph{{Etudes de Philosophie des sciences en hommage \`a Ferdinand
  Gonseth \`a l'occasion de son soixanti\`eme anniversaire}}, pages 11--15.
  Editions du Griffon, Neuchatel, 1950.
\newblock Reprinted in \cite{Bernays1976}, pp. 92-106 and in translation by W.
  Sieg and R. Zach in \textsl{Bernays' Essays on the Philosophy of Mathematics:
  Volume II Reflections on Structural Mathematics and Methodological Frames}
  (Oxford Univress Press, forthcoming).

\bibitem[Bernays(1954)]{Bernays1954}
P.~Bernays.
\newblock {A system of axiomatic set theory: Part VII}.
\newblock \emph{Journal of Symbolic Logic}, 19\penalty0 (2):\penalty0 81--96,
  1954.

\bibitem[Bernays(1957)]{Bernays1957}
P.~Bernays.
\newblock {Betrachtungen zum Paradoxen von Thoralf Skolem}.
\newblock In \emph{Avandlinger utgitt av Det Norske Videnskaps-Akademi i Oslo},
  pages 3--9, Oslo, 1957. Aschehoug.
\newblock Reprinted in \cite{Bernays1976}, pp. 113-118.

\bibitem[Bernays(1967)]{Bernays1967}
P.~Bernays.
\newblock {David Hilbert}.
\newblock In P.~Edwards, editor, \emph{Encyclopedia of Philosophy}, volume~3,
  pages 496--504, New York, 1967. Macmillan Publishing Compan.

\bibitem[Bernays(1976{\natexlab{a}})]{Bernays1976}
P.~Bernays.
\newblock \emph{{Abhandlungen zur Philosophie der Mathematik}}.
\newblock Wiss. Buchgesellschaft, Darmstadt, 1976{\natexlab{a}}.

\bibitem[Bernays(1976{\natexlab{b}})]{Bernays1976d}
P.~Bernays.
\newblock On the problem of schemata of infinity in axiomatic set theory.
\newblock In G.~M{\"u}ller, editor, \emph{{Sets and classes: On the work by
  Paul Bernays}}, volume~84 of \emph{Studies in Logic and the Foundations of
  Mathematics}, pages 121--172. North Holland, 1976{\natexlab{b}}.

\bibitem[Bernays and Sch{\"o}nfinkel(1928)]{Bernays1928c}
P.~Bernays and M.~Sch{\"o}nfinkel.
\newblock Zum entscheidungsproblem der mathematischen logik.
\newblock \emph{Mathematische Annalen}, 99\penalty0 (1):\penalty0 342--372,
  1928.

\bibitem[Bernays(1941)]{Bernays1941d}
Paul Bernays.
\newblock {Sur les questions m\'{e}thodologiques actuelles de la th\'{e}orie
  Hilbertienne de la d\'{e}monstration}.
\newblock In Ferdinand Gonseth, editor, \emph{{Les entretiens de Z\"{u}rich sur
  le fondements et la m\'{e}thode des sciences math\'{e}matiques, 6--9
  D\'{e}cembre 1938}}, pages 144--152, {Z\"{u}rich}, 1941. Leemann \& Co.

\bibitem[Bernays(1970)]{Bernays1970}
Paul Bernays.
\newblock {Die schematische Korrespondenz und die idealisierten Strukturen}.
\newblock \emph{Dialectica}, 24:\penalty0 53--66, 1970.
\newblock Reprinted in \cite{Bernays1976}, pp. 176--188.

\bibitem[Blanchette(1996)]{Blanchette1996}
P.~Blanchette.
\newblock {Frege and Hilbert on consistency}.
\newblock \emph{The Journal of Philosophy}, 93\penalty0 (7):\penalty0 317--336,
  1996.

\bibitem[Blanchette(2012)]{Blanchette2012}
P.~Blanchette.
\newblock \emph{Frege's Conception of Logic}.
\newblock Oxford University Press, 2012.

\bibitem[B{\"o}rger et~al.(2001)B{\"o}rger, Gr{\"a}del, and
  Gurevich]{Borger2001}
E.~B{\"o}rger, E.~Gr{\"a}del, and Y.~Gurevich.
\newblock \emph{The classical decision problem}.
\newblock Springer, Berlin, 2001.

\bibitem[Church(1936)]{Church1936a}
A.~Church.
\newblock {A note on the \textsl{Entscheidungsproblem}}.
\newblock \emph{The Journal of Symbolic Logic}, 1\penalty0 (1):\penalty0
  40--41, 1936.

\bibitem[Curry(1941)]{Curry1941}
H.~Curry.
\newblock {The Paradox of Kleene and Rosser}.
\newblock \emph{Transactions of the AMS}, 50:\penalty0 454--516, 1941.

\bibitem[Curry(1942)]{Curry1942}
H.~Curry.
\newblock The inconsistency of certain formal logic.
\newblock \emph{Journal of Symbolic Logic},  7\penalty0 (3), 115--117, 1942.

\bibitem[Dean(2017)]{Dean2017c}
W.~Dean.
\newblock {Bernays and the Completeness Theorem}.
\newblock \emph{Annals of the Japanese Association for the Philosophy of
  Science}, 25:\penalty0 44--55, 2017.

\bibitem[Dean(2019)]{Dean2019}
W.~Dean.
\newblock Computational complexity theory and the philosophy of mathematics.
\newblock \emph{Philosophia Mathematica}, 27\penalty0 (3):\penalty0 381--439,
  2019{\natexlab{b}}.

\bibitem[Dean(2020{\natexlab{a}})]{Dean2020a}
W.~Dean.
\newblock Incompleteness via paradox and completeness.
\newblock \emph{The Review of Symbolic Logic},  13\penalty0 (2), 541--592, 2020{\natexlab{a}}.

\bibitem[Dean(2020{\natexlab{b}})]{Dean2020b}
W.~Dean.
\newblock Recursive functions.
\newblock In Edward~N. Zalta, editor, \emph{The Stanford Encyclopedia of
  Philosophy}. Metaphysics Research Lab, Stanford University, summer 2020
  edition, 2020{\natexlab{b}}.

\bibitem[Dean and Walsh(2017)]{Dean2017b}
W.~Dean and S.~Walsh.
\newblock The prehistory of the subsystems of second-order arithmetic.
\newblock \emph{The Review of Symbolic Logic}, 10\penalty0 (2):\penalty0
  357--396, 2017.

\bibitem[Dedekind(1888)]{Dedekind1888}
R.~Dedekind.
\newblock \emph{Was sind und was sollen die {Z}ahlen?}
\newblock Vieweg, Braunschweig, 1888.

\bibitem[Detlefsen(1990)]{Detlefsen1990}
M.~Detlefsen.
\newblock {On an alleged refutation of Hilbert's program using G{\"o}del's
  first incompleteness theorem}.
\newblock \emph{Journal of Philosophical Logic}, 19\penalty0 (4):\penalty0
  343--377, 1990.

\bibitem[Detlefsen(2014)]{Detlefsen2014a}
M.~Detlefsen.
\newblock Completeness and the ends of axiomatization.
\newblock In J.~Kennedy, editor, \emph{Interpreting G{\"o}del: Critical
  Essays}, pages 59--77. Cambridge University Press, 2014.

\bibitem[Doherty(2019)]{Doherty2019}
F.~Doherty.
\newblock {Hilbertian Structuralism and the Frege-Hilbert Controversy}.
\newblock \emph{Philosophia Mathematica}, 27\penalty0 (3):\penalty0 335--361,
  2019.

\bibitem[Eastaugh(2019)]{Eastaugh2019}
B.~Eastaugh.
\newblock {Set Existence Principles and Closure Conditions: Unravelling the
  Standard View of Reverse Mathematics}.
\newblock \emph{Philosophia Mathematica}, 27\penalty0 (2):\penalty0 153--176,
  2019.

\bibitem[Eder(2015)]{Eder2015}
G.~Eder.
\newblock {Frege's `On the Foundations of Geometry'} and axiomatic metatheory.
\newblock \emph{Mind}, 125\penalty0 (497):\penalty0 5--40, 09 2015.

\bibitem[Eder and Schiemer(2018)]{Eder2018}
G.~Eder and G.~Schiemer.
\newblock Hilbert, duality, and the geometrical roots of model theory.
\newblock \emph{The Review of Symbolic Logic}, 11\penalty0 (1):\penalty0
  48--86, 2018.

\bibitem[Ewald(1996)]{Ewald1996}
W.~Ewald.
\newblock \emph{From {K}ant to {H}ilbert: {A} {S}ource {B}ook in the
  {F}oundations of {M}athematics.}
\newblock Oxford University Press, New York, 1996.

\bibitem[Ewald and Sieg(2013)]{Hilbert2013}
W.~Ewald and W.~Sieg, editors.
\newblock \emph{{David Hilbert's Lectures on the Foundations of Logic and
  Arithmetic 1917 -- 1933}}.
\newblock Springer, Berlin, 2013.

\bibitem[Feferman(1960)]{Feferman1960}
S.~Feferman.
\newblock Arithmetization of metamathematics in a general setting.
\newblock \emph{Fundamenta Mathematicae}, 49:\penalty0 35--92, 1960.

\bibitem[Feferman(1988)]{Feferman1988}
S.~Feferman.
\newblock Hilbert's {P}rogram relativized: proof-theoretical and foundational
  reductions.
\newblock \emph{The Journal of Symbolic Logic}, 53\penalty0 (2):\penalty0
  364--384, 1988.

\bibitem[Feferman et~al.(1986)]{Godel1986}
S.~Feferman et~al., editors.
\newblock \emph{{Kurt G\"odel Collected Works. Vol. I. Publications
  1929-1936}}.
\newblock Oxford Univeristy Press, Oxford, 1986.

\bibitem[Feferman et~al.(1995)]{Godel1995}
S.~Feferman et~al., editors.
\newblock \emph{{Kurt G\"odel Collected Works. {V}ol. {III}. Unpublished
  Lectures and Essays}}.
\newblock Oxford Univeristy Press, 1995.

\bibitem[Franks(2009)]{Franks2009}
C.~Franks.
\newblock \emph{The Autonomy of Mathematical Knowledge: {H}ilbert's Program
  Revisited}.
\newblock Cambridge University Press, Cambridge, 2009.

\bibitem[Frege()]{Frege1897}
G.~Frege.
\newblock {On Mr. Peano's conceptual notation and my own (1897)}.
\newblock Reprinted in \citep{Frege1984}, pp. 234-248.

\bibitem[Frege(1880)]{Frege1880}
G.~Frege.
\newblock {Booles rechnende Logik und die Begriffsschrift}.
\newblock Translated as ``Boole's Logical Calculus and the Concept-script'' in
  \cite{Frege1979} pp. 9-46, 1880.

\bibitem[Frege(1884)]{Frege1884}
G.~Frege.
\newblock \emph{Die {G}rundlagen der {A}rithmetik}.
\newblock Koebner, Breslau, 1884.

\bibitem[Frege(1903)]{Frege1903}
G.~Frege.
\newblock \emph{{Grundgesetze der Arithmetik}}, volume~2.
\newblock Verlag Hermann Pohle, Jena, 1903.
\newblock Translated as \citep{Frege2013}.

\bibitem[Frege(1906)]{Frege1906}
G.~Frege.
\newblock {\"Uber Schoenflies: Die logischen Paradoxien der Mengenlehre}, 1906.
\newblock Translated as `` On Schoenflies $\ldots$'' in \citep{Frege1979} pp.
  176-183.

\bibitem[Frege(1914)]{Frege1914a}
G.~Frege.
\newblock Logik in der mathematik, 1914.
\newblock Translated as ``Logic in mathematics'' in \citep{Frege1979}, pp.
  219-270.

\bibitem[Frege(2013)]{Frege2013}
G.~Frege.
\newblock \emph{Basic laws of arithmetic (translated and edited by M. Rossberg
  and P. Ebert)}.
\newblock Oxford University Press, Oxford, 2013.

\bibitem[Gabriel and MacGuinness(1980)]{Frege1980a}
G.~Gabriel and B.~MacGuinness, editors.
\newblock \emph{{Gotlob Frege: Philosophical and mathematical correspondence}}.
\newblock Blackwell, Oxford, 1980.

\bibitem[G{\"o}del(1929)]{Godel1929a}
K.~G{\"o}del.
\newblock On the completeness of the calculus of logic.
\newblock Reprinted in \citep{Godel1986}, pp. 64--101, 1929.

\bibitem[G{\"o}del(1930)]{Godel1930c}
K.~G{\"o}del.
\newblock {Vortrag {\"u}ber Vollst{\"a}ndigkeit des Funktionenkalk{\"u}ls}.
\newblock Reprinted in \citep{Godel1995}, pp. 16-29, 1930.

\bibitem[G\"odel(1931)]{Godel1931a}
K.~G\"odel.
\newblock {On} formally undecidable propositions of {{\em {Principia
  Mathematica}\/}} and related systems~{I}.
\newblock Reprinted in \cite{Godel1986}, pp. 144-195, 1931.

\bibitem[Greenberg(2010)]{Greenberg2010}
M.~Greenberg.
\newblock {Old and new results in the foundations of elementary plane Euclidean
  and non-Euclidean geometries}.
\newblock \emph{The American Mathematical Monthly}, 117\penalty0 (3):\penalty0
  198--219, 2010.

\bibitem[H{\'a}jek and Pudl{\'a}k(1998)]{Hajek1998}
P.~H{\'a}jek and P.~Pudl{\'a}k.
\newblock \emph{Metamathematics of First-Order Arithmetic}.
\newblock Springer, Berlin, 1998.
\newblock First edition 1993.

\bibitem[Hallett(1990)]{Hallett1990}
M.~Hallett.
\newblock {Physicalism, reductionism \& Hilbert}.
\newblock In A.~Irvine, editor, \emph{Physicalism in mathematics}, pages
  183--257. Kluwer, Dordrecht, 1990.

\bibitem[Hallett(2008)]{Hallett2008}
M.~Hallett.
\newblock Reflections on the purity of method in hilbert's grundlagen der
  geometrie.
\newblock In P.~Mancosu, editor, \emph{Mancosu, Paolo}. Oxford Univeristy
  Press, Oxford, 2008.

\bibitem[Hallett(2010)]{Hallett2010}
M.~Hallett.
\newblock {Frege and Hilbert}.
\newblock In T.~Ricketts and M.~Potter, editors, \emph{{The Cambridge Companion
  to Frege}}. {C}ambridge {U}niversity {P}ress, Cambridge, 2010.

\bibitem[Hallett and Majer(2004)]{Hilbert2004}
M.~Hallett and U.~Majer, editors.
\newblock \emph{David {H}ilbert's {L}ectures on the {F}oundations of
  {G}eometry, 1891--1902}, volume~1 of \emph{David Hilbert's Foundational
  Lectures}.
\newblock Springer, Berlin, 2004.

\bibitem[Hartshorne(2000)]{Hartshorne2000}
R.~Hartshorne.
\newblock \emph{Geometry: {E}uclid and {B}eyond}.
\newblock Springer, Berlin, 2000.

\bibitem[Hellman(1989)]{Hellman1989}
G.~Hellman.
\newblock \emph{Mathematics without Numbers: Towards a Modal-Structural
  Interpretation: Towards a Modal-Structural Interpretation}.
\newblock Clarendon Press, Oxford, 1989.

\bibitem[Hellman(2005)]{Hellman2005}
G.~Hellman.
\newblock Structuralism.
\newblock In S.~Shapiro, editor, \emph{{The Oxford handbook of philosophy of
  math and logic}}, pages pp. 536--563. Oxford Univeristy Press, Oxford, 2005.

\bibitem[Henkin(1949)]{Henkin1949}
L.~Henkin.
\newblock The completeness of the first-order functional calculus.
\newblock \emph{Journal of Symbolic Logic}, 14\penalty0 (03):\penalty0
  159--166, 1949.

\bibitem[Herbrand(1930)]{Herbrand1930a}
J.~Herbrand.
\newblock \emph{Recherches sur la th{\'e}orie de la d{\'e}monstration}.
\newblock PhD thesis, Universit{\'e} de Paris, 1930.

\bibitem[Hermes et~al.(1979)Hermes, Kambartel, and Kaulbach]{Frege1979}
H.~Hermes, F.~Kambartel, and F.~Kaulbach.
\newblock \emph{Posthumous Writings G. Frege}.
\newblock Blackwell, Oxford, 1979.

\bibitem[Hilbert(1899)]{Hilbert1899}
D.~Hilbert.
\newblock \emph{{Grundlagen der Geometrie}}.
\newblock Leipzig: Teubner, 1899.
\newblock Reprinted in \citep{Hilbert2004}, pp. 407-529.

\bibitem[Hilbert(1900{\natexlab{a}})]{Hilbert1900}
D.~Hilbert.
\newblock Mathematische {P}robleme, 1900{\natexlab{a}}.
\newblock Transalated as ``Mathematical Problems'' in \cite{Ewald1996}, pp.
  1096-1105.

\bibitem[Hilbert(1900{\natexlab{b}})]{Hilbert1900a}
D.~Hilbert.
\newblock {\"Uber den Zahlbegriff}, 1900{\natexlab{b}}.
\newblock Translated as ``On the Concept of Number'' in \citep{Ewald1996}, pp.
  1089-1095.

\bibitem[Hilbert(1904)]{Hilbert1904}
D.~Hilbert.
\newblock {{\"U}ber die Grundlagen der Logik und der Arithmetik}, 1904.
\newblock Translated as ``On the foundations of logic and and arithmetic'' in
  \cite{Heijenoort1967}, 129-139.

\bibitem[Hilbert(1923)]{Hilbert1923}
D.~Hilbert.
\newblock Die logischen {G}rundlagen der {M}athematik, 1923.
\newblock Translated as ``The logical foundations of mathematics'' in
  \cite{Ewald1996}, pp. 1134--1148.

\bibitem[Hilbert(1925)]{Hilbert1925}
D.~Hilbert.
\newblock {\"{U}ber der Unendliche}, 1925.
\newblock Translated as ``On the infinite'' in \cite{Heijenoort1967}, pp.
  367-292.

\bibitem[Hilbert(1927)]{Hilbert1927}
D.~Hilbert.
\newblock {Die Grundlagen der Mathematik}, 1927.
\newblock Translated as ``The foundations of mathematics'' in
  \citep{Heijenoort1967}, p. 464-479.

\bibitem[Hilbert(1971)]{Hilbert1971}
D.~Hilbert.
\newblock \emph{Foundations of geometry}.
\newblock Open Court, La Salle, Illinois, 1971.
\newblock Translated by Leo Unger from the Tenth German Edition.

\bibitem[Hilbert and Ackermann(1928)]{Hilbert1928}
D.~Hilbert and W.~Ackermann.
\newblock \emph{Grundz{\"u}ge der theoretischen Logik}.
\newblock Springer, first edition, 1928.
\newblock Reprinted in \cite{Hilbert2013} pp. 806-916.

\bibitem[Hilbert and Bernays(1934)]{Hilbert1934}
D.~Hilbert and P.~Bernays.
\newblock \emph{{Grundlagen der Mathematik}}, volume~I.
\newblock Springer, Berlin, 1934.
\newblock Second edition 1968. Chapters 1 and 2 translated as
  \citep{HIlbert2011}.

\bibitem[Hilbert and Bernays(1939)]{Hilbert1939}
D.~Hilbert and P.~Bernays.
\newblock \emph{{Grundlagen der Mathematik}}, volume~II.
\newblock Springer, Berlin, 1939.
\newblock Second edition 1970.

\bibitem[Hilbert and Bernays(2011)]{HIlbert2011}
D.~Hilbert and P.~Bernays.
\newblock \emph{{Grundlagen der Mathematik I. Foundations of mathematics I.
  Part A}}.
\newblock College Publications, 2011.
\newblock Translated by C. Wirth et al. from the second German edition of 1968.

\bibitem[Isaacson(2011)]{Isaacson2011a}
D.~Isaacson.
\newblock The {R}eality of {M}athematics and the {C}ase of {S}et {T}heory.
\newblock In Zsolt Novak and Andras Simonyi, editors, \emph{Truth, Reference,
  and Realism}, pages 1--75. Central European University Press, Budapest, 2011.

\bibitem[Jockusch and Soare(1972)]{Jockusch1972}
C.~Jockusch and R.~Soare.
\newblock {$\Pi^0_1$} classes and degrees of theories.
\newblock \emph{Transactions of the AMS}, 173:\penalty0 33--56, 1972.

\bibitem[Kaye(1991)]{Kaye1991}
R.~Kaye.
\newblock \emph{Models of {P}eano {A}rithmetic}, volume~15 of \emph{Oxford
  Logic Guides}.
\newblock Oxford University Press, Oxford, 1991.

\bibitem[Kennedy(2011)]{Kennedy2011a}
J.~Kennedy.
\newblock G\"odel's thesis: An appreciation.
\newblock In Matthias Baaz, editor, \emph{Kurt G{\"o}del and the Foundations of
  Mathematics}, pages 95--109. Cambridge University Press, 2011.

\bibitem[Kleene(1950)]{Kleene1950a}
S.~Kleene.
\newblock A symmetric form of g{\"o}del's theorem.
\newblock \emph{Indagationes mathematicae}, 12:\penalty0 244--24, 1950.

\bibitem[Kleene(1952)]{Kleene1952}
S.~Kleene.
\newblock \emph{Introduction to {M}etamathematics}.
\newblock North-Holland, Amsterdam, 1952.

\bibitem[Kreisel(1967)]{Kreisel1967}
G.~Kreisel.
\newblock Informal {R}igour and {C}ompleteness {P}roofs.
\newblock In I.~Lakatos, editor, \emph{{Problems in the philosophy of
  mathematics}}, pages 138--186. North-Holland, Amsterdam, 1967.

\bibitem[Kreisel(1976)]{Kreisel1976}
G.~Kreisel.
\newblock {What have we learnt from Hilbert's second problem?}
\newblock \emph{Browder [8: 1, pp. 93-130]}, 67:\penalty0 8--12, 1976.

\bibitem[Makowsky(2019)]{Makowsky2019}
J.~Makowsky.
\newblock Can one design a geometry engine?
\newblock \emph{Annals of Mathematics and Artificial Intelligence}, 85\penalty0
  (2-4):\penalty0 259--291, 2019.

\bibitem[Mancosu(1998)]{Mancosu1998}
P.~Mancosu, editor.
\newblock \emph{From {B}rouwer to {H}ilbert: The Debate on the Foundations of
  Mathematics in the 1920s}.
\newblock Oxford University Press, Oxford, 1998.

\bibitem[Mancosu(2003)]{Mancosu2003}
P.~Mancosu.
\newblock The {R}ussellian influence on {H}ilbert and his school.
\newblock \emph{Synthese}, 137:\penalty0 59--101, 2003.

\bibitem[McGuinness(1984)]{Frege1984}
B.~McGuinness, editor.
\newblock \emph{Collected papers on mathematics, logic, and philosophy, ed. B.
  McGuinness}.
\newblock Blackwell, Oxford, 1984.

\bibitem[Parsons(2014)]{Parsons2014}
C.~Parsons.
\newblock {Paul Bernays' later philosophy of mathematics}.
\newblock In \emph{Philosophy of mathematics in the twentieth century :
  selected essays}, pages 67--92. Harvard University Press, Cambridge, MA,
  2014.

\bibitem[Post(1944)]{Post1944}
E.~Post.
\newblock Recursively enummerable sets of postive integers and their decision
  problems.
\newblock \emph{Bulletin of the American Mathematical Society}, 50:\penalty0
  284--316, 1944.

\bibitem[Putnam(1957)]{Putnam1957}
H.~Putnam.
\newblock {Arithmetic models for consistent formulae of quantification theory}.
\newblock \emph{The Journal of Symbolic Logic}, 22:\penalty0 110--111, 1957.

\bibitem[Quine(1970)]{Quine1970a}
W.~Quine.
\newblock \emph{Philosophy of Logic}.
\newblock Harvard University Press, 1970.

\bibitem[Rabin(1958)]{Rabin1958}
M.~Rabin.
\newblock On recursively enumerable and arithmetic models of set theory.
\newblock \emph{Journal of Symbolic logic}, 23\penalty0 (4):\penalty0 408--416,
  1958.

\bibitem[Rasiowa and Sikorski(1953)]{Rasiowa1953}
H.~Rasiowa and R.~Sikorski.
\newblock Algebraic treatment of the notion of satisfiability.
\newblock \emph{Fundamenta Mathematicae}, 40:\penalty0 62--95, 1953.

\bibitem[Rathjen and Sieg(2018)]{Rathjen2018}
M.~Rathjen and W.~Sieg.
\newblock Proof theory.
\newblock In Edward~N. Zalta, editor, \emph{The Stanford Encyclopedia of
  Philosophy}. Metaphysics Research Lab, Stanford University, fall 2018
  edition, 2018.

\bibitem[Reck and Price(2000)]{Reck2000}
E.~H Reck and M.~Price.
\newblock Structures and structuralism in contemporary philosophy of
  mathematics.
\newblock \emph{Synthese}, 125\penalty0 (3):\penalty0 341--383, 2000.

\bibitem[Resnik(1974)]{Resnik1974a}
M.~Resnik.
\newblock {The Frege-Hilbert controversy}.
\newblock \emph{Philosophy and Phenomenological Research}, 34\penalty0
  (3):\penalty0 386--403, 1974.

\bibitem[Rosser(1942)]{Rosser1942}
B.~Rosser.
\newblock {The Burali-Forti Paradox}.
\newblock \emph{The Journal of Symbolic Logic}, 7\penalty0 (1):\penalty0 1--17,
  1942.

\bibitem[Shapiro(1997)]{Shapiro1997}
S.~Shapiro.
\newblock \emph{Philosophy of mathematics: structure and ontology}.
\newblock Oxford University Press, Oxford, 1997.

\bibitem[Shapiro(2005)]{Shapiro2005}
S.~Shapiro.
\newblock {Categories, structures, and the Frege-Hilbert controversy: The
  status of meta-mathematics}.
\newblock \emph{Philosophia Mathematica}, 13\penalty0 (1):\penalty0 61--77,
  2005.

\bibitem[Sieg(1990)]{Sieg1990a}
W.~Sieg.
\newblock Relative consistency and accessible domains.
\newblock \emph{Synthese}, 84:\penalty0 259--297, 1990.

\bibitem[Sieg(2002)]{Sieg2002a}
W.~Sieg.
\newblock {Beyond Hilbert's Reach?}
\newblock In David Malament, editor, \emph{Reading natural philosophy: essays
  in the history and philosophy of science and mathematics}, pages 363--405.
  Open Court, 2002.

\bibitem[Sieg(2013)]{Sieg2013}
W.~Sieg.
\newblock \emph{Hilbert's programs and beyond}.
\newblock Oxford University Press, Berlin, 2013.

\bibitem[Sieg(2014)]{Sieg2014}
W.~Sieg.
\newblock {The ways of Hilbert's axiomatics: Structural and formal}.
\newblock \emph{Perspectives on Science}, 22\penalty0 (1):\penalty0 133--157,
  2014.

\bibitem[Sieg(2020)]{Sieg2020}
W.~Sieg.
\newblock {Methodological frames: Paul Bernays, mathematical structuralism, and
  proof theory}.
\newblock In E.~Reck and G.~Schiemer, editors, \emph{The Prehistory of
  Mathematical Structuralism}, pages 352--382. Oxford University Press, 2020.

\bibitem[Simpson(2009)]{Simpson2009}
S.~Simpson.
\newblock \emph{Subsystems of second order arithmetic}.
\newblock Cambridge University Press, Cambridge, second edition, 2009.

\bibitem[Skolem(1923)]{Skolem1923b}
T.~Skolem.
\newblock {Einige Bemerkungen zur axiomatischen Begriindung der Mengenlehre}.
\newblock In \emph{{Matematikerkongressen i Helsingfors 4~7 Juli 1922, Den
  femte skandinaviska matematikerkongressen, Redogorelse}}, Helsinki, 1923.
  Akademiska Bokhandeln.
\newblock English translation in \citep{Heijenoort1967}, pp. 290-301.

\bibitem[Soare(2016)]{Soare2016}
R.~Soare.
\newblock \emph{Turing Computability: Theory and Applications}.
\newblock Springer, 2016.

\bibitem[Tait(2005)]{Tait2005}
W.~Tait.
\newblock \emph{The {P}rovenance of {P}ure {R}eason}.
\newblock Oxford University Press, Oxford, 2005.

\bibitem[Tarski(1959)]{Tarski1959a}
A.~Tarski.
\newblock What is {E}lementary {G}eometry?
\newblock In \emph{The Axiomatic Method}, pages 16--29. North-Holland,
  Amsterdam, 1959.

\bibitem[Tarski et~al.(1953)Tarski, Mostowski, and Robinson]{Tarski1953}
A.~Tarski, A.~Mostowski, and R.~Robinson.
\newblock \emph{Undecidable {T}heories}.
\newblock North-Holland, Amsterdam, 1953.

\bibitem[Turing(1936)]{Turing1936}
A.~Turing.
\newblock {On computable numbers, with an application to the
  \textsl{Entscheidungsproblem}}.
\newblock \emph{Proceedings of the London mathematical society}, 42\penalty0
  (2):\penalty0 230--265, 1936.

\bibitem[V{\"a}{\"a}n{\"a}nen(2001)]{Vaananen2001}
J.~V{\"a}{\"a}n{\"a}nen.
\newblock Second-{O}rder {L}ogic and {F}oundations of {M}athematics.
\newblock \emph{The Bulletin of Symbolic Logic}, 7\penalty0 (4):\penalty0
  504--520, 2001.

\bibitem[van Heijenoort(1967)]{Heijenoort1967}
J.~van Heijenoort, editor.
\newblock \emph{From {F}rege to {G}\"odel : {A} {S}ource {B}ook in
  {M}athematical {L}ogic, 1879-1931}.
\newblock Harvard University Press, Cambridge, MA, 1967.

\bibitem[Zach(1999)]{Zach1999}
R.~Zach.
\newblock {Completeness before Post: Bernays, Hilbert, and the development of
  propositional logic}.
\newblock \emph{Bulletin of Symbolic Logic}, 5\penalty0 (03):\penalty0
  331--366, 1999.

\bibitem[Ziegler(1982)]{Ziegler1982}
M.~Ziegler.
\newblock Einige unentscheidbare {K}\"orpertheorien.
\newblock \emph{L'Enseignement Math\'ematique}, 28\penalty0 (3-4):\penalty0
  269--280, 1982.

\end{thebibliography}

\renewcommand{\bibpreamble}{Note: In cases where an English translation is available, page references are to the  translation indicated below.}

\end{document}